\documentclass[a4paper,10pt,leqno]{amsart}

\usepackage{defs}
\usepackage{amssymb}
\usepackage{graphicx}
\usepackage[dvips]{hyperref}
\usepackage{enumerate}
\usepackage{graphpap}
\usepackage[all]{xy}
\usepackage{amscd}

\numberwithin{equation}{section}

\newcommand\abrir{\makeatletter\renewcommand{\p@enumii}{}\makeatother}

\theoremstyle{plain}
\newtheorem{teorema}{Theorem}[section]
\newtheorem{corolario}[teorema]{Corollary}
\newtheorem{proposicion}[teorema]{Proposition}
\newtheorem{lema}[teorema]{Lemma}

\theoremstyle{definition}
\newtheorem{nota}[teorema]{Note}
\newtheorem{hipotesis}[teorema]{Hypothesis}
\newtheorem{definicion}[teorema]{Definition}
\newtheorem{notacion}[teorema]{Notation}
\newtheorem{ejemplo}[teorema]{Example}
\newtheorem*{gracias}{Acknoledgement}
\theoremstyle{remark}
\newtheorem{observacion}[teorema]{Remark}

\font\ddpp=msbm10

\def\C{\mathbb{C}}

\def\A{{\hbox{\ddpp A}}}

\def\R{{\hbox{\ddpp R}}}
\def\B{{\hbox{\ddpp B}}}
\def\D{{\hbox{\ddpp D}}}
\def\Z{{\hbox{\ddpp Z}}}
\def\Q{{\hbox{\ddpp Q}}}
\def\I{{\hbox{\ddpp I}}}

\def\N{{\hbox{\ddpp N}}}
\DeclareMathOperator\supp{Supp}
\def\Int{\mbox{Int}}
\def\intt{\mbox{Int}}

\DeclareMathOperator\dis{discrim}
\DeclareMathOperator\res{result}

\title[Discriminant method $\ldots$]{Discriminant method for
the homological monodromy of tame polynomials.}
\author[M. Escario]{Mario Escario Gil}
\address{Escuela de Ingenier\'ia Inform\'atica. Universidad San Jorge.
Parque Tecnol\'ogico Walqa, Edificio 4.
Crta. Zaragoza N-330a, km 566.
22197-Cuarte (Huesca), SPAIN.}
\email{mescario@usj.es}
\thanks{Partially
supported by MTM2004-08080-C02-02}
\date{\today}

\begin{document}

\maketitle

\begin{abstract}
We construct an effective algorithmic method to compute the homological 
monodromy of a complex polynomial which is tame. As an application we show 
the existence of conjugated polynomials in a number field which are not 
topologically equivalent.
\end{abstract}

\section{Introduction}\label{s:intro}
Let $f:\C^n \rightarrow \C$ be a polynomial function. It induces a  locally trivial fibration 
over  $\C\setminus B_f$ with $B_f$ the minimal finite bifurcation set; this is a consequence of Thom's 
work~\cite{th:69} and a proof can be found in~\cite[Appendix A1]{Ph} or \cite[Corollary~5.1]{V}.
Let $T$ be a big disk so that $B_f \subset \Int(T)$ and $\star \in \partial T$. A classical construction 
gives us a geometric basis $(\gamma_i)$ of $\pi_1(T \setminus B_f;\star)$ (see Definition~\ref{def:basegeom})
and  this basis provides a direct sum decomposition of $\tilde{H}_q(f^{-1}(\star))$ (reduced homology over $\Z$)
which depends essentially on the choice of $(\gamma_i)$ (see~\cite{Bro1, DN, DN99,
NN00, NN00:bis, ST95} for varius degrees of generality). With this sum decomposition and if $f$ has only isolated singularities, 
local monodromy $(h_{\gamma_i})_*$ has a block decomposition. There are two kinds of blocks: \emph{local
blocks} which only depend on the local Milnor fibers   and \emph{global blocks}   which depend on the embedings
of the local Milnor fibers into the fixed regular fiber $f^{-1}(\star)$. There are several papers dealing with
the local blocks and how to compute them, for example, Brieskorn singularities by A.~Hefez and F.~Lazzari~\cite{HL}, 
certain singularities and unimodal singularities by A.~M.~Gabri\'elov~\cite{Gab2,Gab74:2} and general methods
like the one which uses real morsifications  due to N.~A'Campo~\cite{AC1,AC2} and S.~M.~Gusein-Zade~\cite{G2,G1} 
and the one which uses an inductive argument due to A.~M.~Gabri\'elov~\cite{Gab}. This is not the situation for the
global blocks. There are some relations  
between local and global blocks (see W.~Neumann and P.~Norbury~\cite[Theorem~3.3]{NN00}) which can give
useful constraints but, nevertheless,  usually these data are 
computed depending  on  the particular polynomial $f$ 
and  a practical complete algorithmic method does not exist in the literature. 

Specially interesting is the case of conjugated polynomials in a number field, i.e., a polynomial
with coefficients in a number field and its image by a Galois isomorphism of the field. 
Due to the Galois isomorphism
both have the same algebraic properties (degree, number of components, global Milnor number,
Alexander polynomials, types and position of singularities, ...).
In order to topologically distinguish these polynomials it seems a good idea to compute 
the global blocks, since they reflect how the Milnor fibers sit in the fixed regular fiber 
and this need not be invariant under Galois isomorphims. 

\begin{figure}[ht]
%%\begin{center}
\vspace*{-1cm}
\hspace*{-3cm}
\includegraphics[scale=.7]{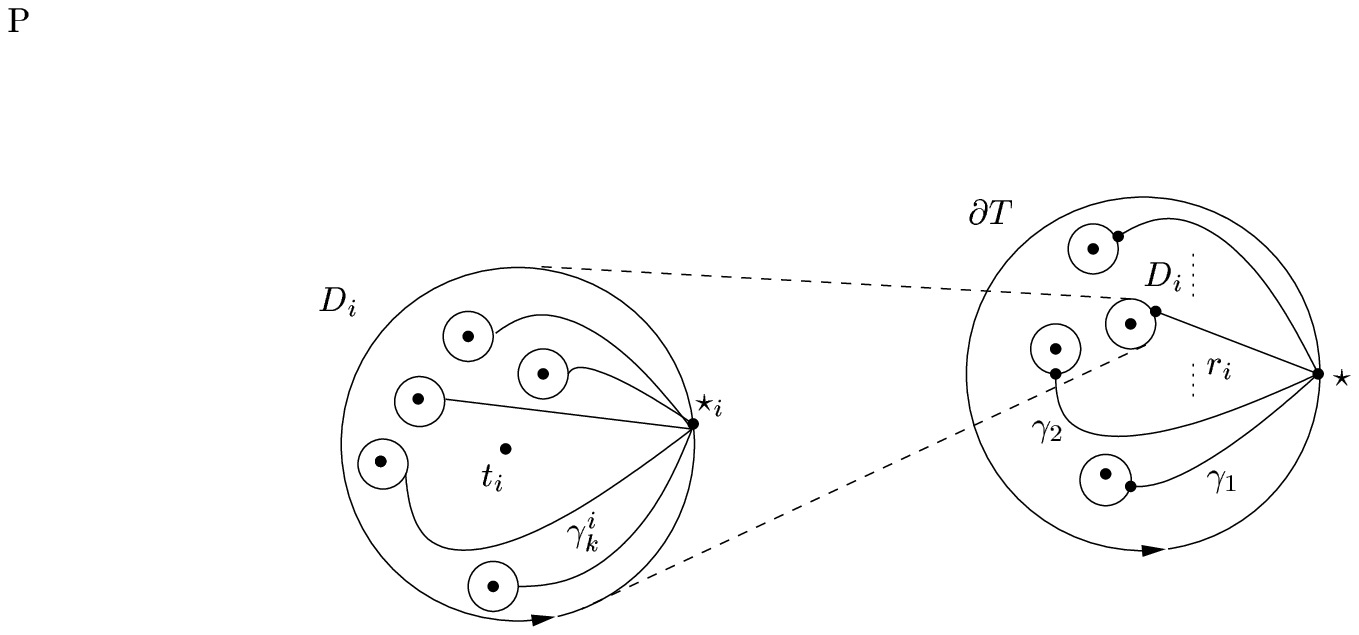}
\caption{}
%%\end{center}
\label{basedistininf}
\end{figure}

We have constructed a practical complete algorithmic method to compute local 
monodromies for a tame polynomial $f$ with $n=2$ (see 
S.~A.~Broughton~\cite[Definition~3.1]{Bro1}). This implies that there exist morsifications 
$\tilde{f}$ of $f$ whose regular fibers are diffeomorphic, and $B_f=\{t_i\}$ contains only critical 
values coming from affine singularities. 
This way we can compute the local monodromies of $f$ from the local monodromies of  $\tilde{f}$.  To obtain the block decomposition of the local monodromies $(h_{\gamma_i})_*$ we need to consider \textit{special geometric bases} of $\pi_1(T\setminus B_{\tilde{f}};\star)$ associated with $(\gamma_i)$.
We can construct this basis as follows. Let $(t_i,\star_i,r_i,S_i,D_i)$ be the elements 
associated with $\gamma_i$ (see Figure~\ref{basedistininf}). 
For a sufficiently small deformation $\tilde{f}$, $B_{\tilde{f}} \cap D_i$ is the set of those critical values which belong to the morsification of the singularities in $f^{-1}(t_i)$. We order the set $B_{\tilde{f}} \cap D_i$ in such a way that the critical values corresponding to the morsification of 
the same critical point in $f^{-1}(t_i)$ are together.
We can consider a geometric basis $(\gamma_k^i)_{k=1,...,k(i)}$ of 
$\pi_1(D_i \setminus  B_{\tilde{f}} \cap D_i;\star_i)$ which respects this order so that 
$(r_i \cdot \gamma_k^i \cdot r_i^{-1})_{k=1,...,k(i)}^{i=1,...,\# B_f}$ is a geometric basis of 
$\pi_1(T \setminus  B_{\tilde{f}};\star)$. Note that, for a fixed $i$, the path 
$r_i \cdot (\prod_{k} \gamma_k^i )\cdot r_i^{-1}$ is homotopic to $\gamma_i$ and hence the ordered 
product of the associated local monodromies gives the block decomposition of $(h_{\gamma_i})_*$ (see Remark~\ref{o:submatriz}). Moreover, the ordered product of the $(h_{\gamma_i})_*$ produces the
monodromy at infinity.

We will then construct a method to compute the local monodromies for a tame Morse function. 
Hence we can assume $f$ to be this tame Morse function. 
Let $\mathfrak{D}_f:=\{(t,x) \in \C^2 \mid \dis_y(f(x,y)-t)=0\}$ be the 
curve discriminant of the polar mapping $\phi_{f,x}:\C^2 \rightarrow 
\C^2,(x,y) \mapsto (t,x):=(f(x,y),x)$. The method needs two data which 
depend on $\mathfrak{D}_f$. The first one is the classical monodromy 
$m$ of the projection $\pi \vert:f^{-1}(\star)\rightarrow \C,(x,y) \mapsto 
x$, in particular we need to know the transpositions given by a 
geometric basis $(\mu_n)$ associated with the ramification points of 
$\pi \vert_{f^{-1}(\star)}$ (the set of $k$ points given by $\mathfrak{D}_f \cap 
\{t=\star\}$).  The second ingredient is the braid monodromy $\nabla$ 
of the  projection $\pi \vert:\mathfrak{D}_f \rightarrow \C,(t,x) \mapsto 
t$, in particular we need to know the braids determined by the 
geometric basis $(\gamma_i)$. The computation of these data can be 
done with the help of computer programs such as 
\cite{bessis:prog} and
\cite{Car03}.
Since our method uses strongly the discriminant curve $\mathfrak{D}_f$ 
we call it \textit{the discriminant method}.

In the following we sketch our method. Let $\mathbb{F}_k$ be the free 
group of $k$ generators: $\mu_1,\ldots,\mu_k$ and let $\B_k$ be the 
braid group on $k$ strings given by the following representation 
$$\left\langle \sigma_1, \ldots, \sigma_{k-1} \mid 
[\sigma_i,\sigma_j]=1 \ \mbox{if} \ \vert i-j \vert \geq 2, 
\sigma_{i+1} \sigma_{i} \sigma_{i+1}= \sigma_{i} \sigma_{i+1} 
\sigma_{i},  i=1,\ldots,k-2\right\rangle . $$ 
Since $f$ is a Morse polynomial, the braid $\nabla(\gamma_i)$ can be 
written as a conjugate of any  basis element $\sigma_j$. Let $\beta_i$ 
be the element which conjugates, for example, $\sigma_1$ (this can be 
obtained from the braid monodromy calculated above). There is a 
natural right action $\Phi:\mathbb{F}_k \times \B_k \rightarrow 
\mathbb{F}_k$ such that 
\begin{equation*}\mu_{j}^{\sigma_i}=\begin{cases}\mu_{i+1} 
&\mbox{if} \  j=i \\\mu_{i+1} \,\cdot \, \mu_i \, \cdot  
\mu_{i+1}^{-1} &\mbox{if} \ j=i+1 \\\mu_j &\mbox{if} \  j \neq 
i,i+1,\end{cases}\end{equation*}
in exponential notation. We compute the elements 
$\mu_1^{\beta_i}$ and $\mu_2^{\beta_i}$ and denote by $\delta_i$ 
their product. Note that $\delta_i$ 
(resp. $\delta_i^{\nabla(\gamma_j)}$) is the image by $\pi \vert_{f^{-1}(\star)}$ of 
a path in $f^{-1}(\star)$ defining a vanishing 1-cycle $\Delta_i$ 
(resp. the Picard-Lefschetz element $(h_{\gamma_j})_*(\Delta_i)$).
We construct an appropriate model of the fiber $f^{-1}(\star)$ 
that enables us to perform algebraic computations. This model 
allows us to express the 1-homology of the fiber $f^{-1}(\star)$ in 
terms of the generators $(\mu_n)$ via the classical monodromy 
calculated above. Finally using the expression of $\delta_i$ and  
$\delta_i^{\nabla(\gamma_j)}$ in the basis $(\mu_n)$ and by means 
of the 1-homology model we can express the cycles $\Delta_i$ and 
$(h_{\gamma_i})_*(\Delta_j)$ algebraically. 

\medskip

We develop the discriminant method in Section~\ref{s:dis_met}. In Section~\ref{modtopyhom}  we construct 
the aforementioned model of the fiber $f^{-1}(\star)$. Afterwards, in Section~\ref{s:dis_met_morse} we explain the method in the Morse function case and state the main Theorem~\ref{t:met_discrim} on their
local monodromy. Finally in Section~\ref{S:pol_tame} we apply the method using generic morsifications for the tame case. Sections~\ref{s:dis_prop} and \ref{s:proof} are mainly technical in order to prove the main theorem and its application to tame polynomials.

In Section~\ref{C:top_conj} we  prove the existence of non-topologically equivalent 
tame polynomials whose equations are conjugated in the number field $\Q(\sqrt{3})$ 
(we thank E.~Artal and A.~Bodin for suggesting these polynomials as worth studying). 
To do this we prove that under certain conditions, if two polynomials are topologically equivalent 
then their local monodromies should be simultaneously conjugated by the same invertible matrix.
Using the discriminant method
we completly compute   local and global blocks in a distinguished basis of vanishing cycles and 
we show that such simultaneous conjugation does not exist. 

Finally Section~\ref{s:base} is devoted to setting definitions, notations and 
properties of geometric bases, braids and Hurwitz moves.

\begin{gracias}
We deeply appreciate the support and valuable ideas suggested by E.~Artal and J.I.~Cogolludo.
\end{gracias}

\section{Geometric bases, braids  and Hurwitz moves}\label{s:base}

All results and definitions in this section can be found,  more extensively,  in~\cite{Artin47} and~\cite{Bir}.
Notations follow~\cite{KJJ}.

\begin{figure}[ht]
%\begin{center}
\vspace*{-1cm}
\hspace*{-3cm}
\includegraphics[scale=.8]{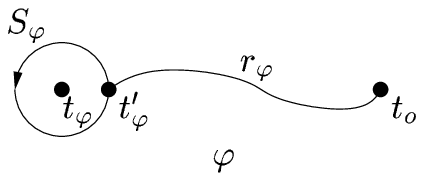}
\caption{}
%%\end{center}
\label{meridiano}
\end{figure}

\begin{definicion}\label{defmeridiano}
Let $S \subset \C$ be a finite set of points and $t_o \in \C \setminus S$. We say  
the path $\varphi$ is a  \emph{meridian in $\C \setminus S$ of the point 
$t_{\varphi} \in S$ and based at $t_o$} (Figure~\ref{meridiano}) if $\varphi$ is the 
union of paths $r_{\varphi} \cdot S_{\varphi}  \cdot r_{\varphi}^{-1}$ with  
$S_{\varphi}$ the positive oriented boundary  of a disk $D_{\varphi}$ with 
center $t_{\varphi}$  such that $D_{\varphi} \cap S =\{t_{\varphi}\}$ and 
$r_{\varphi}$ path from  $t_o$ to  $t'_{\varphi} \in 
S_{\varphi}$ such that $\supp(r_{\varphi}) \cap S=\emptyset$. If $r_{\varphi}$ has not self intersections $\varphi$ is said  
\emph{simple}.\end{definicion}

Let   $\Omega \subset \C$ be a simple compact region, 
$S \subset \intt (\Omega)$ set with $s$ points and  $t_o \in \partial \Omega$. 
We consider special families of bases in the fundamental group
 $\pi_{1}(\Omega \setminus S;\,t_o)$.
\begin{definicion}\label{def:basegeom}
A \emph{geometric basis}   
of the group 
$ \pi_1(\Omega \setminus S;\,t_o)$ is an ordered list $(\varphi_{1},\dots,\varphi_{s})$  such that:
\begin{enumerate}[(i)]
\item\label{def:basegeom.1}  $ \varphi_{i}$ is a simple meridian in $\Omega 
\setminus S$ based at $t_o$ defined by $(t_{\varphi_{i}}, t'_{\varphi_{i}}, 
r_{\varphi_{i}}, S_{\varphi_{i}}, D_{\varphi_{i}})$. 
\item\label{def:basegeom.2}  $\supp(\varphi_{i}) \cap \supp(\varphi_{j})= 
\{t_o\}$ for all $i$, $j$ with $i \neq j$.
\item\label{def:basegeom.3} $\varphi_s \cdot \, \dots \, \cdot \varphi_1$ 
is homotopic to $\partial \Omega$ which is positively oriented (product from left to right).
\end{enumerate}
\end{definicion}
 
\begin{ejemplo}\label{ejemplobase}
Let $X$ be a geometric disk  centered at $0 \in \C$ such that
$\{1,2, \ldots, k\} \subset \Int(X)$. Let $* \in \partial X \cap \R^+$.
Figure~\ref{basefundamental} shows a list of meridians
$(\mu_1,\dots,\mu_k)$ defining a geometric basis of $\pi_1(X \setminus \{1,2, \ldots, k\};*)$. 
\begin{figure}[ht]
%\begin{center}
\vspace*{-3cm}
\hspace*{-2.5cm}
\includegraphics[scale=.8]{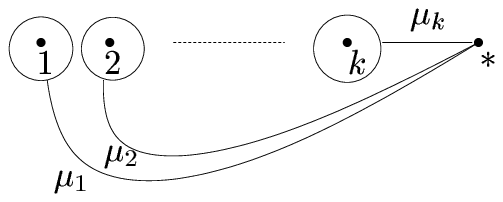}
%\vspace*{-5.5cm}
\caption{}
%\end{center}
\label{basefundamental}
\end{figure} 
\hfill
$\blacksquare$
\end{ejemplo}

\begin{definicion}\label{d:grup_tren}
Let $\mathcal{V}$ be the topological space of sets with $k$ different points in~$\C$. We write $\B({\bf s}^{1},{\bf s}^{2}), {\bf s}^{1},{\bf s}^{2} \in 
\mathcal{V}$, the set of braids which start  at ${\bf s}^{1}$ and end at 
${\bf s}^{2}$. As it is usual we denote simply by $\B_k$ the group $\B(\{1, \ldots, 
k\},\{1, \ldots, k\})$.\end{definicion}

\begin{proposicion}\label{p:grupo_tren}
The group  $\B_{k}$ is generated by elements 
$\sigma_j,j=1,\ldots,k-1$ (Figure~\ref{sigmaj}) with relations
$ [\sigma_i,\sigma_j]=1 \ \mbox{if} \ \vert i-j \vert \geq 2 \ \mbox{and} \ 
\sigma_{i+1} \sigma_{i} \sigma_{i+1}= \sigma_{i} \sigma_{i+1} \sigma_{i},$  $i=~1,\ldots,k-2. $
\begin{figure}[ht]
%\begin{center}
\vspace*{-1.5cm}
\hspace*{-3cm}
\includegraphics[scale=.7]{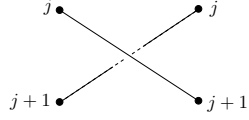}
\caption{Braid $\sigma_j$.}
%\end{center}
\label{sigmaj}
\end{figure} 
\end{proposicion}

\begin{definicion}\label{d:conj_base_tren}
A  \textit{base set of braids}  $\mathfrak{O}$ is a set such that for every
$\mathbf{s} \in \mathcal{V}$ there exists a unique braid 
$\tau_{{\bf{s}},\mathfrak{O}} \in \mathfrak{O}$ 
with $\tau_{{\bf{s}},\mathfrak{O}} \in \B(\{1, \dots, k\}, {\bf{s}})$.
\end{definicion}

\begin{notacion}\label{not:orden}
Let $\psi \in (0,\pi)$.
\begin{itemize}
\item Define by  $
\pi_{\psi}(z):=\Re{(z)} -(\Im{(z)}/\tan(\psi))
$
a certain projection of $\C$ onto $\R$. Let
 the lexicographic order $<^{\psi}$ in $\C$ (Figure~\ref{f:orden})
 given by 
\begin{equation*}%\label{orden}
z <^{\psi} w \Longleftrightarrow 
\pi_{\psi}(z) < \pi_{\psi}(w) \ \mbox{ or} \
\pi_{\psi}(z) = \pi_{\psi}(w) \ \mbox{and} \ 
\Im{(z)}< \Im{(w)}.
\end{equation*}

\begin{figure}[ht]
\hspace*{-4cm}
\includegraphics[scale=.7]{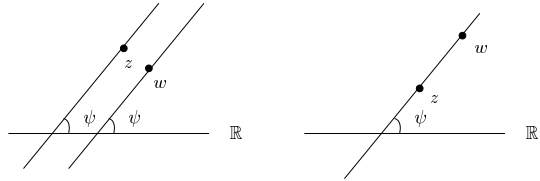}
\caption{The two cases $z<^{\psi} w$.}
\label{f:orden}
\end{figure}

\item For ${\bf s} \in \mathcal{V}$ we consider the $k$-tuple ${\bf s}(\psi):=(s_1, \dots, s_{k}), s_i \in {\bf s},$   
such that it verifies $s_1 <^{\psi} s_2 <^{\psi} \ldots <^{\psi} s_{k}$.

\item We denote by $\mathcal{V}(\psi)$ the set of  lists ${\bf s}(\psi),{\bf s} \in \mathcal{V}$.
 
 \item We denote by $\mathfrak{O}(\psi)$ a base set of braids  such that every $\tau_{{\bf 
s}, \mathfrak{O}(\psi)}$ has a set of strands $\{\alpha_1, \ldots, \alpha_k\}$ 
with $(\alpha_1(u), \dots,\alpha_{k}(u)) \in  \mathcal{V}(\psi), \forall u \in 
[0,1],$  and 
$\alpha_i(0)=~i,$ $\alpha_i(1)=s_i$. For instance 
$\alpha_i(u)=s_i \cdot u+i \cdot (1-u), i=1,\ldots,k$.
\end{itemize}
\end{notacion}

Let ${\bf{s}}^{1},{\bf{s}}^{2} \in \mathcal{V}$ and $\mathfrak{O}$ a base set of braids.
One can define the bijection\begin{equation}\label{iso1}
\begin{CD}
\B({\bf{s}}^{1},{\bf{s}}^{2}) @>\I_{{\bf s}^{1},{\bf s}^{2}, \mathfrak{O}}>> 
\B_{k}, \ \tau  \mapsto \tau_{{\bf{s}}^{1},\mathfrak{O}} \, \tau \, 
 \tau_{{\bf{s}}^{2},\mathfrak{O}}^{-1}=: \tau(\mathfrak{O})
\end{CD}
\end{equation}

\begin{proposicion}\label{ident_trenzas}
\hfill

\begin{enumerate}[(a)]
\item\label{ident_trenzas.1}  $\I_{\bullet,\bullet, 
\mathfrak{O}}$  is an isomorphism of grupo\"ids.

\item\label{ident_trenzas.2} Let  $\mathfrak{O}$ and $ \bar{\mathfrak{O}}$ 
base sets of braids  and $\tau \in \B({\bf{s}}^{1},{\bf{s}}^{2})$
then
\begin{equation*}
\tau(\bar{\mathfrak{O}})=  ( \tau_{{\bf{s}}^{1},\bar{\mathfrak{O}}} \, 
\tau_{{\bf{s}}^{1},\mathfrak{O}}^{-1} ) \, \tau(\mathfrak{O}) \,
( \tau_{{\bf{s}}^{2},\bar{\mathfrak{O}}} \, 
\tau_{{\bf{s}}^{2},\mathfrak{O}}^{-1} )^{-1}.
\end{equation*}
In fact if $\bf{s}^{1}=\bf{s}^{2}$, braids 
$\tau(\mathfrak{O})$ and $\tau(\bar{\mathfrak{O}})$
are conjugated in $\B_{k}$.
\end{enumerate}
\end{proposicion}

Let   $\tau \in \B({\bf s}^{1},{\bf s}^{2})$,   $X$ a disk such that
${\bf s}^{1} \cup \, {\bf s}^{2}  \subset \Int(X)$ and $* \in \partial X$.
\begin{definicion}\label{d:mov_H}
 A map $H:\C \times [0,1] \rightarrow \C$ is said 
 a \emph{$({\bf s}^{1},{\bf s}^{2})$-special isotopy associated with $\tau$}  
if  $H(\bullet,u)$, $u$ fix, is a preserving orientation homomorphism so that 
is the identity in $\C \setminus \intt(X)$ and
 $H \vert_{ {\bf s}^{1} \times [0,1]}$ defines the braid $\tau$. 
\end{definicion}

\begin{definicion}
We define the \emph{Hurwitz move associated with $\tau$} as the isomorphism 
$$
\Psi_{\tau}:\pi_1(X \setminus {\bf s}^{1};*) \rightarrow \pi_1(X \setminus 
{\bf s}^{2};*), \varrho\mapsto \Psi_{\tau}(\varrho)=:\varrho^{\tau}
$$
 which defines a $({\bf s}^{1},{\bf 
s}^{2})$-special isotopy associated with $\tau$ and we denote by 
\begin{equation*}
 \Phi_{{\bf s}^{1},{\bf s}^{2}}:  \pi_1(X \setminus {\bf s}^{1};*) 
\times \B({\bf s}^{1},{\bf s}^{2}) \rightarrow \pi_1(X \setminus {\bf 
 s}^{2};*) \,
\end{equation*} 
 the right grupo\"id action defined by Hurwitz moves.
\end{definicion}

\begin{ejemplo}\label{ej:fund}
Let us consider the situation given in Example~\ref{ejemplobase}. The action $\Phi_{{\bf 
x},{\bf x}}$ is defined by	
\begin{equation}\label{formulaexp}
\mu_{j}^{\sigma_i}=
\begin{cases}
\mu_{i+1} &\mbox{if} \  j=i \\
\mu_{i+1} \,\cdot \, \mu_i \, \cdot  \mu_{i+1}^{-1} &\mbox{if} \ j=i+1 \\
\mu_j &\mbox{if} \  j \neq i,i+1.
\end{cases}
\end{equation}
Also note that  for any base set of braids  $\mathfrak{O}$  we have the commutative 
diagram \begin{equation}\label{diagramaprimario}
\begin{CD}
\pi_1(X\setminus {{\bf s}^{1}};*)  \times  \B({{\bf s}^{1}},{{\bf s}^{2}}) 
@>\Phi_{{\bf s}^{1},{\bf s}^{2}}>> \pi_1(X\setminus {{\bf s}^{2}};*) 
\\@V\Psi_{\tau_{{\bf s}^{1},\mathfrak{O}}}^{-1} \times  \I_{{\bf s}^{1},{\bf 
s}^{2}, \mathfrak{O}}VV @AA\Psi_{\tau_{{\bf  s}^{2},\mathfrak{O}}}A \\\pi_1(X 
\setminus \{1,\ldots,k\}; *)  \times  \B_{k} @>\Phi_{{\bf x},{\bf  x}}>>\pi_1(X 
\setminus \{1,\ldots,k\}; *)\end{CD}\begin{picture}(0,0)
%\graphpaper(-150,-50)(100,75)
\put(-115,-4){$\circlearrowleft$}
\end{picture}
\end{equation}
\hfill $\blacksquare$
\end{ejemplo}

From Proposition~\ref{p:grupo_tren} and action~\eqref{formulaexp}  we have the following lemma.

\begin{lema}\label{l:reduc_1}
Let $\alpha_s:=(\sigma_2 \sigma_1) (\sigma_3 \sigma_2) \dots (\sigma_s \sigma_{s-1}) \in \B_k, s=2, \dots,k-1$. Then:
\begin{enumerate}[(a)]
\item\label{l:reduc_1.a} $\sigma_s=\alpha_s^{-1} \sigma_1 \alpha_s$.

\item\label{l:reduc_1.b} $\mu_{s}=\mu_1^{\alpha_s}$ and $\mu_{s+1}=\mu_2^{\alpha_s}$.
\end{enumerate}
\end{lema}

\section{Discriminant properties}\label{s:dis_prop}

Let $f \in \C[x,y]$ be a monic polynomial of degree $N$ in the variable $y$ with 
isolated singularities. Let
\begin{equation*}
\phi_{f,l}:\C^2 \longrightarrow   \C^{2}, \ (x,y)  \mapsto (f(x,y),l(x,y))
\end{equation*}
the polar map associated with $f$. We take the linear form $l(x,y)=x$.
In this case the discriminant curve $\mathfrak{D}_f$ of $\phi_{f,x}$ (i.e., the image of the critical locus of $\phi_{f,x}$) has the 
equation\begin{equation*}%\label{e:def_df}
d_f(t,x) := \, \mbox{discrim}_{y}(f(x,y)-t)  \in \C[t,x].
\end{equation*}

\begin{notacion}
We denote by $\mu_{t_o,x_o}$ the sum of  Milnor numbers of the singularities in $f^{-1}(t_o)$ 
whose projection onto the first variable $x$ is $x_o \in \C$.
\end{notacion}
\begin{lema}\label{desc_disc}
Let $r$ be the number of different roots of $f(0,y)$. Then the discriminant $d_f(t,x)$ is $x$-regular 
of degree $\mu_{0,0}+ N - r$ and $t$-regular of degree  $N-r$.
\end{lema}
\begin{proof}
Let $y_j \in \C, j=1, \dots,r,$ root of $f(0,y)$ with multiplicity $N_j$. 
By Hensel's Lemma~\cite[page 44]{GR} we write 
$
f(x,y)= \prod_{j=1}^{r} P_j(x,y), x\in \D, 
$
with $\D$  small disk 
centered at $0 \in \C$ and  $P_j(x,y)$  Weierstrass polynomial centered at $(0,y_j)$ of degree 
$N_j$ in the variable $y$. 
As $f^{-1}(0)$ has isolated singularities then $P_j$ is reduced and so 
$\dis_y(P_j)$ is not zero and we deduce that $t$ does not divide   
$\dis_y(P_j-t)$. Then $\dis_y(P_j-t)$ is $x$-regular 
of degree~\cite[page~317]{Te}\begin{equation}\label{e:P_j}
\mu(P_j)+N_j-1.
\end{equation} 
By  properties of the resultant we have
\begin{equation*}
 \dis_y(P_iP_j) = (\res_y(P_i,P_j))^2 \dis_y(P_i) \dis_y(P_j) 
\end{equation*}
where $\res_y(P_i,P_j) \in \C\{x\}$ is a unit because $P_i$ and $P_j$ do not 
share their Puiseux series. Then we deduce that the degree of $d_f(0,x)$ is the sum of 
the degrees of $\dis_y(P_j), j=1,  \dots,r$, and finally by~\eqref{e:P_j}   
$d_f(t,x)$ is $x$-regular of degree\begin{equation*}
\sum_{j=1}^{r}(\mu(P_j)+ N_j - 1) = \mu_{0,0}+ N - r.
\end{equation*}

For the second part of the lemma we note that $f(0,y)$ and $f'(0,y)$ share 
$y_j$ as root of multiplicity $N_j-1$. Then
$
d_f(t,0)
$
has  $0$ as a root 
of multiplicity
$
\sum_{j=1}^{r} (N_j-1)=N-r.
$
Then $d_f(t,x)$ is $t$-regular of degree $N-r$.
\end{proof}

\begin{corolario} \label{propdiscrim}
Under the notation of Lemma~\ref{desc_disc} the germ $(\mathfrak{D}_{f},0)$ is empty if and only if   
$N=r$, it is smooth if and only if  $N=r+1$ and it  is singular if and only if  $N>r+1$. 
Moreover, in the case $N \geq r+1$ the line $\{t=0\}$ belongs to its tangent cone if and only 
if $\mu_{0,0}>0$.
\end{corolario}
\begin{corolario}\label{c:iden_vc}
The point $0 \in \C$ is a critical value of 
$f$ if and only if $\{t=0\}$ belongs to the tangent cone of 
$\mathfrak{D}_f$.
\end{corolario}
\begin{proof}The point $(p_x,p_y) \in f^{-1}(0)$ is a critical 
point if and only if $\mu_{0,p_x}>0$ and by Corollary~\ref{propdiscrim} (at the 
point $(0,p_x)$ by change of variable) if and only if $\{t=0\}$ belongs to the tangent cone of 
$\mathfrak{D}_f$.
\end{proof}
\begin{corolario}\label{c:Morse}
Let assume that the line $\{t=0\}$ is tangent at only one point $p$ of $\mathfrak{D}_f$.
If $p$ is smooth and $(\mathfrak{D}_f,\{t=0\})_p=2$ then  $f^{-1}(0)$ only has one  singular point which is  no degenerate.
\end{corolario}
\begin{proof}
Without loss of generality we can suppose that 
$p=0 \in \C^2$.  
As multiplicity of 
$\{t=0\}$ in $\mathfrak{D}_f$ is~2 then  $d_f(t,x)$  is $x$-regular of degree 2. Also
as  $0 \in 
\mathfrak{D}_f$ is a smooth no inflection point by Corollary~\ref{propdiscrim} we have $N=r+1$ then  by Lemma~\ref{desc_disc}  we have
$
\mu_{0,0}+N-r=2 
$
so $\mu_{0,0}=1$. Then there is 
only one  no degenerate singular point in $f^{-1}(0)$ with projection onto
the first variable  $0 \in \C$.  Finally, by Corollary ~\ref{propdiscrim}, 
there is not other singular point in $f^{-1}(0)$ because $\{t=0\}$ is only 
tangent at  $0 \in \mathfrak{D}_f$.
\end{proof}

\begin{lema}\label{p:buenas_coor_abier}
Let us denote by
$
f_{\lambda}(x,y):=f(x+ \lambda y, y).
$
For almost every $\lambda \in \C$ the discriminant  
$\mathfrak{D}_{f_{\lambda}}$ is a reduced curve in $\D \times \C$ with $\D$ 
small disk centered at~$0 \in \C$. 
\end{lema}
\begin{proof}
The polynomials $f_{\lambda}(x,y)$ have only  isolated 
singularities in the zero fiber for every $\lambda \in \C$. We consider  
$
d_{f_{\lambda}}(t,x) \in \C[t,x,\lambda]  
$
as a polynomial in the variable~$x$. By 
Lemma~\ref{desc_disc} this  polynomial is not zero for any $(t,\lambda) \in 
\D \times \C$ with $\D$ small disk centered at~$0 \in \C$. We reduce $\D$ 
if necessary such 
that $f^{-1}(t')$ is smooth for every $t' \in \D \setminus \{0\}$.
We consider
$
\dis_{x}(d_{f_{\lambda}}(t,x)) \in \C[t, \lambda]
$
as a polynomial in the variable $\lambda$. If this polynomial is not zero in  $\D$ 
then $\lambda_o \in \C$ is a root (for every $t\in \D$) if and only if the 
polynomial $d_{f_{\lambda_o}}(t,x)$ is not reduced. Then for $\lambda \neq  
\lambda_o$ the discriminant  $\mathfrak{D}_{f_{\lambda}}$ is reduced in 
$\D \times \C$. Therefore lemma   will be proved if  it exists a 
point $t_o \in \D \setminus \{0\}$ such that $\dis_{x}(d_{f_{\lambda}}(t_o,x))$ is 
not zero.
 
Let $F_{\lambda,t_o}$ homogeneous polynomial  associated with 
$f_{\lambda}-t_o, t_o \in \D \setminus \{0\}$. The projective curve 
$\mathcal{Z}(F_{\lambda,t_o})$ is  smooth except for maybe at  a finite number of points on the 
line at infinity $\{z=0\}$.
We can decompose $\mathcal{Z}(F_{0,t_o})=\mathcal{L} \cup \mathcal{C}$ where $\mathcal{L}$ is the line components in 
 $\mathcal{Z}(F_{0,t_o})$.
The projective curve $\mathcal{C}$ (without lines as components) has a finite number $\mathcal{L}'$ of 
tangent lines  at inflection points and bitangent 
at smooth points
(see, 
for example,~\cite[Theorem 1, Section~7.3]{BK86}). 
Let $\{a_ix+b_iy+c_iz=0\}$  the union $\mathcal{L}' \cup \mathcal{L}$. Then 
only the lines $\{a_i(x+\lambda y) +b_iy+c_iz=0\}$ are tangent lines at 
inflection points or bitangent at smooth points of 
$\mathcal{Z}(F_{\lambda,t_o})$ for every  $\lambda \in \C$. 
Let $\lambda \neq 
-\frac{b_i}{a_i}$. The lines $\{ax+cz=0\}$ are not tangent lines at inflection 
points or bitangent at smooth points of $\mathcal{Z}(F_{\lambda,t_o})$, so the 
line $\{ax+c=0\}$ cut to the affine curve  $f_{\lambda}^{-1}(t_o)$ in a traversal 
way or it is tangent at only one smooth no inflection point. Then the projection 
of $f_{\lambda}^{-1}(t_o)$ onto the  variable $x$ has as  ramification points  
$x_o \in \C$  such as the line $\{x=x_o\}$ is tangent at only one smooth 
point of $f_{\lambda}^{-1}(t_o)$ and so  $d_{f_{\lambda}}(t_o,x)$
is reduced and  $\dis_{x}(d_{f_{\lambda}}(t_o,x))\neq 0$.
\end{proof}

\section{Discriminant method}\label{s:dis_met}
Let $f$ be as  at the beginning of Section~\ref{s:dis_prop}. 
Let $P$ be  the set of isolated singularities of $f$. We  fix  $T$ a geometric disk in 
$\C$ such that $f(P) \subset \Int(T)$ and  
 $
 M:=X \times Y 
 $
a 2-dimensional polydisk with $X$ and $Y$ geometric disks such as:
 \begin{enumerate}
\item\label{s:dis_met_cond1} $P \subset \Int(X) \times \Int(Y) \ \mbox{and} \  
\partial M \cap f^{-1}(T) \subset \partial X  \times \intt(Y).$
\item\label{s:dis_met_cond2} For every $t \in T$  the projection
\begin{equation}\label{proy_prim}
\pi \vert_{V_{t}(f)}: V_{t}(f) 
\longrightarrow  X,  \ (x,y) \mapsto x,
\end{equation}
with $V_{t}(f):=f^{-1}(t) \cap M$, is a covering map with $N$ sheets and ramification 
points $ {\bf x}^{t}:=\{x^{t}_1, \dots,x^{t}_{k(t)}\} \subset \Int(X), k(t) 
\in \N$.
\item \label{s:dis_met_cond3}
${\bf x}:=\{1,\ldots,k\} \subset \Int(X)$ where $k$ is the degree of the projection 
\begin{equation}\label{e:proy_Df}
\pi \vert_{{\mathfrak{D}_{f}}} : {\mathfrak{D}_{f}}  \longrightarrow  T,\ (t,x) 
\mapsto  t,
\end{equation} 
 ramified in a set $\mathcal{T}$ (note that $k(t)=k$ if $t \in T \setminus \mathcal{T}$).
\end{enumerate}

\begin{nota}\label{n:modif_hip2cap2} 
We are interested   in the discriminant curve  in
 $T \times X$ then as~\eqref{proy_prim} has its ramification points in 
$\Int(X)$, we write~$\mathfrak{D}_{f}$  as  $\mathfrak{D}_{f} \cap (T \times X)$.
Also  ${\bf x}\subset \Int(X)$ is a technical hypothesis  due to taking $\B_k$ as
our base group of braids.
\end{nota}

By Corollary~\ref{c:iden_vc} the projection of the set of points 
in ~$\mathfrak{D}_f$ with vertical tangent line is~$f(P)$ then $f(P) \subset \mathcal{T}$. 
Finally by Note~\ref{n:modif_hip2cap2} 
$\mathfrak{D}_f$ has no vertical asymptotes so
$
\mathcal{T}= f(P) \cup \mathcal{T}_{Sing}
$
with $\mathcal{T}_{Sing}$ 
projection onto the first variable of Sing$(\mathfrak{D}_f)$.

Let $t \in  T \setminus \mathcal{T}$. Associated with~\eqref{e:proy_Df} we can define the action
\begin{equation}\label{n:pi_bar}
(\bar{\pi}\vert_{\mathfrak{D}_{f}})_*:\pi_1(T\setminus \mathcal{T}; t)
\longrightarrow  \B({\bf x}^{t},{\bf x}^{t})  
 \end{equation}
given by  $\pi\vert^{-1}_{\mathfrak{D}_f}(t')=\{t'\} \times 
\bar{\pi}\vert_{\mathfrak{D}_f}(t'), t'\in \C,$. If we fix $\mathfrak{O}$ a base set of braids
we have by composition the braid monodromy
\begin{equation}\label{n:mono}
 \begin{CD} 
\nabla_{t,\mathfrak{O}}:\pi_1(T\setminus \mathcal{T}; t)
 @>(\bar{\pi}\vert_{\mathfrak{D}_{f}})_*>> 
 \B({\bf x}^{t},{\bf x}^{t}) 
 @>\I_{\mathfrak{O},{\bf x}^{t},{\bf x}^{t}}>> 
  \B_k
 \end{CD}
 \end{equation}
The indices $\mathfrak{O}$ and $t$ will be dropped if they are clear from the context.

\subsection{A model of the regular fiber.}\label{modtopyhom}
Let us consider $t \in T \setminus  \mathcal{T}$ and $(\varphi_1, \ldots, \varphi_k)$ 
geometric basis of $\pi_1(X \setminus {\bf x}^t;*), * \in \partial X$. 
We lift in $V_t(f)$ (according to~\eqref{proy_prim}) the sets
\begin{equation*}
\bigcup_{1 \leq i \leq k} \supp(\varphi_i) \cup 
D_{\varphi_i}\setminus \{ x_{\varphi_i}\} \ \mbox{and} \ X \setminus \bigcup_{1 
\leq i \leq k} ( \supp(\varphi_i) \cup D_{\varphi_i})
\end{equation*}    
and we glue them with the lifting of $\varphi_k \cdot \, \dots \, \cdot 
\varphi_1$. With this we obtain a topological model of 
$V_{t}(f) \setminus \pi\vert^{-1} _{V_{t}(f)}({\bf x}^t)$ that can be 
extended to $V_{t}(f)$ by gluing the lifting of~${\bf x}^t$.

\begin{notacion}
We  denote by  $m_t:\pi_1(X\setminus {\bf x}^t;*) \rightarrow \Sigma_N$ 
the classical monodromy associated with~\eqref{proy_prim} and by 
$l^{m_{t}(\alpha)}$ the image of $l \in \{1, \ldots, N\}$ by the permutation 
$m_{t}(\alpha), \alpha \in \pi_1(X \setminus {\bf x}^t;*)$. 
\end{notacion}

In this subsection we consider the following hypothesis.

\begin{hipotesis}\label{h:tras}
The permutations $m_{t}(\varphi_i)$ are transpositions $(a_i, 
b_i)$ with $a_i \neq b_i$ and $1 \leq a_i, b_i \leq N$. 
\end{hipotesis}

\begin{notacion}\label{notelevacion}
For every $i=1, \dots,k$ and $l=1, 
\dots, N$ we denote by:
\begin{enumerate}[(1)]
\item $(*,y_1), \dots, (*,y_N)$ the lifting of $*\in 
\partial X$ in $V_{t}(f)$. 
\item  $(\alpha)_l$ the lifting  based at~$(*,y_l)$ 
of the path~$\alpha$   with support in $ X \setminus {\bf x}^t$ and 
based at $* \in \partial X$.
\item $(D_{i})_l$ the compact component of the lifting of
$\supp(\varphi_i) \cup D_{\varphi_i}$ in $V_{t}(f)$
which has  the point 
$(*,y_l),  l \neq a_i, b_i,$ and $(D_{i})_{\bar{l}}$ with 
$\bar{l}:=\min\{a_i, b_i\}$ (minimum by convention) which has  the 
points $(*,y_{a_i})$ and  $(*,y_{b_i})$.
\end{enumerate}
\end{notacion}

The topological model of $V_{t}(f)$ defines  a 2-dimensional CW-complex $\mathfrak{K}$
 with chain groups 
 \begin{equation*}
C_0(\mathfrak{K}) :=<[1], 
\dots, [N]>, 
C_1(\mathfrak{K}):=<(\varphi_{i})_l \mid i=1,\dots,k \, 
, l=1, \dots, N> \ \mbox{and}
\end{equation*} 
\begin{equation*}
C_2(\mathfrak{K}) :=<(D_{i})_l \mid i=1,\dots, k, 
\, , l=1, \dots, N, \, l \neq \max\{a_i, 
b_i\}>
\end{equation*} 
and  boundary operators $\partial_1$ and $\partial_2$ defined by
\begin{equation}\label{opborde1}
\partial_1((\varphi_i)_l) := 
[l^{m_{t}(\varphi_i) }]-[l]=
\begin{cases}
0 &\mbox{if $l \neq a_i,b_i$} \\
[b_i]- [a_i] &\mbox{if $l=a_i$} \\
[a_i]- [b_i] &\mbox{if $l=b_i$}
\end{cases}
\end{equation}
\begin{equation}\label{opborde2}
\partial_2(D_{i,l}) :=
\begin{cases}
(\varphi_i)_l &\mbox{if $l \neq \min\{a_i, b_i\}$} \\
(\varphi_i)_{a_{i}}+(\varphi_i)_{b_{i}} &\mbox{if $l=\min\{a_i, 
b_i\}$} \end{cases}
\end{equation}

We obtain the homology  
$H_1(V_{t}(f)):=Z_1(\mathfrak{K})/B_1(\mathfrak{K})$ where $Z_1(\mathfrak{K}):= \mbox{Ker} \, \partial_1$ and 
$B_1(\mathfrak{K}):= \mbox{Im} \, \partial_2$. By~\eqref{opborde2} it has the following  relations
\begin{equation}\label{rel_grup_homolog}
(\varphi_i)_l=
0  \ \mbox{if $l \neq a_i, b_i$} \ \mbox{and} \
(\varphi_i)_{a_{i}}+(\varphi_i)_{b_{i}}=0.  
\end{equation}

\begin{lema}\label{lem:elevcaminos}
The path $ (\varphi_{i}^{\epsilon})_{l}, \epsilon \in\{1,-1\},	$ is such 
that $\partial_1((\varphi_{i}^{\epsilon})_{l})=[l^{m_t(\varphi_i)}]- [l]$ and 
its class in $C_1(\mathfrak{K}) /  B_1(\mathfrak{K})$ is $(\varphi_{i})_l$.
\end{lema}
\begin{proof} 
By uniqueness of lifting 
  \begin{equation*}%\label{e:menos}
(\varphi_i^{-1})_l=
\begin{cases}
(\varphi_{i})_l^{-1} &\mbox{if} \ l \neq a_i, b_i\\
(\varphi_{i})_{b_i}^{-1} &\mbox{if} \ l= a_i\\
(\varphi_{i})_{a_i}^{-1} &\mbox{if} \ l= b_i\\
\end{cases}
\end{equation*} 
Then the path  $ 
(\varphi_{i}^{\epsilon})_{l} $ as element in 
$C_1(\mathfrak{K})$ is
\begin{equation*} 
(\varphi_{i}^{\epsilon})_{l}=
\begin{cases}
\epsilon \cdot (\varphi_{i})_{b_{i}} &\mbox{if}  \ \epsilon=-1 \ \mbox{and} \ l= a_i\\
\epsilon \cdot (\varphi_{i})_{a_{i}} &\mbox{if}  \ \epsilon=-1 \ \mbox{and} \ l= b_i\\
\epsilon \cdot (\varphi_{i})_{l} & \mbox{other cases.}
\end{cases}
\end{equation*}
The lemma is deduced by computing $\partial_1$ (as described in~\eqref{opborde1}), 
the equality 
$m_{t}(\varphi_{i}^{\epsilon})= m_{t}(\varphi_{i})$ 
and using the relations~\eqref{rel_grup_homolog}.
\end{proof}

\begin{proposicion}\label{propcamhomolog}
Let $\alpha=\varphi_{i_1}^{\epsilon_{1}} \, 
\cdot \, \dots \, \cdot \, \varphi_{i_n}^{\epsilon_{n}}$,  
$\epsilon_{s} \in \{1, -1 \}$ expressed in the geometric basis and $l \in \{1, \dots, 
N\}$. If we write in a recursive way
\begin{equation*}%\label{def_l}
l(s+1):=
\begin{cases}
l &\mbox{if} \  s=0 \\
l(s)^{m_{t}(\varphi_{i_s}^{\epsilon_{s}})}=
l(s)^{ m_{t}(\varphi_{i_s})} &\mbox{if} \   1 \leq s \leq n
\end{cases}
\end{equation*}
then 
\begin{enumerate}
[(a)]\item\label{propcamhomolog.a} $(\alpha)_l=(\varphi_{i_1})_{l(1)} + (\varphi_{i_2})_{l(2)} +  \dots  + 
(\varphi_{i_n})_{l(n)}   
$   in $C_1(\mathfrak{K}) /  B_1(\mathfrak{K})$.
\item\label{propcamhomolog.b} $(\alpha^{-1})_{l^{m_{t}(\alpha) }}=-(\alpha)_l$ 
 in $C_1(\mathfrak{K}) /  B_1(\mathfrak{K})$.
\item\label{propcamhomolog.c}  $(\alpha)_l \in Z_1(\mathfrak{K})$ if and only if
$l^{m_{t}(\alpha) }=l$.\end{enumerate}
\end{proposicion}

\begin{proof}
Since $(\alpha)_l  =  \prod_{s=1}^{n}(\varphi_{i_s}^{\epsilon_{s}})_{l(s)}$, then 
by the second part of Lemma~\ref{lem:elevcaminos} we have
$(\alpha)_l= \sum_{s=1}^{n} (\varphi_{i_s})_{l(s)}$  as an element of 
$C_1(\mathfrak{K}) /  B_1(\mathfrak{K})$ so we have~\eqref{propcamhomolog.a}.

To prove~\eqref{propcamhomolog.b} we proceed by induction. If $n=1$ 
then~$\alpha=\varphi_{i}^{\epsilon}$.  By the second part of 
Lemma~\ref{lem:elevcaminos}  and  the equality 
$m_{t}(\varphi_{i}^{\epsilon})= m_{t}(\varphi_{i})=(a_i,b_i)$, 
we have   in $C_1(\mathfrak{K}) /  B_1(\mathfrak{K})$
\begin{equation*}
(\varphi_{i}^{-\epsilon})_{l^{m_{t}(\varphi_{i}^{\epsilon})  }}= 
(\varphi_{i})_{l^{m_{t}(\varphi_{i}) }}=\begin{cases}
0=-(\varphi_{i})_{ 
l}, l \neq a_i,  b_i \\
(\varphi_{i})_{b_i}=-(\varphi_{i})_{a_i}=-(\varphi_{i}^{\epsilon})_{a_i}, 
 l=a_i \\
(\varphi_{i})_{a_i}=-(\varphi_{i})_{b_i}=-(\varphi_{i}^{\epsilon})_{b_i},  
l=b_i \end{cases}
\end{equation*}
 Then \eqref{propcamhomolog.b} is true for~$n=1$. We suppose 
\eqref{propcamhomolog.b} is true for~$n-1$ 
and we denote $\tilde{\alpha}:=\varphi_{i_1}^{\epsilon_{1}} \, \cdot \, \dots 
\, \cdot \, \varphi_{i_{n-1}}^{\epsilon_{{n-1}}}$. Then
\begin{equation*}
\begin{array}{rcl}
\hspace*{-.5cm} (\alpha^{-1})_{l^{ 
m_{t}(\alpha) }} &=& 
(\varphi_{i_n}^{-\epsilon_{n}})_{l(n)^{  
m_{t}(\varphi_{i_n})}} \cdot \left(  
\tilde{\alpha}^{-1}\right)_{l^{ m_{t}(\tilde{\alpha})   
}}	\\
\end{array}
\end{equation*}  
and by hypothesis of induction we have  in
$C_1(\mathfrak{K}) /  B_1(\mathfrak{K})$
\begin{equation*}
\begin{array}{rcccc}
\hspace*{-.5cm} (\alpha^{-1})_{l^{ m_{t}(\alpha) }} &=&    
(\varphi_{i_n}^{-\epsilon_{n}})_{l(n)^{ 
m_{t}(\varphi_{i_n})  }} + \left( 
\tilde{\alpha}^{-1}\right)_{l^{ m_{t}(\tilde{\alpha})  }}
&=& 
-(\varphi_{i_n}^{\epsilon_{n}})_{l(n)} - ( \tilde{\alpha})_{l}=
-(\alpha)_l.
\end{array}
\end{equation*}
Hence  \eqref{propcamhomolog.b} is true.

To prove \eqref{propcamhomolog.c}, from the first part of 
Lemma~\ref{lem:elevcaminos}, we have
\begin{equation*}
\partial_1((\alpha)_l)= \sum_{s=1}^{n} 
\partial_1((\varphi_{i_s}^{\epsilon_{s}})_{l(s)})= \sum_{s=1}^{n} ([l(s+1)]- 
[l(s)])=[l(n+1)]- [l(1)]
\end{equation*}
 then $(\alpha)_l \in Z_1(\mathfrak{K})$ if 
and only if $l(n+1)=l(1)=l$. Finally since $
l(n+1)  =  l^{m_{t}(\alpha) }$ we 
proof~\eqref{propcamhomolog.c}.\end{proof}

\begin{notacion}\label{n:convencion}
Let $\mathfrak{O}$ a 
base set of braids. Since $\Psi_{\tau_{{\bf x}^t,\mathfrak{O}}}$ is an 
isomorphism, for $\alpha 
\in \pi_1(X \setminus {\bf x};*)$ we will denote $(\alpha)_l$ and 
 $m_t(\alpha)$ the lifting   
$(\alpha^{\tau_{{\bf x}^t,\mathfrak{O}}})_l$ 
and permutation $m_t(\alpha^{\tau_{{\bf x}^t,\mathfrak{O}}})$, if no
ambiguity seems likely to arise. 
\end{notacion}

\subsection{Discriminant method for a Morse polynomial.}\label{s:dis_met_morse}

In this section we consider the following hypothesis.
\begin{hipotesis}\label{hip2cap2}
\hfill
\begin{enumerate}[i.]
\item\label{buenamorsif.0} $f$ is a Morse function with critical points 
$P=\{(x_i,y_i)\} \subset \Int(X) \times \Int(Y)$
 and~$\mu(f)$ critical values $f(P)=\{t_i\} \subset \Int(T)$.
\item\label{buenamorsif.2}  $\mathfrak{D}_{f}$ is reduced.
\item\label{buenamorsif.3}   If $(t_o,x_o) \in \mbox{Sing}(\mathfrak{D}_{f})$ 
then	 $t_o \notin f(P) \cup \partial T$.
\end{enumerate}
\end{hipotesis}

By Hypothesis~\ref{hip2cap2} we have $\mathcal{T} \cap \partial T= \emptyset$ 
and the set $\mathcal{T}_{Sing}$ is finite and disjoint with~$f(P)$.

\begin{definicion}\label{d:base_dis_cam}
A list of paths $(\gamma_1,\ldots,\gamma_{\mu(f)})$ is said  a
\textit{distinguished basis of $f$ based at $\star \in \partial T$}
if it is a basis of $\pi_1(T \setminus f(P);\star)$ with 
 $\supp(\gamma_i) \subset \Int(T)  \setminus \mathcal{T}$ and
there exist   homotopy  representatives $\tilde{\gamma}_i$ of $\gamma_i$ such as 
 $(\tilde{\gamma}_1, \ldots,\tilde{\gamma}_{\mu(f)})$ is a geometric basis
of $\pi_1(T \setminus f(P);\star)$. 
\end{definicion}

\begin{teorema}\label{t:met_discrim}
Let $f\in \C[x,y]$ be a monic polynomial of degree $N$ in the variable $y$ which satisfies 
Hypothesis~\ref{hip2cap2} and let $(\gamma_1, \dots, \gamma_{\mu(f)})$ be a distinguished basis 
of~$f$ based at $\star \in \partial T$. Let $\mathfrak{O}$ be a base set of braids and
$(\mu_1, \ldots, \mu_k)$ as in Example~\ref{ejemplobase} with $k$
the degree of $d_f(t,x)$ in the variable $x$. Then for any $j \in \{1,\ldots,k\}$:
\begin{enumerate}
\item\label{met_dis.2}  There exists a braid
$\beta_i \in \B_k$, depending on $j$, such 
that~$ \nabla_{\star, \mathfrak{O}}(\gamma_i)=~\beta_i^{-1} \sigma_j \beta_i$.
\item\label{met_dis.1} With Notation~\ref{n:convencion} for the isomorphism 
$\Psi_{\tau_{{\bf x}^{\star},
\mathfrak{O}}}$ we have:
\begin{enumerate}
\makeatletter\renewcommand{\p@enumii}{\theenumi.}\makeatother
\item\label{met_dis.1.a} 
The permutation $m_{\star}(\mu_{j}^{\beta_i })$ is a transposition $(a_i,b_{i}), 1\leq a_i<b_{i}\leq N$. 
\item\label{met_dis.1.b} 
If we denote 
$\delta:=\mu_{j+1} \cdot \mu_{j} \subset X \setminus {\bf x}$, then
\begin{enumerate}
\makeatletter\renewcommand{\p@enumiii}{\theenumi.\theenumii.}\makeatother
\item\label{met_dis_Delta}  $(\Delta_1, \dots, \Delta_{\mu(f)})$
 with $
\Delta_i:=(\delta^{\beta_i} )_{a_{i}}
$
defines a distinguished basis of 
vanishing 1-cycles
in $H_1(V_{\star}(f))$. Also  the cycle $\Delta_i$   is the first non homologically zero lifting 
in $((\delta^{\beta_{i}})_{1}, \ldots, (\delta^{\beta_{i}})_{N})$.

\item\label{met_dis_h(Delta)} If we denote $\beta_{s,i}:=
\beta_s  \nabla_{\star, \mathfrak{O} }(\gamma_i) \in \B_k$, the Picard-Lefschetz operators
$
(h_{\gamma_i})_*
$
satisfy
$
(h_{\gamma_i})_{*}(\Delta_s)= 
(\delta^{\beta_{s,i}})_{a_{s}}.
$ Also  the cycle $(h_{\gamma_i})_{*}(\Delta_s)$ is  the unique lifting
in $((\delta^{\beta_{s,i}})_{1},\ldots, (\delta^{\beta_{s,i}})_{N})$
that can be written in homology as $\Delta_s-n \Delta_i$ for some $n\in \Z$.
\end{enumerate}
\end{enumerate}
\end{enumerate}
\end{teorema}

We will prove Theorem~\ref{t:met_discrim} in Section~\ref{s:proof}. 
Note that the paths 
$\delta^{\beta_i}$ give us the vanishing cycles and 
$\delta^{\beta_{s,i}}=\Psi_{\nabla_{\star,\mathfrak{O} }(\gamma_i)}(\delta^{\beta_s})$ 
describes the Picard-Lefschetz operators. Hence we define:
\begin{definicion}\label{d:cam_evan}
For $i=1,\ldots,\mu(f)$  paths $\delta^{\beta_i}  $ are said \textit{vanishing paths} and isomorphisms 
$\Psi_{\nabla_{ \star,\mathfrak{O}}(\gamma_i)}$ \textit{Picard-Lefschetz transformations}.
\end{definicion}

Theorem~\ref{t:met_discrim} shows us how to compute, in an algebraic way, the 
elements $\Delta_i$ and~$(h_{\gamma_s})_{*}(\Delta_i)$. To do this we need some extra
data. To simplify notation we denote by $\gamma_0$ the identity path onto $\star$ and so 
we can extend the definition of the $\beta_{s,i}$ to 
$\beta_{s,0}:=\beta_s  \nabla_{ \star,\mathfrak{O}}(\gamma_0)$. 
The extra data are:
\begin{enumerate}[(a)]
\item\label{cons_cicl_seg}  
The list of permutations 
$(m_{\star}(\mu_1^{\tau_{{\bf x}^{\star},\mathfrak{O}}}), 
\dots,m_{\star}(\mu_k^{\tau_{{\bf x}^{\star},\mathfrak{O}}}))$
as elements of the permutation group $\Sigma_N$.

\item\label{cons_cicl_cuar} If we fix, for example,  $j=1$, a representative 
of $\nabla_{ \star,\mathfrak{O}}(\gamma_i),i=1,\ldots, \mu(f),$ 
in the braid group  $\B_k$  which is written as
$
\beta_i^{-1} \sigma_{1} \beta_i.
$
\end{enumerate}

With these data we compute $\Delta_i$ and $(h_{\gamma_i})_{*}(\Delta_s)$  
algebraically as follows:
\begin{enumerate}[1)]
\item By~\eqref{formulaexp} we compute the paths  
$\mu_1^{\beta_{s,i}}$ and  $\mu_2^{\beta_{s,i}}$ expressed in the 
geometric basis~$(\mu_1, \dots, \mu_k )$. 
\item We compute the transposition $(a_{s},b_{s})=m_{\star}(\mu_{1}^{\beta_{s,0} \tau_{{\bf x}^{\star},\mathfrak{O}} })$ 
using~\eqref{cons_cicl_seg}.
\item We lift 
$
\delta^{\beta_{s,i}}=\mu_2^{\beta_{s,i}} \cdot \mu_1^{\beta_{s,i}}
$
at the point $a_s$ as it is showed by Proposition~\ref{propcamhomolog}\eqref{propcamhomolog.a}
(using the geometric basis 
$(\mu_1^{\tau_{{\bf x}^{\star},\mathfrak{O}}}, \dots,
\mu_k^{\tau_{{\bf x}^{\star},\mathfrak{O}}})$) to obtain
$\Delta_i$ and $(h_{\gamma_s})_*(\Delta_i)$.
\end{enumerate}
Items (2) and (3) above can be merged into an unique point:
\begin{itemize}
\item[(2')] For $l=1, \dots, N,$ with Proposition~\ref{propcamhomolog}\eqref{propcamhomolog.a}
we lift $(\delta^{\beta_{s,0}})_{l}$ until we obtain a non-zero 
element in homology and we lift $(\delta^{\beta_{s,i}})_{l}$
until it can be written in homology as $\Delta_s-n \Delta_i$ for some $n\in \Z$.
\end{itemize}

\begin{observacion}\label{o:prac} 
To use Theorem~\ref{t:met_discrim} in an effective way we can observe that:
\begin{enumerate}[1)]
\item\label{o:prac_1} If $f$ verifies~\eqref{buenamorsif.0} and~\eqref{buenamorsif.2} 
in Hypothesis~\ref{hip2cap2} then one can check  the condition
$(t_o,x_o) \in \mbox{Sing}(\mathfrak{D}_{f}) \Rightarrow t_o \notin f(P),$ 
as in~\eqref{buenamorsif.3}, if either  every $t_i$ is a simple root of the polynomial
$\dis_x(d_f(t,x)) \in \C[t]$ or if every $\nabla_{ \star,\mathfrak{O}}(\gamma_i)$ 
is a conjugated of  $\sigma_{1} \in \B_k$. Finally, condition 
$(t_o,x_o) \in \mbox{Sing}(\mathfrak{D}_{f}) \Rightarrow t_o \notin \partial T,$
is true with only a small perturbation of $T$.
\item\label{o:prac_2} 
It is not necessary to find $X$ in Hypothesis~\ref{hip2cap2} explicitly. 
Let $\tilde{X} \subset \C$ be a
disk with ${\bf x}^{\star} \subset \Int(\tilde{X})$ and $\tilde{*} \in \partial 
\tilde{X}$. Let $(\tilde{\mu}_1, \ldots, \tilde{\mu}_k)$ be
a geometric basis of $\pi_1(\tilde{X} \setminus {\bf x}^{\star}; \tilde{*})$ 
such that there exists a path $r$ that begins at $*$ and ends at $\tilde{*} \in \C$  
so that
$$
\mu_j^{\tau_{{\bf x}^{\star},\mathfrak{O}}}=r \cdot \tilde{\mu}_j \cdot r^{-1}
$$
in the group $\pi_1(\C \setminus {\bf x}^{\star}; *)$. 
Hence $m_{\star}(\tilde{\mu}_j)$ and 
$m_{\star}(\mu^{\tau_{{\bf x}^{\star},\mathfrak{O}}}_j)$ 
are conjugated by a permutation and if we define the braid
$\tilde{\beta}_{s,i}:=\tau_{{\bf x}^{\star},\mathfrak{O}}^{-1} 
\beta_{s,i} \tau_{{\bf x}^{\star},\mathfrak{O}}$
and we write 
$
\mu_j^{\beta_{s,i}}=~\mu_{j_1}^{\epsilon_{1}}\cdot~\ldots
\cdot~\mu_{j_m}^{\epsilon_{m}},$ $\epsilon_{n} \in \{-1,1\},
$
then one has
\begin{equation}\label{e:con_tilde}
\begin{array}{rcl}
\vspace*{.1cm}
\mu_j^{\beta_{s,i} \tau_{{\bf x}^{\star},\mathfrak{O}}}&=&
\Psi_{\tau_{{\bf x}^{\star},\mathfrak{O}}}(\mu_j^{\beta_i})=
\Psi_{\tau_{{\bf x}^{\star},\mathfrak{O}}}(\mu_{j_1}^{\epsilon_{1}}\cdot \ldots
\cdot  \mu_{j_m}^{\epsilon_{m}})
\\
\vspace*{.1cm}
&=& r \cdot \tilde{\mu}_{j_1}^{\epsilon_{1}}\cdot \ldots
\cdot  \tilde{\mu}_{j_m}^{\epsilon_{m}}\cdot r^{-1}=
r \cdot \tilde{\mu}_{j}^{\tilde{\beta}_{s,i}}\cdot r^{-1}.
\end{array}
\end{equation}
It is easy to prove that permutation
$
m_{\star}(\tilde{\mu}_{j}^{\tilde{\beta}_{i,0}})$ is a transposition, say  
$(\tilde{a}_i,\tilde {b}_i)$, $1 \leq \tilde{a}_i < \tilde {b}_i \leq N,$
and if we define $\tilde{\delta}:=\tilde{\mu}_{j+1} \cdot \tilde{\mu}_{j}$, 
the elements $(\tilde{\delta}^{\tilde{\beta}_{s,i}})_{\tilde{a}_{i}},$
are, up to sign, $(h_{\gamma_i})_*(\Delta_s)$.
Furthermore, by~\eqref{e:con_tilde}, computing $(\tilde{\delta}^{\tilde{\beta}_{s,i}})_{\tilde{a}_{i}}$
only requires computing $\delta^{\beta_{s,i}}$ using~\eqref{formulaexp}, changing 
every $\mu_{j_n}^{\epsilon_{n}}$ by $\tilde{\mu}_{j_n}^{\epsilon_{n}}$ and, with  
the list of permutations 
$(m_{\star}(\tilde{\mu}_1), \dots, m_{\star}(\tilde{\mu}_k)) \subset \Sigma_N^k$,
computing the lifting at the point $\tilde{a_i}$ as it is showed by
Proposition~\ref{propcamhomolog}\eqref{propcamhomolog.a}.
In other words, up to sign, $(h_{\gamma_i})_*(\Delta_s), i=0,1,\ldots, \mu(f),$ 
only needs $\beta_{s,i}$ and $(m_{\star}(\tilde{\mu}_1), \dots, m_{\star}(\tilde{\mu}_k))$.

\item\label{o:prac_3} If we decompose 
$\nabla_{ \star,\mathfrak{O}}(\gamma_i)=\beta_i^{-1} \sigma_{j(i)} \beta_i$ 
for some $j(i), 1 \leq j(i) \leq k,$ depending on~$i$, by Lemma~\ref{l:reduc_1} we obtain
$\nabla_{ \star,\mathfrak{O}}(\gamma_i)=(\beta_i \alpha_{j(i)})^{-1} \sigma_{1}
\alpha_{j(i)} \beta_i$ and the equality $\mu_1^{\alpha_{j(i)} \beta_i}=\mu_{j(i)}^{\beta_i}$. Then
$(a_i,b_i)=m_{\star}(\mu_{j(i)}^{\beta_i})$ and, 
if we define $\delta_i:=\mu_{j(i)+1} \cdot \mu_{j(i)}$, 
$\delta^{\alpha_{j(i)} }=\delta_{i}$.
So $(\delta_{i}^{\beta_{i}})_{a_i}$ and $(\delta_{s}^{\beta_{s} \nabla(\gamma_i)})_{a_s}$ are $\Delta_i$ and
$(h_{\gamma_i})_*(\Delta_s)$, again.
\end{enumerate}
\end{observacion}

Finally note that, by Picard-Lefschetz Theorem, the element $n$ defined in 
Theorem~\ref{t:met_discrim}\eqref{met_dis_h(Delta)} 
is the intersection number between the 
vanishing cycles $\Delta_s$ and $\Delta_i$, i.e., $n=(\Delta_s,\Delta_i)$. 
As a corollary we have:
\begin{corolario}\label{c:inter_0}
Let $E(\sigma_{j})$ be the stabilizer of $\sigma_{j}$ by the 
action~$\Phi_{{\bf x}, {\bf x}}$ given in~\eqref{formulaexp}. If~$\delta^{ \beta_{s} \beta_i^{-1}} \in E(\sigma_{j})$ then $(\Delta_s,\Delta_i)=0$.
\end{corolario}
\begin{proof}
In homotopy we have 
$\delta^{\beta_{s,i}}=\delta^{\beta_s  \beta_i^{-1} \sigma_{j} \beta_i}=\delta^{\beta_s },$
then by Theorem~\ref{t:met_discrim} we have $(h_{\gamma_i})_{*}(\Delta_s)= \Delta_s$.
\end{proof}

\begin{observacion}
Since $\delta^{ \beta_i \beta_i^{-1}}= \delta \in E(\sigma_{j})$ 
by Corollary~\ref{c:inter_0} we have $(\Delta_i,\Delta_i)=0$ which is already known in
the case of two variables (see for example~\cite[Lemma~1.4]{AGV}).
\end{observacion}

\subsection{Discriminant method for tame polynomials.}\label{S:pol_tame}

In this section we consider the following hypothesis.
\begin{hipotesis}
The polynomial application $f$ is tame.
\end{hipotesis}

\begin{definicion}
We will call a morsification of~$f$ which satisfies Hypothesis~\ref{hip2cap2} a 
\textit{good morsification}.  
\end{definicion}

\begin{proposicion}\label{p:bu_mor_pol}
Almost every perturbation~$f_{\lambda,a}(x,y):=f(x+\lambda y,y)+ax$ of $f$, 
$\lambda, a~\in \C, \vert \lambda \vert, \vert a \vert$ sufficiently small, 
is a good morsification of~$f$.
\end{proposicion}
\begin{proof}
Since $f$ has isolated singularities, by Lemma~\ref{p:buenas_coor_abier} 
the discriminant~$\mathfrak{D}_{f_{\lambda,0}}$ is reduced in~$T \times \C$ for almost 
every $\lambda \in \C$. Let $\lambda \in \C$ be so that
$\mathfrak{D}_{f_{\lambda,0}}$ is reduced in~$T \times \C$.
The polynomial $f_{\lambda,0}$ is tame because $f$ is. Then~$f_{\lambda,0}$ 
is regular at infinity and hence the leader coefficient of~$d_{f_{\lambda,0}}(t,x)$, 
as polynomial in the variable~$x$, is a non-zero constant and therefore
$\mathfrak{D}_{f_{\lambda,0}}$ has no vertical asymptotes.

Let $\vert \lambda \vert$ be small enough such that the critical points of $f_{\lambda,0}$ 
are inside $M$. Since the discriminant curve $\mathfrak{D}_{f_{\lambda,0}}$ has a 
finite number of singularities and inflection points in $T \times X$, a generic 
linear transformation
$$\{ (t,x) \in T \times \C \mid d_{f_{\lambda,0}}(t-ax,t)=0\}, a \in \C,$$ 
with $\vert a \vert$ sufficiently small, is again reduced and has no 
vertical asymptotes and also vertical lines $\{t=t_o\}, t_o \in \C$, 
are not tangent to singular points or inflection points.
This curve is the discriminant $\mathfrak{D}_{f_{\lambda,a}}$.

Discriminant $\mathfrak{D}_{f_{\lambda,a+\tilde{a}}}$, with  
$\vert\tilde{a} \vert$ sufficiently small, has the same properties in
$T \times \C$ that $\mathfrak{D}_{f_{\lambda,a}}$ and also
every vertical line $\{t=t_o\}$ is at most tangent at a smooth (non-inflexion) point 
and if it is tangent then there is no critical point on it.
 
From the properties of $\mathfrak{D}_{f_{\lambda,a+\tilde{a}}}$, 
Corollary~\ref{c:iden_vc} and Lemma~\ref{c:Morse} one obtains that the polynomial $f_{\lambda,a+\tilde{a}}$ has no degenerated singularities in 
$f_{\lambda,a+\tilde{a}}^{-1}(T)$. Also, if  $\vert a + \tilde{a} \vert$ 
is sufficiently small, $f_{\lambda,a+\tilde{a}}$ has no other singularities, 
because $f_{\lambda}$ is tame, and they are in $\Int(X) \times \Int(Y)$.
\end{proof}

By Proposition~\ref{p:bu_mor_pol} we deduce the existence of good morsifications of $f$
and thus we have the following theorem.
\begin{teorema}\label{t:moa_pol}
Local monodromies of a tame polynomial $f(x,y)$ which is monic in the variable $y$ 
can be computed in an algebraic effective way using a good morsification $\tilde{f}$ of $f$ and Theorem~\ref{t:met_discrim} for $\tilde{f}$.
\end{teorema}
\begin{observacion}
If $f$ is good at infinity but not tame we can use Theorem~\ref{t:moa_pol} with a good morsification $\tilde{f}$ but we must distinguish between the critical values of 
$\tilde{f}$ that come from the morsification of affine critical points of $f$ and the 
ones that come from critical points at infinity of $f$.
\end{observacion}

\section{Topology of conjugated polynomials}\label{C:top_conj} 
At the present section we show two polynomial applications  $f^+$ and $f^-$, 
conjugated in a number field that are not topologically equivalent.

\subsection{A uniparametric family of polynomials.}\label{s:prop_top_f+-}

We consider the family of polynomials 
$f^s(x,y):=c(x,y) r^s(x,y), s \in \C,$
given by the product of a cubic curve
$c(x,y):= {y}^{2}x- \left( y+1 \right) ^{3}$ 
and a line
$r^s(x,y):=  {s}^{2} \left( 2\,s-3 \right) y+x-3\,{s}^{2}.$

The polynomial $f^s(x,y)$, with $s$ generic, is such that:
\begin{enumerate}
\item\label{i:1} The zero fiber is singular with singularities $\A_1$ and $\A_3$.
\item There are critical values $t^s_1, t^s_2\neq 0 $   with  only one
singularity $\A_1$ over their fibers.
\item\label{i:3} The other fibers are smooth.
\item\label{i:4} It is tame.
\end{enumerate}
Then if $s_1,s_2 \in \C$ are generic, the polynomial applications
$f^{s_i}: \C^2 \rightarrow \C, i=1,2,$
are topologically equivalent because they are in a uniparametric equisingular 
family and they have the same link at infinity.

We write the homogeneous polynomial $F^s(x,y,z)$ associated with $f^s(x,y)$. 
If we compute a Gr\"obner basis (with respect to the lexicographic 
order given by $(x,y,t)$) of the ideal
$\langle F^s_x(x,y,1),F^s_y(x,y,1),F^s(x,y,1)-t \rangle,$
the polynomial
$$
(256s^2+256s+64)t^2+(-1647s^4-1836s^3+2430s^2+2916s+729)t+2916s^6
$$
belongs to this ideal. Its discriminant with respect to the variable $t$ is
$$
-2187(s-1)^2(s-3-2\sqrt{3})^2(s-3+2\sqrt{3})^2(5s+3)^2.
$$
Then the critical values $t^s_1$ and $t^s_2$ are equal if  $s=3\pm2\sqrt{3}$. 
We write
$f^+:=~f^{3+2\sqrt{3}}$ and $f^-:=f^{3-2\sqrt{3}}.$
These polynomials satisfy again the properties~\eqref{i:1}, \eqref{i:3} and~\eqref{i:4} 
and they have only one non-zero critical value
($t^+:=t^{3+2\sqrt{3}}$ and $t^-:=t^{3-2\sqrt{3}}$ respectively), with two 
singularities $\A_1$ in the singular fiber $V_{t^{\pm}}(f^{\pm})$. 

Since the polynomials $f^+$ and  $f^-$ are conjugated in the number field $\Q(\sqrt{3})$,
they have the same algebraic properties. The reason is the existence of the Galois
isomorphism between the two of them. Also the pairs of quartics $V_0(f^{\pm})$,
$V_{t^{\pm}}(f^{\pm})$, and $V_{t}(f^{\pm}), t\neq 0,t^{\pm},$ are diffeomorphic because
they have the same combinatorial properties (see, for example, classification 
in~\cite[Section 2.2]{N}, \cite{Degt89}). In what follows we will prove that,
the polynomial applications $f^+$ and $f^-$, are not topologically equivalent.

If $f^{+}$ and $f^{-}$ were topologically equivalent then there would exist
$\Theta$ and $\theta$ homeomorphisms such that the diagram
\begin{equation}\label{conm_fpm}
\begin{CD}
\C^{2}	@>\Theta>>	\C^{2} \\
@Vf^{+}VV	 @VVf^{-}V \\
\C @>\theta>> \C
\end{CD}
\begin{picture}(0,0)
%\graphpaper(-90,-50)(100,75)
\put(-37,-2){$\circlearrowleft$}
\end{picture}
\end{equation}
is commutative. Then $\Theta$ sends singular fibers to singular fibers and $\theta$ 
critical values to critical values. Since both singular fibers are not homeomorphic, 
$\theta$ must send $0$ to $0$ and  $t^+$ to $t^-$ respectively. Finally note that, 
using the complex conjugation $\ \bar{ }:z \mapsto \bar{z}$ if necessary (changing
$\Theta$ by $\bar \Theta$, and $\theta$ by $\bar \theta$), we can assume that 
$\theta$ preserves the orientation.

Let $T$ be a geometric disk containing $0$ and $t^{\pm}$, also let $\star^{\pm} \in \partial T$ 
be regular values of $f^{\pm}$ and $r:[0,1] \rightarrow \C$ a path with support
in $T \setminus \{0,t^-\}$  such that $r(0)=\theta(\star^+)$ and $r(1)=\star^-$. 
The pair of homeomorphism $\theta$ and path $r$ induces the group isomorphism
\begin{equation*}
\begin{matrix}
\theta_*:&\pi_1(T \setminus \{0,t^{+}\};\star^+) &\longrightarrow& 
\pi_1(T \setminus \{0,t^-\};\star^-)\\
 & \gamma& \mapsto & r^{-1} \cdot \theta(\gamma) \cdot r.
 \end{matrix}
\end{equation*}
Let $h_r$ be a diffeomorphism given by lifting the path $r$. The pair of 
homeomorphism $\Theta$ and diffeomorphism $h_r$ induces the isomorphism
\begin{equation*}
\begin{matrix}
\Theta_*:&H_1(V_{\star^+}(f^{+})) &\longrightarrow  &H_1(V_{\star-}(f^{-})) \\
& \Delta & \mapsto & (h_r)_*(\Theta(\Delta))
\end{matrix}
\end{equation*}
of free Abelian groups of rank 
$\mu(f^{\pm})=\mu(\A_1)+\mu(\A_3) +\mu(\A_1)+\mu(\A_1)=6.$

Let $\mathcal{B}^{\pm}$ be a basis in $H_1(V_{\star^{\pm}}(f^{\pm}))$,
then $\Theta_*$ is given by an invertible matrix $P_{\Theta} \in \mbox{GL}(6, \Z)$ 
which induce the group isomorphism
\begin{equation*}
\begin{matrix}
\mbox{Aut}H_1(V_{\star^+}(f^{+}))  &\stackrel{\mathcal{B}^+}{\cong}&  
\mbox{GL}(6, \Z) &\stackrel{\Theta^*}{\longrightarrow}
&\mbox{GL}(6, \Z) &\stackrel{\mathcal{B}^-}{\cong} 
&\mbox{Aut}H_1(V_{\star^-}(f^{-})) \\ & & P & \mapsto & P_{\Theta}^{-1} P 
P_{\Theta} & &\end{matrix}
\end{equation*}
The commutativity of the diagram~\eqref{conm_fpm} and the isomorphisms 
of groups $\theta_*$ and $\Theta^*$ induce the commutative diagram given by the 
homologic monodromies of $f^{\pm}$, $\rho_*^{\pm}:\gamma \rightarrow (h_{\gamma})_*$:
\begin{equation}\label{isofinf}
\begin{CD}
\pi_1(T \setminus \{0,t^{+}\};\star^+) @>\rho_{*}^{+}>> 
\mbox{Aut}H_1(V_{\star^+}(f^{+}))  \stackrel{\mathcal{B}^+}{\cong}  
\mbox{GL}(6, \Z) \\
@V\theta_{*}VV  \hspace*{3cm} @VV\Theta^{*}V   \\
\pi_1(\C \setminus \{0,t^-\};\star^-) @>\rho_*^->> 
\mbox{Aut}H_1(V_{\star^-}(f^{-})) \stackrel{\mathcal{B}^-}{\cong}  \mbox{GL}(6, 
\Z)\\\end{CD}
\begin{picture}(0,0)
%\graphpaper(-90,-50)(100,75)
\put(-155,-2){$\circlearrowleft$}
\end{picture}
\end{equation}

Let $(\gamma_1^{\pm},\gamma_0^{\pm})$ be geometric bases of 
$\pi_1(T \setminus \{0,t^{\pm}\};\star^{\pm})$
where $\gamma_0^{\pm}$ are meridians of~$0$ and $\gamma_1^{\pm}$ are meridians of~$t^{\pm}$. 
Since $\theta(0)=0$ and $\theta(t^+)=t^-$, $\theta_*(\gamma_0^+)$ is
conjugated to~$\gamma_0^-$, $\theta_*(\gamma_1^+)$ is  conjugated to~$\gamma_1^-$, and 
$\theta_*(\gamma_0^+ \cdot \gamma_1^+)$ is conjugated to $\gamma_0^- \cdot \gamma_1^-$.
Therefore there exists $\gamma \in \pi_1(T \setminus \{0,t^{+}\};\star)$ such that
$$
\theta_{*}(\gamma) \cdot \theta_{*}(\gamma_0^+) \cdot \theta_{*}(\gamma_1^+) 
\cdot  \theta_{*}^{-1}(\gamma)=\gamma_0^- \cdot \gamma_1^-.
$$
Thus if we consider
\begin{equation}\label{e:cam_base_gamma+}
\tilde{\gamma}_i^+:= \gamma \cdot \gamma_i^+ \cdot \gamma^{-1},
\end{equation}
the new basis $(\tilde{\gamma}_1^+,\tilde{\gamma}_0^+)$ is such that
$\theta_*(\tilde{\gamma}_0^+)$ is conjugated to $\gamma_0^-$,
$\theta_*(\tilde{\gamma}_1^+)$ is conjugated to~$\gamma_1^-$, and
$\theta_*(\tilde{\gamma}_0^+ \cdot \tilde{\gamma}_1^+)= 
\gamma_0^- \cdot \gamma_1^-$. Therefore
$(\theta_*(\tilde{\gamma}_1^+), \theta_*(\tilde{\gamma}_0^+))$ is a geometric 
basis of the group $\pi_1(T \setminus \{0, t^-\}; \star^-)$.

By Artin Theorem, there exists a Hurwitz move that sends the geometric basis 
$(\theta_*(\tilde{\gamma}_1^+), \theta_*(\tilde{\gamma}_0^+)) $ to the   
geometric basis $(\gamma_1^-, \gamma_0^-)$. This move must have 
$\sigma_1^{2n} \in \B_{2}$, $n \in \Z$, as its associated braid since 
$\theta_*(\tilde{\gamma}_i^+)$ and $\gamma_i^-$ are conjugated, that is, 
since fiber order must be respected. Hence
\begin{equation}\label{e:mov_sigma_2n}
\Psi_{\sigma_1^{2n}}(\theta_*(\tilde{\gamma}_i^+))=\gamma_i^-,
\end{equation}
where $\Psi_{\sigma_1^{2n}}$ is the Hurwitz move associated with $\sigma_1^{2n}$.

\begin{observacion}
The fact that $f^{\pm}$ has only two singular non-diffeomorphic fibers makes the 
study of the commutative diagram~\eqref{conm_fpm} easier. For example, it was
key in order to deduce that the Hurwitz move had to be of type $\sigma_1^{2n}$.
\end{observacion}

\begin{notacion}\label{n:H+-}
We denote by $H_0^{\pm}$ and $H_1^{\pm}$ the monodromy matrices associated 
with the operators $(h_{\gamma_0^{\pm}})_*$ and $(h_{\gamma_1^{\pm}})_*$ in 
the bases $\mathcal{B}^{\pm}$, respectively.
\end{notacion}

By~\eqref{e:cam_base_gamma+} we have
$(h_{\tilde{\gamma}_i^+})_*=(h_{\gamma})_*^{-1} 
\circ (h_{\gamma_i^+})_* \circ (h_{\gamma})_*.$
Thus, if $P_{\gamma} \in \mbox{GL}(6,\Z)$ is the matrix of $(h_{\gamma})_*$ 
in the basis $\mathcal{B}^+$, then the matrices 
$P_{\gamma}^{-1} H_0^{+} P_{\gamma} \ \, \mbox{and}  \, 
\ P_{\gamma}^{-1} H_1^{+} P_{\gamma}$
define the monodromy operators $(h_{\tilde{\gamma}_0^+})_*$ and 
$(h_{\tilde{\gamma}_1^+})_*$ in the basis $\mathcal{B}^+$. 
Then, using the commutative diagram~\eqref{isofinf}, the matrices 
$P_{\Theta}^{-1} P_{\gamma}^{-1} H_0^{+} P_{\gamma} P_{\Theta}$ and 
$P_{\Theta}^{-1} P_{\gamma}^{-1} H_1^{+} P_{\gamma} P_{\Theta}$
define the monodromy operators
$(h_{\theta_*(\tilde{\gamma}_0^+)})_*$ and 
$(h_{\theta_*(\tilde{\gamma}_1^+)})_*$ in the basis $\mathcal{B}^-$.

We define 
$$
\alpha_1:=\gamma_0^- \, \cdot \, \gamma_1^-, 
\alpha_n:=\alpha_{n-1}^{\sigma_1^2} \, \cdot \, \alpha_1 \ \mbox{if} \ n>1. 
$$	
The isomorphism $ \Psi_{\sigma_1^{2n}} $ of 
$\pi_1(T \setminus \{0, t^-\}; \star^-)$ is such that
\begin{equation*}
\Psi_{\sigma_1^{2n}}(\varrho)=
\begin{cases}
\alpha_1 \, \cdot \, \varrho \, \cdot \, \alpha_1^{-1} &  \mbox{if} \ n=1 \, \\
\alpha_n \, \cdot \, \varrho \, \cdot \, \alpha_n^{-1} &  \mbox{if} \ n>1. 
\end{cases}
\end{equation*}
Therefore it is an inner automorphism of 
$\pi_1(T \setminus \{0,t^{-}\};\star^-)$. By~\eqref{e:mov_sigma_2n},
$\gamma_i^-,i=0,1,$ are conjugated to $\theta_*(\tilde{\gamma}_i^+),i=0,1,$
by the same element $\alpha_n \in \pi_1(T \setminus \{0,t^{-}\};\star^-) $. 
If we denote by $P_{\alpha_n} \in GL(6, \Z)$ the associated matrix 
to $(h_{\alpha_n})_*$ in the basis $\mathcal{B}^-$ and define
$P:=P_{\gamma} P_{\Theta}  P_{\alpha_n}$, then the matrices 
$P^{-1} H_0^{+} P$ and $ P^{-1} H_1^{+} P$ must define the monodromy operators  
$(h_{\gamma_0^-})_*$ and $(h_{\gamma_1^-})_*$ in the basis~$\mathcal{B}^-$.

By Notation~\ref{n:H+-}, if $f^+$ and $f^-$ were
topologically equivalent, then they should satisfy
$H_0^{-}=  P^{-1} H_0^{+} P $ and $H_1^{-}= P^{-1} H_1^{+} P$
where $P \in GL(6, \Z)$, or analogously
\begin{equation}\label{igualmatriz_0}
P \, H_0^{-} =    H_0^{+} \, P, \ 
P \,  H_1^{-}  =    H_1^{+} \, P \ \, \mbox{and} \ \mbox{det}(P)\neq 0. 
\end{equation}

\subsection{Calculus of the homological monodromy.}\label{s:grup_mon_f+-}
The polynomials $f^{\pm}$ are monic, 
of degree 4 in the variable $y$, and tame. In order to compute their homological 
monodromy we can use Theorem~\ref{t:moa_pol}.

Since $d_{f^{\pm}}(t,x)$ is reduced (since $\dis_x(d_{f^{\pm}}(t,x))\neq 0$) 
it is enough to consider $f^{\pm}_{0,a}(x,y):=f^{\pm}(x,y)+ax, a \in \C,$
a morsification of $f^{\pm}$ verifying Hypothesis~\ref{hip2cap2} for 
appropriate disks $T, X$ and $Y$. We denote such morsifications by $f_a^{\pm}$.

\begin{observacion}
From now on we will use in this section the lexicographic order $<^{\pi/2}$ (see Notation~\ref{not:orden}), but for
the sake of simplicity we will omit any reference to it in the text. 
Also as coefficients are in $\Q(\sqrt{3})$, we will use a computer program based on a method showed 
by J.~Carmona~\cite{Car03} in order to compute both the braid monodromy associated 
with the discriminant of $f_a^{\pm}$ and the classical monodromy $m_{\star^{\pm}}$ 
of $\pi \vert_{V_{\star^{\pm}}(f_a^{\pm})}$.
\end{observacion}

\begin{figure}[ht]
%\begin{center}
\vspace*{-1.5cm}
\hspace*{-3cm}
\includegraphics[scale=.8]{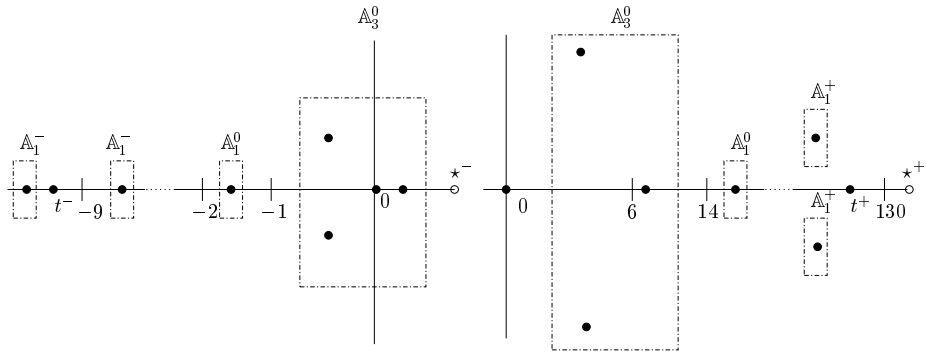}
\caption{Critical values of $f^{\pm}$ and $f_a^{\pm}$.}
%\end{center}
\label{f:val_crit_fa^+}
\end{figure}

\paragraph{\textbf{The discriminant method for $f^{\pm}$.}}
Let $a:=10^{-1}$ be a sufficiently small constant.
For this value $f_a^{\pm}$ is a Morse function with a reduced discriminant and 
without any vertical asymptotes (since $f^{\pm}_a$ is also tame).
The critical values of $f^{\pm}$ and $f_a^{\pm}$ are shown in Figure~\ref{f:val_crit_fa^+}
(note that the critical values associated with the morsification of the singularities 
$\A_3^0$ and $\A_1^0$ over the zero fiber and the singularities~$\A_1^{\pm}$ over
the~$t^{\pm}$ fiber are enclosed each in a rectangle).

\begin{figure}[ht]
%\begin{center}
\hspace*{-3cm}
\includegraphics[scale=.7]{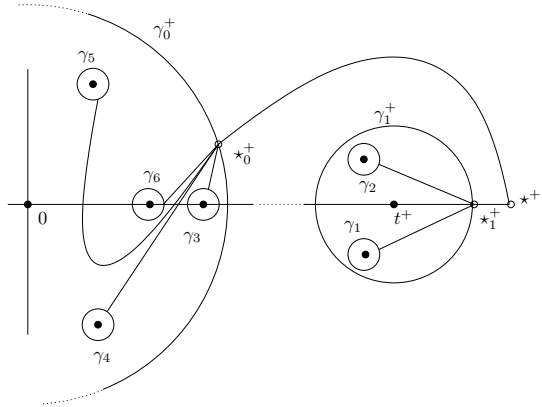}
\caption{Distinguished basis for $f^+_a$ with~$a=~10^{-1}$ and $\star^+=132$.}
%\end{center}
\label{f:base_f+}
\end{figure}

\begin{figure}[ht]
%\begin{center}
\vspace*{-1cm}
\hspace*{-3cm}
\includegraphics[scale=.7]{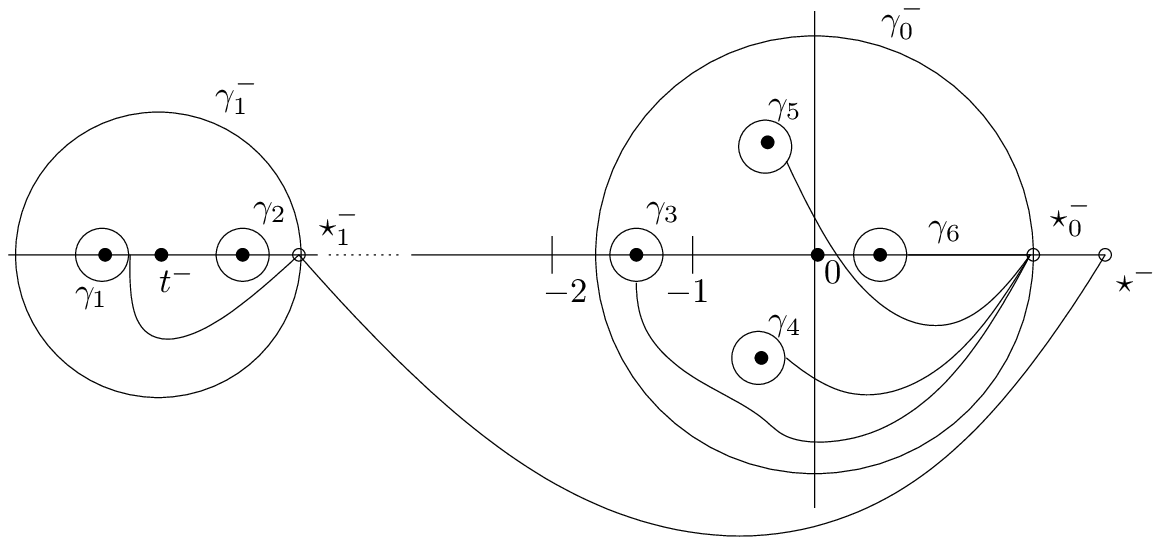}
\caption{Distinguished basis for $f^-_a$ with~$a=~10^{-1}$ and $\star^-=2$.}
%\end{center}
\label{f:base_f-}
\end{figure}

The points $\star^+:=132$ and $\star^-:=2$ are regular values of $f^{\pm}$ and $f_a^{\pm}$. 
We consider the paths given in Figures~\ref{f:base_f+} and~\ref{f:base_f-}.

The list $(\gamma_1, \dots, \gamma_6)$, in each figure, is a distinguished basis 
of $f^{\pm}_a$ based at $\star^{\pm}$ whereas the list $(\gamma_1^{\pm},\gamma_0^{\pm})$, 
is a geometric basis of $\pi_1(T \setminus \{0, t^{\pm}\};\star^{\pm})$ as considered 
in Section~\ref{s:prop_top_f+-}. 

We denote $b*a:=a^{-1}ba$ in $\B_9$. The braid monodromy for $\mathfrak{D}_{f_{a}^+}$ is:
\begin{equation}\label{e:mon_tren_f+}
\begin{array}{rl}
\nabla_{\star^+}(\gamma_1)=&    \sigma_{2}^{} * \sigma_{1}^{-1},\\
\nabla_{\star^+}(\gamma_2)=&   \sigma_{2}^{} * \sigma_{3}^{-1},\\
\nabla_{\star^+}(\gamma_3)=& \sigma_{8}^{} * \sigma_{7}^{} \sigma_{6}^{-1} 
\sigma_{5}^{-1} \sigma_{4}^{-1} \sigma_{3}^{-1} \sigma_{2}^{-1} \sigma_{3}^{-1} 
\sigma_{7}^{-1},\\\nabla_{\star^+}(\gamma_4)=& \sigma_{5} * \sigma_{6}^{} 
\sigma_{4}^{-1} \sigma_{3}^{-1} \sigma_{4}^{} \sigma_{7}^{} \sigma_{6}^{-1} 
\sigma_{2}^{-1} \sigma_{3}^{} \sigma_{1}^{-1} \sigma_{2}^{} \sigma_{5}^{-1} 
\sigma_{4}^{-1}   \sigma_{3}^{-1} \sigma_{2}^{-1} \sigma_{3}^{-1} 
\sigma_{7}^{-1},\\\nabla_{\star^+}(\gamma_5)=&   \sigma_{5} * \sigma_{6}^{} 
\sigma_{4}^{} \sigma_{7}^{} \sigma_{6}^{-1} \sigma_{3}^{} \sigma_{2}^{} 
\sigma_{5}^{-1} \sigma_{4}^{-1} \sigma_{3}^{-1} \sigma_{2}^{-1} \sigma_{3}^{-1} 
\sigma_{7}^{-1} \ \mbox{and}\\\nabla_{\star^+}(\gamma_6)=& \sigma_{6} * 
\sigma_{7}^{} \sigma_{6}^{-1} \sigma_{5}^{-1} \sigma_{4}^{-1} \sigma_{3}^{-1} 
\sigma_{2}^{-1} \sigma_{3}^{-1} \sigma_{7}^{-1}.\\\end{array}
\end{equation}
The braid monodromy for $\mathfrak{D}_{f_{a}^-}$ is:
\begin{equation}\label{e:mon_tren_f-}
\begin{array}{rl}
\nabla_{\star^-}(\gamma_1)=&  \sigma_{2}^{} * \sigma_{3}^{} \sigma_{5}^{} 
\sigma_{4}^{} \sigma_{6}^{} \sigma_{3}^{},\\\nabla_{\star^-}(\gamma_2)=& 
\sigma_8 * \sigma_{7}^{-1} \sigma_{8}^{} \sigma_{7}^{} \sigma_{6}^{} 
\sigma_{7}^{-1} \sigma_{5}^{} \sigma_{4}^{} \sigma_{5}^{} \sigma_6 
\sigma_3,\\\nabla_{\star^-}(\gamma_3)=&  
\sigma_{1}^{},\\\nabla_{\star^-}(\gamma_4)=& \sigma_{4}^{} * \sigma_{5}^{} 
\sigma_{6}^{} \sigma_{3}^{},\\\nabla_{\star^-}(\gamma_5)=& \sigma_{4}  * 
\sigma_{5}^{-1} \sigma_{6}^{-1} \sigma_{5}^{-1} \sigma_{3}^{2}  \ 
\mbox{and}\\
\nabla_{\star^+}(\gamma_6)=& \sigma_{7}^{}. \end{array}
\end{equation}
Then we obtain the elements $\beta_i, \sigma_{j(i)} \in \B_9$ of
Remark~\ref{o:prac}\eqref{o:prac_3} and we deduce,
by Remark~\ref{o:prac}\eqref{o:prac_1}, that $f_a^{\pm}$ is a good 
morsification for $X$ and $Y$ geometric disks sufficiently big.

The ordered set of roots of $d_{f_a^{\pm}}(\star^{\pm},x)=0$
is the list  ${\bf x}^{\star^{\pm}}=(x_1^{\pm}, \dots, x_9^{\pm})$ with
\begin{itemize}
\item $x_9^+ >0$ and $\vert x_9^+ \vert > \vert x_i^+ \vert,i=1,\dots,8$.
\item $x_1^-<0$ and  $\vert x_1^- \vert > \vert x_i^- \vert,i=2,\dots,9$.
\end{itemize}
We want to use Remark~\ref{o:prac}\eqref{o:prac_2}. Let us consider two points $\tilde{*}^+:=x_9+1$ and $\tilde{*}^-:=-x_1+1$
and two geometric disks $\tilde{X}^{\pm}$ with center $0\in \C$ and radius
$\tilde{*}^{\pm}$. Let us consider  Example~\ref{ejemplobase} 
with $X$ the big disk fixed above and $r^{\pm}:=[*,\tilde{*}^{\pm}]$.
By Remark~\ref{o:prac}\eqref{o:prac_2} it is enough to know the classical 
monodromy $m_{\star^{\pm}}$ of the geometric basis  
$(\tilde{\mu}_1^{\pm}, \dots, \tilde{\mu}_9^{\pm})$ 
of the group $\pi_1(\tilde{X} \setminus {\bf x}^{\star^{\pm}};\tilde{*}^{\pm})$ where
$\tilde{\mu}_j^{\pm}=(r^{\pm})^{-1} \cdot \mu_j^{\tau_{{\bf x}^{\star}}} 
\cdot r^{\pm}$.

The classical monodromy $m_{\star^+}$ of the geometric basis  
$(\tilde{\mu}_1^+, \dots, \tilde{\mu}_9^+)$ is the list
\begin{equation}\label{e:mon_clas_f+}
((1,3),(2,3),(1,3),(1,2),(1,4),(1,4),(1,3),(2,3),(1,2)),
\end{equation}
and the classical monodromy $m_{\star^-}$ of the geometric basis  
$(\tilde{\mu}_1^-, \dots, \tilde{\mu}_9^-)$ 
is the list
\begin{equation}\label{e:mon_clas_f-}
((3, 4), (3, 4), (1, 2), (2, 4), (3, 4), (2, 4), (1, 4), (1, 4), (2, 3)).
\end{equation}

Combining the data~\eqref{e:mon_tren_f+},~\eqref{e:mon_clas_f+} and 
~\eqref{e:mon_tren_f-},~\eqref{e:mon_clas_f-} we can compute the invariants
associated with the homological monodromy of $f^{\pm}$. In particular, let
$\mathcal{B}^{\pm}$ the distinguished 
basis of vanishing cycles obtained by lifting the 
distinguished bases~$(\gamma_1, \dots, \gamma_6)$ in 
Figures~\ref{f:base_f+} and~\ref{f:base_f-}, then 
the intersection matrix $\mathbb{I}_{f^{\pm}}$ with respect~to~$\mathcal{B}^{\pm}$~ is:
\begin{equation}\label{matriz1}
\mathbb{I}_{f^+}=
%\begin{center}
\begin{scriptsize}\left[ 
\begin{tabular}{l  l | l  l l l}
{ }0&{ }0&-1&{ }1&{ }1&{ }1\\
{ }0&{ }0&-1&{ }0&{ }0&{ }1\\
\hline
{ }1&{ }1&{ }0&{ }0&{ }0&{ }0\\
-1&{ }0&{ }0&{ }0&{ }1&{ }1\\
-1&{ }0&{ }0&-1&{ }0&{ }1\\
-1&-1&{ }0&-1&-1&{ }0
\end{tabular}
\right] \end{scriptsize} \ \ \mbox{and} \ \ \ 
\mathbb{I}_{f^-}=\begin{scriptsize}
\left[ 
\begin{tabular}{l l | l l l l}
{ }0&{ }0&-1&{ }1&{ }0&{ }0\\
{ }0&{ }0&{ }0&-1&-2&-1\\
\hline
{ }1&{ }0&{ }0&{ }0&{ }0&{ }0\\
-1&{ }1&{ }0&{ }0&{ }1&{ }1\\
{ }0&{ }2&{ }0&-1&{ }0&{ }1\\
{ }0&{ }1&{ }0&-1&-1&{ }0
 \end{tabular}\right]\end{scriptsize}, 
\end{equation}
where number on the $i$-th column and the $j$-th arrow is $(\Delta_i,\Delta_j)$.

\begin{ejemplo}
We will show in detail how to obtain $(\Delta_1,\Delta_3)=1$ from matrix 
$\mathbb{I}_{f^+}$. The rest of intersection numbers follow the same pattern.

By Remark~\ref{o:prac}\eqref{o:prac_3} we compute 
$\mu_2^{\beta_1}=\mu_1 \  \mbox{and} \ \mu_8^{\beta_3}=\mu_7\cdot \mu_2 
\cdot \mu_7^{-1}$. Transpositions $(\tilde{a}_1,\tilde{b}_1)$ and 
$(\tilde{a}_3,\tilde{b}_3)$, which are defined in 
Remark~\ref{o:prac}\eqref{o:prac_2}, are
$$
m_{\star^+}(\tilde{\mu}_1)=(1,3
) \ \mbox{and} \
m_{\star^+}(\tilde{\mu}_7\cdot 
\tilde{\mu}_2 \cdot \tilde{\mu}_7^{-1})=(1,2).
$$
We define $\tilde{a}_i$ to be the lowest in $\{\tilde{a}_i,\tilde{b}_i\}, i=1,3$. 
This choice only produces a change of orientation in the vanishing cycles. 
Thus we have $\tilde{a}_1= \tilde{a}_3=1$. 

By Remark~\ref{o:prac}\eqref{o:prac_3} again, vanishing paths are
\begin{equation*}
\begin{array}{rcccl}
\delta_1^{\beta_1}&=&\mu_2^{\beta_1} \cdot \mu_3^{\beta_1}
&=&\mu_3 \cdot \mu_1\\
\delta_3^{\beta_3}&=&\mu_8^{\beta_3} \cdot \mu_9^{\beta_3}
&=&\mu_9 \cdot \mu_7 \cdot \mu_2 \cdot \mu_7^{-1}\\
\end{array}
\end{equation*}
Therefore, by Remark~\ref{o:prac}\eqref{o:prac_2} and 
Proposition~\ref{propcamhomolog}\eqref{propcamhomolog.a}, we have
\begin{equation*}
\begin{array}{rclll}
\Delta_1&=&(\tilde{\mu}_{3})_{1}+ 
(\tilde{\mu}_{1})_{3}\\
\Delta_3&=&(\tilde{\mu}_{9})_{1} + (\tilde{\mu}_{7})_{2} + (\tilde{\mu}_{2})_{2} + 
(\tilde{\mu}_7^{-1})_3&=&(\tilde{\mu}_{9})_{1} 
+(\tilde{\mu}_{2})_{2}+(\tilde{\mu}_{7})_{3}.\\
\end{array}
\end{equation*}
The Picard-Lefschetz transformation $\Psi_{\nabla(\gamma_1)}$ of the vanishing 
path $\delta_3^{\beta_3}$ is
$$
 \Psi_{\nabla(\gamma_1)}(\delta_3^{\beta_3})=\delta_3^{\beta_3 
\nabla(\gamma_1)}=\mu_9 \cdot \mu_7 \cdot \mu_3 \cdot \mu_1^{-1} \cdot \mu_2  
\cdot \mu_1 \cdot \mu_3^{-1} \cdot \mu_7^{-1}.
$$
Finally, by Remark~\ref{o:prac}\eqref{o:prac_2} and
Proposition~\ref{propcamhomolog}\eqref{propcamhomolog.a}, we have
\begin{equation*}
\begin{array}{rcl}
h_{\gamma_1}(\Delta_3)&=& 
(\tilde{\mu}_{9})_{1} +
(\tilde{\mu}_{7})_{2} + 
(\tilde{\mu}_{3})_{2} + 
(\tilde{\mu}_1^{-1})_2 + 
(\tilde{\mu}_{2})_{2}  + 
(\tilde{\mu}_{1})_{3} + 
(\tilde{\mu}_3^{-1})_1 + 
(\tilde{\mu}_7^{-1})_3\\
&=&
(\tilde{\mu}_{9})_{1} + 
(\tilde{\mu}_{2})_{2} +
(\tilde{\mu}_{1})_{3} + 
(\tilde{\mu}_{3})_{1} + 
(\tilde{\mu}_{7})_{3}\\
&=&
[ (\tilde{\mu}_{9})_{1} + 
(\tilde{\mu}_{2})_{2} - 
(\tilde{\mu}_{7})_{1}]  + 
[ -(\tilde{\mu}_{1})_{1} + 
(\tilde{\mu}_{3})_{1}] = \Delta_3+\Delta_1,
\end{array}
\end{equation*}
and hence $(\Delta_1,\Delta_3)=1$.
\hfill
$\blacksquare$
\end{ejemplo}
\begin{observacion}
Note that the matrices in~\eqref{matriz1} have been decomposed into four blocks. 
The diagonal ones can be divided as 
\begin{equation*}
%\begin{center}
\begin{scriptsize}\left[ 
\begin{tabular}{l  | l}
0&0\\
\hline
0&0\\
\end{tabular}
\right] \end{scriptsize} \  \mbox{and} \
\begin{scriptsize}
\left[ 
\begin{tabular}{ l | l l l}
{ }0&{ }0&{ }0&{ }0\\
\hline
{ }0&{ }0&{ }1&{ }1\\
{ }0&-1&{ }0&{ }1\\
{ }0&-1&-1&{ }0
\end{tabular}
\right]. 
\end{scriptsize}
\end{equation*}
representing the intersection matrices of the germ $f^{\pm}$ in the singular 
points with singularities $\A_1^{\pm},\A_1^0$ and $\A_3^0$ respectively.
\end{observacion}

We have (see Figures~\ref{f:base_f+} and~\ref{f:base_f-})
$
\gamma_0^{\pm}= \gamma_6 \, \cdot \, \gamma_5 \, \cdot \, \gamma_4 \, \cdot \, 
\gamma_3 \ \mbox{and} \ \gamma_1^{\pm}=  \gamma_2 \, \cdot \, \gamma_1
$
then
\begin{equation}\label{e:op_gordos_f+}
\begin{matrix}
(h_{\gamma_0^{\pm}})_*&=& (h_{\gamma_3})_* \circ (h_{\gamma_4})_* \circ 
(h_{\gamma_5})_* \circ (h_{\gamma_6})_* \ \mbox{and}\\(h_{\gamma_1^{\pm}})_*&=& 
\hspace{-2.6cm} (h_{\gamma_1})_* \circ (h_{\gamma_2})_*.\end{matrix}
\end{equation}

By Picard-Lefschetz Theorem, the intersection matrices $\mathbb{I}_{f^{\pm}}$ give 
us the Picard-Lefschetz operators $(h_{\gamma_i})_*, i=1, \dots, 6,$ in the bases 
$\mathcal{B}^{\pm}$. Finally, by~\eqref{e:op_gordos_f+} and Notation~\ref{n:H+-}, 
we have \begin{equation}
\label{e:H0+} 
H_0^+=\begin{scriptsize}
\left[ 
\begin{tabular}{ll|llll} { }1&{ }0&{ }0&{ }0&{ }0&{ }0\\
{ }0&{ }1&{ }0&{ }0&{ }0&{ }0\\
\hline
-1&-1&{ }1&{ }0&{ }0&{ }0\\
{ }0&{ }0&{ }0&{ }0&-1&{ }0\\
{ }0&-1&{ }0&{ }0&{ }0&-1\\
{ }1&{ }1&{ }0&{ }1&{ }1&{ }1
\end{tabular} 
\right] \end{scriptsize} \ \mbox{and} \
H_1^+=\begin{scriptsize}
\left[ \begin{tabular}{ll|llll} { }1&{ }0&{ }1&-1&-1&-1\\
{ }0&{ }1&{ }1&{ }0&{ }0&-1\\
\hline
{ }0&{ }0&{ }1&{ }0&{ }0&{ }0\\
{ }0&{ }0&{ }0&{ }1&{ }0&{ }0\\
{ }0&{ }0&{ }0&{ }0&{ }1&{ }0\\
{ }0&{ }0&{ }0&{ }0&{ }0&{ }1
\end {tabular} \right] \end{scriptsize}
\begin{picture}(0,0)
%\graphpaper(-200,-50)(200,100)
\put(-236,3){\line(90,0){60}}
\put(-208,8){\line(0,-90){30}}
\put(-106,19){\line(90,0){27}}
\put(-91,26){\line(0,-90){15}}
\end{picture}
\end{equation}
in the basis $\mathcal{B}^+$ and
\begin{equation}
\label{e:H0-} 
H_0^-=\begin{scriptsize}
\left[ 
\begin{tabular}{ll|llll} 
{ }1&{ }0&{ }0&{ }0&{ }0&{ }0\\
{ }0&{ }1&{ }0&{ }0&{ }0&{ }0\\
\hline
-1&{ }0&{ }1&{ }0&{ }0&{ }0\\
{ }1&{ }1&{ }0&{ }0&-1&{ }0\\
{ }0&-1&{ }0&{ }0&{ }0&-1\\
{ }0&-1&{ }0&{ }1&{ }1&{ }1
\end {tabular} 
\right] \end{scriptsize} \ \mbox{and} \
H_1^-=\begin{scriptsize}
\left[ 
\begin{tabular}{ll|llll} 
1&{ }0&{ }1&-1&{ }0&{ }0\\
0&{ }1&{ }0&{ }1&{ }2&{ }1\\
\hline
0&{ }0&{ }1&{ }0&{ }0&{ }0\\
0&{ }0&{ }0&{ }1&{ }0&{ }0\\
0&{ }0&{ }0&{ }0&{ }1&{ }0\\
0&{ }0&{ }0&{ }0&{ }0&{ }1
\end {tabular} 
\right]\end{scriptsize}
\begin{picture}(0,0)
%\graphpaper(-200,-50)(100,75)
\put(-235,3){\line(90,0){60}}
\put(-205,8){\line(0,-90){30}}
\put(-106,19){\line(90,0){27}}
\put(-91,26){\line(0,-90){15}}
\end{picture}
\end{equation}
in the basis $\mathcal{B}^-$.

\begin{observacion}\label{o:submatriz}
Let $f:\C^n\rightarrow \C$ be a polynomial function with isolated singularities 
(no further hypotheses). As we said in Section~\ref{s:intro}, given a geometric basis
$(\gamma_1, \ldots,\gamma_{\# B_f})$ of $\pi_1(\C \setminus B_f; \star)$,
we obtain a decomposition of $H_j(V_{\star}(f))$. 
To do this let $F_{i,\alpha} \subset V_{\star}(f)$ 
($1\leq i \leq \# B_f $) be Milnor fibers of affine singularities for 
$1\leq \alpha \leq \alpha_i$ and Milnor fibers at infinity  for $\alpha_i+1\leq 
\alpha \leq \alpha_i'$. Then for  $j>0$:
\begin{equation}\label{e:desc_homo}
H_j(V_{\star}(f))\cong 
\bigoplus_{i=1}^{\# B_f} \left( \bigoplus_{\alpha=1}^{\alpha_i'} 
(v_i)_j \left( H_j \left(F_{i,\alpha}, \partial F_{i,\alpha} \right)  \right)  
\right) 
\end{equation}
where $(v_i)_j$ is the composition of the variation and inclusion operators
\begin{equation*}
\begin{matrix} (v_i)_j :& H_j(F_{i,\alpha}, \partial 
F_{i,\alpha}) & \stackrel{\text{var}}{\longrightarrow} & H_j(F_{i,\alpha})& 
\stackrel{i}{\hookrightarrow} & H_j(V_{\star}(f))\\ 
&a& \mapsto & 
a-(h_{\gamma_i})_*(a) & \mapsto & a-(h_{\gamma_i})_*(a).
\end{matrix}
\end{equation*}
Using this decomposition $(h_{\gamma_i})_*$ can be written as
\begin{equation*}
\begin{scriptsize}
\left[ 
\begin{array}{cccccccc} I_{1}&0& \cdots&0&0&0& \cdots&0\\
0&I_{2}& \cdots&0&0&0& \cdots&0
\\ \noalign{\medskip} \vdots& \vdots& & \vdots& \vdots& \vdots
& & \vdots\\\noalign{\medskip}*_{i,1}&*_{i,2}& 
\cdots&*_{i,i-1}&h_{{i}}&*_{i,i+1}&\cdots&*_{i,\# B_f}\\\noalign{\medskip} \vdots&\vdots& & \vdots& \vdots& \vdots& & 
\vdots\\\noalign{\medskip}0&0& \cdots&0&0&0& \cdots&I_{\# B_f }
\end{array} 
\right], \ 
h_i=\left[ 
\begin{array}{cccc} 
h_{{i,1}}&0&\cdots&0\\\noalign{\medskip}0&h_{{i,2}}&\cdots&0\\\noalign{\medskip}
\vdots&\vdots& &\vdots\\\noalign{\medskip}0&0&\cdots&h_{{i,\alpha_i'}}
\end{array} 
\right],
\end{scriptsize}
\end{equation*}
where $I_{j}$ is the identity matrix, $h_i$ is the restriction of $(h_{\gamma_i})_*$ to
the subspace 
$\bigoplus_{1 \leq \alpha \leq \alpha_i'} (v_i)_j(H_j(F_{i,\alpha}, \partial F_{i,\alpha}))$
(respecting the sum decomposition) and $h_{i,\alpha}$, for  $1 \leq \alpha \leq \alpha_i'$,  
is the restriction of $(h_{\gamma_i})_*$ 
to~$ (v_i)_j(H_j(F_{i,\alpha}, \partial F_{i,\alpha}))$ 
(see~\cite[Lemma~2.4]{DN}, \cite[Theorem~3.1]{NN00} or~\cite[Theorem~1.4]{NN00:bis}). In particular
$h_{i,\alpha}$ and $*_{i,s}$ are, respectively, the local blocks and the global blocks that we mentioned in Section~\ref{s:intro}. 
\end{observacion}

The distinguished bases $(\gamma_1, \ldots,\gamma_6)$ in  Figures~\ref{f:base_f+} 
and~\ref{f:base_f-} define a special geometric basis as we considered in Section~\ref{s:intro}. Then associated distinguished bases of
cycles $\mathcal{B}^{\pm}$ give the decomposition~\eqref{e:desc_homo}. We have shown 
the aforementioned submatrices $I_{j}, *_{i,s}, h_i$ and $h_{i,\alpha}$ 
in~\eqref{e:H0+} and~\eqref{e:H0-}.

\subsection{The polynomials $f^{\pm}$ are not topologically equivalent.} 
As we mentioned in Section~\ref{s:prop_top_f+-}, if $f^{\pm}$ were 
topologically equivalent then there should exist  $P$
satisfying~\eqref{igualmatriz_0}. In Section~\ref{s:grup_mon_f+-} we have 
computed the matrices $H_0^{\pm}$ and $H_1^{\pm}$ shown
in~\eqref{e:H0+} and~\eqref{e:H0-}. The problem becomes a linear system 
with 72 equations and 36 variables. The system of solutions is 
\[P(a,b):=\begin{scriptsize}
\left[ \begin {array}{cccccc} 0&0&a&a&0&a\\\noalign{\medskip}0&0&-a&-a&0&-a
\\\noalign{\medskip}0&0&a+b&a+b&0&a+b\\\noalign{\medskip}0&0&b&b&0&b
\\\noalign{\medskip}0&0&-b&-b&0&-b
\\\noalign{\medskip}0&0&a+b&a+b&0&a+b\end {array} \right] 
\end{scriptsize}\]
with $a,b \in \C$ (solved with \texttt{Maple}). 
Since $P(a,b)$ is a singular matrix, the system~\eqref{igualmatriz_0}
has no solution and hence the polynomials $f^{\pm}$ are not topologically equivalent.

\section{Proof of Theorem~\ref{t:met_discrim}}\label{s:proof}

This last section is dedicated to proving Theorem~\ref{t:met_discrim}. 
Hence we should assume that $f\in \C[x,y]$ is a monic polynomial of degree $N$ 
in the variable $y$ and satisfying Hypothesis~\ref{hip2cap2}.

\begin{lema}\label{l:simp_1}
It is enough to prove Theorem~\ref{t:met_discrim} for $j=1$.
\end{lema}
\begin{proof}
Let us assume that Theorem~\ref{t:met_discrim} is true for $j=1$. 
By Lemma~\ref{l:reduc_1}\eqref{l:reduc_1.a} and 
Theorem~\ref{t:met_discrim}\eqref{met_dis.2} (case $j=1$) we have  
$\nabla_{\star,\mathfrak{O}}(\gamma_i)=
(\alpha_j^{-1} \beta_i)^{-1}\sigma_j(\alpha_j^{-1} \beta_i)$
and thus Theorem~\ref{t:met_discrim}\eqref{met_dis.2} is true for $j>1$. 

By Lemma~\ref{l:reduc_1}\eqref{l:reduc_1.b} we have
$\mu_{j}^{(\alpha_j^{-1} \beta_i) }=\mu_1^{\beta_i } 
\ \mbox{and} \ \mu_{j+1}^{(\alpha_j^{-1} \beta_i) }=\mu_2^{ \beta_i }$
and hence,  since 
Theorem~\ref{t:met_discrim}\eqref{met_dis.1} it is true for $j=1$, 
it is also true for $j>1$.
\end{proof}

\subsection{Discriminant properties of a Morse polynomial.}
\begin{lema}\label{l:proy_fa}
Let $t_o \in T  \setminus \mathcal{T}_{Sing}$ and let $x^{t_o} \in {\bf x}^{t_o}$ be 
a ramification point of~\eqref{proy_prim}. In a neighborhood of $x^{t_o}$ the function $f(x,y)-t_o$ is the product of $N-1$ Weierstrass polynomials in the variable $y$ 
centered in the $N-1$ preimages of the point $x^{t_o}$. Also note that all those 
polynomials have degree 1 but one which has degree 2 and has either non-degenerate 
singularities if $(t_o,x^{t})=(t_i,x_i)$, for some $i$, or is
smooth if $(t_o,x^{t_o})\neq(t_i,x_i)$.
\end{lema}
\begin{proof}
Let $\bar{y}_j \in \C$, $j=1,\dots, r,$ be a root of $f(x^{t_o},y)-t_o.$
Then $(x^{t_o},\bar{y}_j)$, for $j=1,\dots, r,$ are all preimages of $x^{t_o}$. Since $x^{t_o}$
is a ramification point, $r \leq N-1$. If $r < N-1$, then by Corollary~\ref{propdiscrim} 
the discriminant germ $(\mathfrak{D}_{f},(t_o,x^{t_o}))$ is singular and hence  
$t_o \in \mathcal{T}_{Sing}$ which is a contradiction. 

Therefore $r=N-1$ and hence there exists a root $\bar{y}_1$ of order 2, whereas 
$\bar{y}_j$, for $j =2, \dots, N-1,$ are simple roots. By Hensel's Lemma $f(x,y)-t_o$ is, 
locally, the product of $N-1$ Weierstrass polynomials $P_j(x,y)$ in the variable $y$
centered in $(x^{t_o},\bar{y}_j)$ and whose degree is the multiplicity of $\bar{y}_j$. 
Thus $P_j(x,y)$, $j=2,\dots, N-1,$ have degree~1 and $P_1(x,y)$ has degree 2 where:

\begin{itemize}
\item If $(t_o,x^{t_o}) = (t_i,x_i)$, for some $i$, then $y_i$ is root of 
$f(x^{t_i},y)-t_i$, thus $\bar{y}_1=y_i$ and $P_1(x,y)$ has a non-degenerate singularity
at $(x_i,y_i)$.
\item If $t_o \notin \mathcal{T}$, then the fiber $V_{t_o}(f)$ is smooth. Finally, if $t_o=t_i$ 
but $x^{t_i} \neq x_i$, then the point $(x^{t_i},\bar{y}_1)$ is not singular.
In either case, $P_1(x,y)$ is smooth.
\end{itemize} 
\end{proof}

\begin{corolario}\label{Eculocalgen}
In a neighborhood of~$t_o \in T \setminus \mathcal{T}_{Sing}$ the function~$d_f(t,x)$ is:
\begin{enumerate}[(a)]
\item 
$a \prod_{\alpha=1}^k \ell_{\alpha}(x,t),$ if $t_o \notin f(P)$.
\item\label{Eculocalgen.b} 
$a \ q(x,t) \prod_{\alpha=2}^{k-1} \ell_{\alpha}(x,t)$,
if $t_o=t_i$, for some $i$.
\end{enumerate}
where 
$$
\ell_{\alpha}(x,t):=(x-x^{t_o}_{\alpha})+c_{\alpha}(t)(t-t_o),
q(x,t):=(x-x_i)^2+b_1(t) (x-x_i)+ c_1(t)(t-t_o)
$$
for $x_{\alpha}^{t_o} \in {\bf x}^{t_o}$, $a\in \C^*, b_1,c_{\alpha}\in \C[t]$, $b_1(t_o)=0$ and $c_{\alpha}(t_o)\neq 0, \alpha=1,...,k,$
\end{corolario}
\begin{proof}
By Lemma~\ref{l:proy_fa} and Lemma~\ref{desc_disc} the function  $d_f(\tilde{t}+t_o,\tilde{x}+x_{\alpha}^{t_o})$ is:
\begin{enumerate}[-]
\item a polynomial $\tilde{x}$-regular of degree $0+N-(N-1)=1$ and
$\tilde{t}$-regular of degree $N-(N-1)=1$, if $(t_o,x_{\alpha}^{t_o}) \neq (t_i,x_i)$.
\item a polynomial $\tilde{x}$-regular of degree $1+N-(N-1)=2$ and 
$\tilde{t}$-regular of degree $N-(N-1)=1$, if $(t_o,x_{\alpha}^{t_o}) = (t_i,x_i)$, for some $i$.
\end{enumerate}
The result follows by uniqueness of Weierstrass Theorem, Hensel's Lemma and the fact 
that $d_f(t,x)$ as polynomial in the variable $x$ has degree~$k$.
\end{proof}

\begin{corolario}\label{o:fa_hip_mod_top} 
The image of each member of a geometric basis by 
the classical monodromy $m_{t}:\pi_1(X \setminus {\bf x}^{t}; *) 
\rightarrow \Sigma_N, t \in T \setminus \mathcal{T}, * \in \partial X,$ 
associated with~\eqref{proy_prim} is a transposition.
\end{corolario}
\begin{proof}
Since $t \in T \setminus \mathcal{T}$, Lemma~\ref{l:proy_fa} gives us the existence
of a smooth Weierstrass polynomial of degree~2, say $P(x,y)$. If we lift the path 
$x^t+\rho \, e^{2 \pi u}, u \in [0,1],$ for a sufficiently small $\rho$, we obtain, 
again by Lemma~\ref{l:proy_fa}, a transposition given by the move of the two roots of 
$P(x^t+\rho,y)=0$. So the corollary follows.
\end{proof}

\subsection{Braid monodromy.}\label{s:meridiano}
In this section we prove Theorem~\ref{t:met_discrim}\eqref{met_dis.2}
\begin{notacion}\label{n:simplificacion}
To simplify notation, we will denote 
$\tau_{{\bf x}^{t}, \mathfrak{O}(\psi)},t \in \C$, $\tau(\mathfrak{O}(\psi))$ and
$\nabla_{ \star,\mathfrak{O}(\psi)}$ simply by  $\tau_{t,\psi}$, $\tau(\psi)$ and
$\nabla_{ t,\psi}$. 
\end{notacion}

\begin{proposicion}\label{deltaij}
Let $S$ be a circle with support in $T \setminus \mathcal{T}$, based at 
$\diamond \in \C$ and centered at $t_i \in f(P)$. If the radius of $S$, 
$\xi:=||t_i - \diamond||$, is sufficiently small, then for almost all $\psi \in (0,\pi)$,
$\nabla_{\diamond,\psi}(S)=\sigma_{j}$, where $1\leq j\leq k$ and $j$ depends of $\psi$.
\end{proposicion}

\begin{proof} 
Let $\theta:=(\bar{\pi} \vert_{\mathfrak{D}_f})_*(S)$. The braid $\theta$ is the 
homotopy class of
\begin{equation*}
\begin{aligned}
\{(t_i+\xi e^{2 \pi \sqrt{-1} u}, x) \, \vert \, x \in {\bf x}^{t_i+\xi e^{2 \pi \sqrt{-1} u}}, u \in [0,1]\}
\end{aligned}
\end{equation*}
By Corollary~\ref{Eculocalgen}\eqref{Eculocalgen.b} the set ${\bf x}^{t_i}$ has $k-1$ 
points. Almost every $\psi \in  (0,\pi)$ has the property that all lines with slope
$\tan(\psi)$ have no more than one point in ${\bf x}^{t_i}$. Let~$\psi$ be one of 
such angles. Let also be $j \in \{1, \dots, k-1 \}$ such that 
$x_i=x_{j}^{t_i} \in {\bf x}^{t_i}(\psi)$. 

By Corollary~\ref{Eculocalgen} (see Figure~\ref{discrimlocal}), for a sufficiently small $\xi$, the points of the 
list~${\bf x}^{t_i+\xi e^{2 \pi \sqrt{-1} u} }(\psi)$ satisfy 
\begin{equation*}%\label{e:ec_casi_def}
x_j^{t_i+\xi e^{2 \pi \sqrt{-1} u}} \stackrel{ \xi \rightarrow 0}{\longrightarrow} 
\begin{cases}
x_j^{t_i} & \mbox{for} \ j < j \\
x_{j}^{t_i}=x_i & \mbox{for} \ j = j, j+1 \\
x_{j-1}^{t_i} & \mbox{for} \ j > j+1
\end{cases}
\end{equation*}
and the local equation of $\mathfrak{D}_{f}$ at $(t_{i},x_i)$ is 
%\begin{equation}\label{e:ec_def}
$(x-x_i)^{2}+b_1(t) (x-x_i) +c_1(t)(t-t_i),$ 
$b_1(t_i)=0,c_1(t_i)\neq 0.$
%\end{equation}
Thus the braid 
$\nabla_{\diamond,\psi}(S)=\theta(\psi)=\tau_{\diamond,\psi} \, \theta \, \tau^{-1}_{\diamond,\psi}$ is equal to $\sigma_{j} \in \B_{k}$ by definition 
of $\tau_{\diamond,\psi}$ (see Notations~\ref{not:orden}
and~\ref{n:simplificacion}).
\begin{figure}[ht]
%\begin{center}
\hspace*{-2.5cm}
\includegraphics[scale=.7]{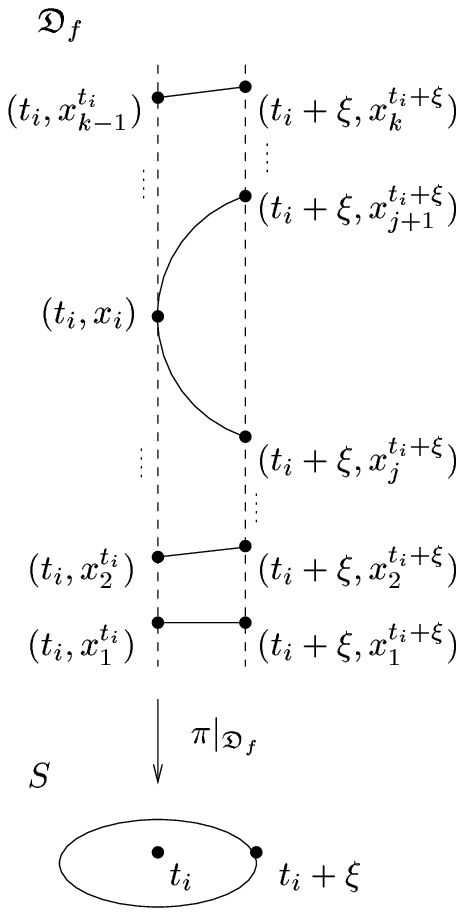}
\caption{}
%\end{center}
\label{discrimlocal}
\end{figure}
\end{proof}

Let $\gamma_i$ be an element of the distinguished basis of the Theorem~\ref{t:met_discrim} 
and let $\gamma$ be a meridian homotopic to $\gamma_i$ as in Definition~\ref{d:base_dis_cam}
such that $\supp \gamma \subset T \setminus  \mathcal{T}$. According to
Definition~\ref{defmeridiano} the following elements are associated with ${\gamma}$
\begin{equation}\label{e:ele_asoc}
(t_{\gamma}(=t_i), t'_{\gamma}(=: \diamond ), 
r_{\gamma}(=:r), S_{\gamma}(=:S), D_{\gamma}(=:D)).
\end{equation}
Let $\psi \in (0,\pi)$. We need the following technical hypothesis:
\begin{hipotesis}\label{o:buena_base_distin}
Proposition~\ref{deltaij} is satisfied for the element $S$ of $\gamma$ and $\psi$.
\end{hipotesis}
 
\begin{notacion}\label{n:mon}
We denote $\theta:=(\tilde{\pi}\vert_{\mathfrak{D}_{f}})_*(S)$ and 
$\beta:=(\tilde{\pi}\vert_{\mathfrak{D}_{f}})_*(r^{-1}) $ 
the braids in~$\B({\bf x}^{\diamond},{\bf x}^{\diamond})$ and $\B({\bf x}^{\diamond}, 
{\bf x}^{\star}),$ respectively.
\end{notacion} 
By the first part of Proposition~\ref{ident_trenzas}\eqref{ident_trenzas.1} and 
since braid monodromy is a homotopy invariant we have
\begin{equation}\label{e:desc_nabla}
\nabla_{ \star,\mathfrak{O}}(\gamma_i)=\nabla_{ \star,\mathfrak{O}}(\gamma)=
\beta^{-1}(\mathfrak{O}) \theta(\mathfrak{O}) \beta(\mathfrak{O})
\end{equation}
and, by Proposition~\ref{deltaij}, 
\begin{equation}\label{e:desc_nabla.2}
\nabla_{ \star,\psi}(\gamma_i)=\nabla_{\star,\psi}(\gamma)= \beta (\psi )^{-1}  \theta(\psi)  \beta(\psi), 
\ \ \mbox{where} \ \theta(\psi)=\sigma_{j}.
\end{equation}
By the second part of Proposition~\ref{ident_trenzas}\eqref{ident_trenzas.2} 
the braid $\theta(\mathfrak{O})$ is conjugated to $\theta(\psi)=\sigma_{j}$. 
Then Theorem~\ref{t:met_discrim}\eqref{met_dis.2} follows from Lemma~\ref{l:reduc_1}\eqref{l:reduc_1.a}.

\subsection{The vanishing cycle of a critical value.}\label{s:cons_ciclos_evan}

In this section we prove 
Theorem~\ref{t:met_discrim}\eqref{met_dis.1.a} and \eqref{met_dis_Delta}.
To do this we will first prove it in the special case when $\gamma_i$ is in fact $\gamma$
and $\mathfrak{O}$ is the base set of braids $\mathfrak{O}(\psi)$ with $\gamma$ and $\psi$
as in Section~\ref{s:meridiano}. This is done in the next three subsections. 
Finally, in the last subsection, we prove the theorem in the general case. 

\subsubsection{Construction of the vanishing cycle.}\label{cons_ciclos_evan}

\begin{figure}[ht]
%\begin{center}
\vspace*{-.1cm}
\hspace*{-3.5cm}
\includegraphics[scale=.8]{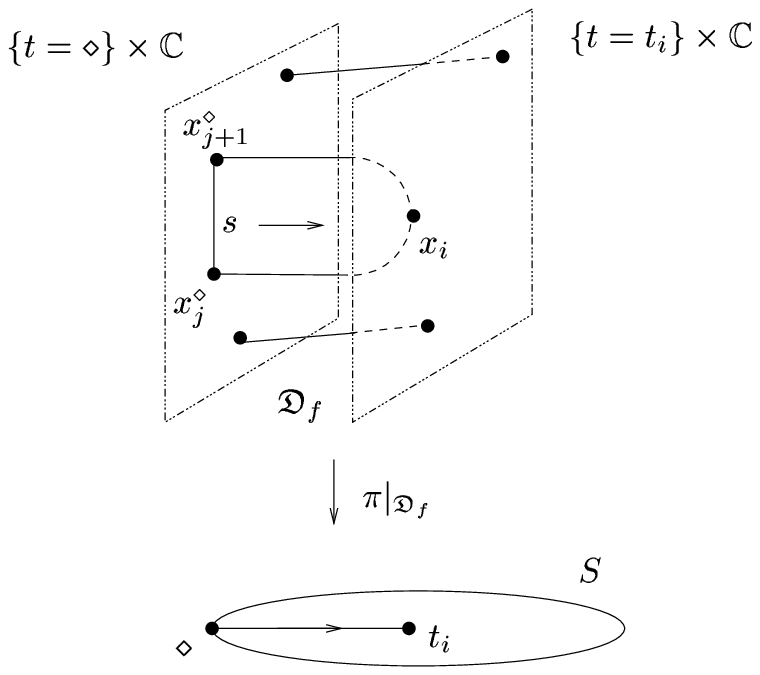}
\caption{}
%\end{center}
\label{consevans1}
\end{figure}

Let us parametrize the segment $[\diamond,t_i]$ as 
$\lambda(u):=\left(1-u\right)(\diamond-t_i)+t_i,  u \in [0,1].$
Let $s$ be the segment in $\{t=\diamond\} \times \C$ between the points 
$x_{j}^{\diamond}$ and $x_{j+1}^{\diamond}$ of ${\bf x}^{\diamond}({\psi})$.
By Corollary~\ref{Eculocalgen}, while $u$ moves towards $1$, the points 
$x_{j}^{\gl(u)}, x_{j+1}^{\gl(u)} \in {\bf x}^{\gl(u)}({\psi})$ 
converge to $x_i$ then the path $s$ degenerates into a point $(t_i,x_i)$ as we move towards $t_i$
on the interval $[\diamond,t_i]$ (see Figure~\ref{consevans1}). 

For every $\gl(u) \in [\diamond,t_i)$ let $H_{\tau_{\gl(u),\psi}}$ be an 
$({\bf x}, {\bf x}^{\gl(u)})$-special isotopy associated with $\tau_{\gl(u),\psi}$, 
which is the identity in $\C \setminus \Int(X)$ such that the map
\begin{equation}\label{e:H}
\begin{matrix}
H: & [0,1) \times (\C \times [0,1]) &\longrightarrow & \C \\
& (u,x,v) & \mapsto & H_{\tau_{\gl(u),\psi}}(x,v)
\end{matrix}
\end{equation}
is continuous in the variable $u$.  By continuity  we can extend $H$ to 
$[0,1]\times (\C \times [0,1])$. 

We define a topological cylinder $\mathfrak{C}$ with 
$
\mathfrak{C}_{u}:=  H \left(u, \D_1\left(\frac{2j+1}{2}\right),1 \right), 0 \leq u \leq 1,
$ 
its fibers, 
where $j$ is the element given by Proposition~\ref{deltaij} for the fixed $\psi$.
The cylinder $\mathfrak{C}$ satisfies that
$$
x_{j}^{\gl(u)},  x_{j+1}^{\gl(u)} \in \intt(\mathfrak{C}_{u}), \
\mathfrak{C}_{u} \cap ({\bf x}^{\gl(u)} \setminus \{x_{j}^{\gl(u)}, 
x_{j+1}^{\gl(u)} \})= \emptyset, \quad 0 \leq u < 1,
$$
$$
x_i \in \intt(\mathfrak{C}_{1})  \ \mbox{and} \
 \mathfrak{C}_{1}  \cap  
({\bf x}^{t_i} \setminus \{x_i \})= \emptyset,
$$

By Lemma~\ref{l:proy_fa} the set 
$\phi_{f,x}^{-1}(\{t=\gl(u)\} \times \mathfrak{C}_{u})$ 
has $N-1$ components where $N-2$ of them are diffeomorphic to a topological disk 
and the remaining one is diffeomorphic to the cylinder 
$\mathbb{S}^{1} \times [-1,1]$ for $0 \leq u <1$, and to two disks intersecting 
at the point $(x_i,y_i)$ for~$u=1$.

\begin{figure}[ht]
%\begin{center}
\vspace*{-0.7cm}
\hspace*{-3cm}
\includegraphics[scale=.8]{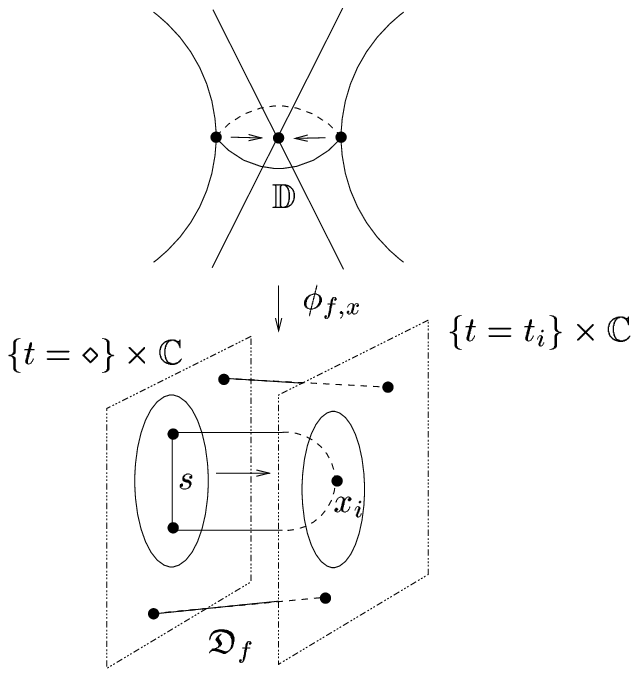}
\caption{}
%\end{center}
\label{esquemasimple}
\end{figure}

Let $\mathbb{S}(0)$ be the preimage of $s$ in 
$\phi_{f,x}^{-1}(\{t=\diamond\} \times \mathfrak{C}_{0})$. Note that $\mathbb{S}(0)$
is a closed path according to~\eqref{proy_prim}.
Analogously, one obtains a family of closed paths $\mathbb{S}(u)$, $0<u<1$ as 
the path $s$ is moved along the segment $[\diamond,t_i]$. We denote also 
$\mathbb{S}(1):=\{(x_i,y_i)\}$. Hence there is a 2-dimensional closed topological 
disk $\D$ satisfying (see Figure~\ref{esquemasimple}): 
\begin{enumerate}[i.]
\item $\mathbb{S}(u)=\D \cap V_{\lambda(u)}(f)$, $0 \leq u <1$, is diffeomorphic to the circle.
\item $\mathbb{S}(1)=\D \cap V_{t_i}(f)=\{(x_i,y_i)\}$.
\item $\D=\bigcup \mathbb{S}(u)$.
\end{enumerate}
Therefore $\mathbb{S}(0)$ defines a vanishing cycle in~$V_{\diamond}(f)$.

\subsubsection{Representation of the vanishing cycle in $H_1(V_{\diamond}(f))$.}
\label{s:repcicloevan}
Since $\supp(\gamma) \subset  T \setminus \mathcal{T}$ we can use the topological
model described in Section~\ref{modtopyhom} to give a representation of $\mathbb{S}(0)$ in~$H_1(V_{\diamond}(f))$.

Let us consider the geometric basis $(\mu_1,...,\mu_k)$ as in Example~\ref{ejemplobase} 
associated with the disk~$X$ fixed above and the set ${\bf x}=\{1,2,...,k\}$
(by~\eqref{s:dis_met_cond3} in Section~\ref{s:dis_met} ${\bf x}\subset X$). By Artin Theorem~\cite{Artin47} the 
Hurwitz move $\Psi_{\tau_{\diamond,\psi}}$ transforms this geometric basis 
into the geometric basis 
$(\mu_{1}^{\tau_{\diamond,\psi}}, \dots, \mu_{k}^{\tau_{\diamond,\psi}})$
of $\pi_1( X \setminus {\bf x}^{\diamond}; *)$. 

For the sake of simplicity we will skip the subindices $\diamond,\psi$ from
$\tau_{\diamond,\psi}$. By Definition~\ref{defmeridiano} the elements associated with
$\mu_{n}^{\tau}$ are
\begin{equation}\label{e:el_tau_i}
(x_{\mu_{n}^{\tau}}(=x_{n}^{\diamond}),
x_{\mu_{n}^{\tau}}',
r_{\mu_{n}^{\tau}},
S_{\mu_{n}^{\tau}},
D_{\mu_{n}^{\tau}}).
\end{equation} 

By construction of the vanishing cycle, the classical monodromy $m_{\diamond}$ satisfies
\begin{equation}\label{tras.originaria}
m_{\diamond}(\mu_{j}^{\tau})=m_{\diamond}(\mu_{j+1}^{\tau}).
\end{equation}
Moreover, by Corollary~\ref{o:fa_hip_mod_top} $m_{\diamond}(\mu_{j}^{\tau})$ is a transposition, say $(a,b), a,b \in \{1, \dots, N\},$ $ a \neq b$. 

Let $\delta:=\mu_{j+1} \cdot \mu_{j}$ be a path in $X \setminus {\bf x}$. We want to
lift $\Psi_{\tau}(\delta)=\delta^{\tau}=\mu_{j+1}^{\tau} \, \cdot \, \mu_{j}^{\tau}$
by means of the projection $\pi\vert_{V_{\diamond}(f)}$~\eqref{proy_prim}.
This lifting is a union of $N$ closed paths in the fiber $V_{\diamond}(f)$:  
\begin{equation*}
(\delta^{\tau})_{l}=(\mu_{j+1}^{\tau})_{l} \cdot 
(\mu_{j}^{\tau})_{l}  \, \, 
\mbox{with $l \in \{1,  \dots,N\}$}, l \neq a, b,
\end{equation*}
\begin{equation*}
(\delta^{\tau})_{a}=(\mu_{j+1}^{\tau})_{a} \cdot 
(\mu_{j}^{\tau})_{b} \ \mbox{and}
 \  
(\delta^{\tau})_{b}=(\mu_{j+1}^{\tau})_{b} \cdot (\mu_{j}^{\tau})_{a}. 
\end{equation*}
 
The following lemma allows us to find paths that are homologous to the vanishing cycle~$\mathbb{S}(0)$ (Section~\ref{cons_ciclos_evan}). These paths will let us 
transport easily~$\mathbb{S}(0)$ to the fiber~$V_{\star}(f)$.

\begin{lema}\label{l:ciclo_evan_bas}
The liftings $(\delta^{\tau})_{l}, l=1, \dots, N,$ in $V_{\diamond}(f)$ 
are homologous -up to sign- to the vanishing cycle $\mathbb{S}(0)$ 
if  $l = a,b$ and to the zero cycle if $l \neq a,b$. Also in homology
$(\delta^{\tau})_{b}=-(\delta^{\tau})_{a}$.
\end{lema}
\begin{proof}
Let us deform the path $s$ continuously in 
\begin{equation*}
\{t=\diamond\}\times ( X  \setminus \bigcup_{n \neq j, \, j+1} 
x_{n}^{\diamond} ) 
\end{equation*}
as shown in Figure~\ref{esquematodoarriba} (a), 
so as to turn it into a path with support
$$
[x_{j}^{\diamond},x_{\mu^{\tau}_{j}}'] \cup 
\supp(r_{\mu^{\tau}_{j}}^{-1}) \cup \supp(r_{\mu^{\tau}_{j+1}}) 
\cup [x_{\mu^{\tau}_{j+1}}', x_{j+1}^{\diamond}],
$$
where these paths are the ones associated with~$\mu_j^{\tau}$ (see~\eqref{e:el_tau_i}).
The support of $\delta^{\tau}$ is 
$$
\supp(S_{\mu^{\tau}_{j}}) \cup \supp(r_ {\mu^{\tau}_{j}}^{-1}) \cup 
\supp(r_ {\mu^{\tau}_{j+1}}) \cup \supp(S_{\mu^{\tau}_{j+1}}). 
$$
Therefore the cycles $\mathbb{S}(0)$, $(\delta^{\tau})_{a}$ and 
$(\delta^{\tau})_{b}$ in $V_{\diamond}(f)$ (Figure~\ref{esquematodoarriba} (b))
are homologous up to orientation.
\begin{figure}[ht]
%\begin{center}
\vspace*{-5cm}
\hspace*{-2.5cm}
\includegraphics[scale=.8]{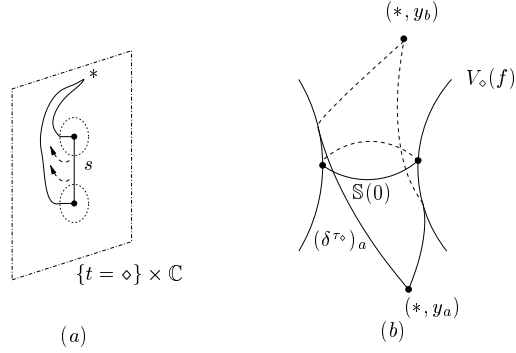}
\caption{(a) Deformation of $s$. (b) Paths $\mathbb{S}(0)$ and $(\delta^{\tau})_{a}$.}
%\end{center}
\label{esquematodoarriba}
\end{figure}

By Proposition~\ref{propcamhomolog}\eqref{propcamhomolog.c} and~\eqref{tras.originaria}
the $N$ paths $(\delta^{\tau})_{l}$ are elements in $Z_1(\mathfrak{K})$ and 
hence, by Proposition~\ref{propcamhomolog}\eqref{propcamhomolog.a} as cycles in~$H_1(V_{\diamond}(f))$ (see~\eqref{rel_grup_homolog}) they look like:
\begin{equation*}
\begin{aligned}
(\delta^{\tau})_{l}&=(\mu^{\tau}_{j+1})_{l} + (\mu^{\tau}_{j})_{l}=
\ 0  
\, \, \mbox{with $l \in \{1,  \dots,N\}$}, l \neq a, b, \\
(\delta^{\tau})_{a}&= (\mu^{\tau}_{j+1})_{a} + (\mu^{\tau}_{j})_{b}
 \, \, \mbox{and}\\
(\delta^{\tau})_{b}&= (\mu^{\tau}_{j+1})_{b} + (\mu^{\tau}_{j})_{a}= 
- (\mu^{\tau}_{j+1})_{a} - (\mu^{\tau}_{j})_{b}=-(\delta^{\tau})_{a}.
\end{aligned}
\end{equation*}
\end{proof}

\subsubsection{Representation of the vanishing cycle 
in $H_1(V_{\star}(f))$.}\label{s:cons.difeo}

In this subsection we will compute the homological class of the local 
vanishing cycle $\mathbb{S}(0)$ in the global fiber $H_1(V_{\star}(f))$. 
In order to do so we will use Lemma~\ref{l:ciclo_evan_bas} which allows us
to work with $(\delta^{\tau})_{a}$ or $(\delta^{\tau})_{b}$ and move them
from $V_{\diamond}(f)$ to $V_{\star}(f)$ by the diffeomorphism~$h_{r^{-1}}$ 
given by lifting the path~$r^{-1}$. Let us first study the 
diffeomorphism $h_{r^{-1}}$ more closely.

\paragraph{\textit{The diffeomorphism $h_{r^{-1}}$.}}

Let $\beta \, \in \B({\bf x}^{\diamond},{\bf x}^{\star})$ the braid defined by $r^{-1}$ as in 
 Notation~\ref{n:mon} and $H_{\beta}$ a 
$({\bf x}^{\diamond},{\bf x}^{\star})$-special isotopy associated to $\beta$ 
which is the identity in $\C \setminus \Int(X)$. Then we can define
\begin{equation*}
\array{rccl}
\bar{H}: & \{\diamond\} \times \C \times [0,1] &
\longrightarrow &  \supp(r^{-1}) \times  \C. \\

 & (\diamond, x,u) & \mapsto & (r^{-1}(u), H_{\beta}(x,u))
\endarray
\end{equation*}
Let
$ \mathcal{X}:=\{(r^{-1}(u),x) \in T \times X \mid u \in [0,1], x \notin {\bf x}^{r^{-1}(u)} \}
 \  \mbox{and}  
\
 \mathcal{R}:=\phi_{f,x}^{-1}(\mathcal{X})
$
with fibers
\begin{equation}\label{X_F}
\begin{matrix} 
\mathcal{X}_u:= \mathcal{X} \cap (\{ r^{-1}(u)\} \times \C) \  \mbox{and} \

 \mathcal{R}_u:=\phi_{f,x}^{-1}(\mathcal{X}_u), u \in [0,1].
\end{matrix}
\end{equation}
Then we have the following $N$-folded covering maps
$$\phi_{f,x} \vert : \mathcal{R} \longrightarrow  \mathcal{X} \ \ \mbox{and}
\ \
(\phi_{f,x} \vert,1) : \mathcal{R}_0 \times [0,1] \longrightarrow  \mathcal{X}_0 \times [0,1] 
$$
and a diffeomorphism between the base spaces
$\bar{H}\vert: \mathcal{X}_0 \times [0,1] \rightarrow \mathcal{X}.$
Let $\tilde{*} \in  \mathcal{R}_0$ be such that
$\phi_{f,x}(\tilde{*})=(\diamond,*) \in \mathcal{X}_0$. 
Since $r^{-1}(u)$ ($u \in [0,1])$ are all regular values of $f$ and 
$\pi \vert_{\mathfrak{D}_{f}}$, one has
$\pi_1( \mathcal{X};*) \simeq \pi_1( \mathcal{X}_0;*) \  \mbox{and}  \
\pi_1( \mathcal{R};\tilde{*}) \simeq \pi_1( \mathcal{R}_0;\tilde{*}).$
Therefore, by Lifting Theorem, the diffeomorphism between the base spaces gives us a unique 
diffeomorphism $\Gamma$ between the covering spaces such as $\Gamma(\tilde{*},0)=\tilde{*}$ 
and such that the following diagram is commutative
\begin{equation}\label{diagrama_super}
\begin{CD}
\mathcal{R}_0 \times  [0,1] @>\Gamma>> \mathcal{R} \\
@V(\phi_{f,x},1)VV @VV\phi_{f,x}V \\
\mathcal{X}_0 \times [0,1] @>\bar{H}\vert>>  \mathcal{X}.
\end{CD}
\begin{picture}(0,0)
%\graphpaper(-90,-50)(100,75)
\put(-36,-2){$\circlearrowleft$}
\end{picture}
\end{equation}

By restriction we have the following families of diffeomorphisms 
\begin{equation}\label{e:Hu}
\begin{matrix}
\bar{H}_u: &\{\diamond\} \times  \C & \longrightarrow & \C\\
 & (\diamond,x) & \mapsto &(r^{-1}(u),H_{\beta}(x,u)) \\
\Gamma_u: & \mathcal{R}_0 & \longrightarrow & \mathcal{R}_u\\
 & (x,y) & \mapsto & \Gamma(x,y,u).
\end{matrix}
\end{equation}
Each $\bar{H}_u$ is such that $\bar{H}_u(\mathcal{X}_0)=\mathcal{X}_u$ and 
each $\Gamma_u$ is such that $f(\Gamma_u(\bullet))=r^{-1}(u)$
and it can  be extended continuously to 
$\phi_{f,x}^{-1}(\{t=\diamond\}\times {\bf x}^{\diamond})$. 
Finally we define the diffeomorphism $h_{r^{-1}}:=\Gamma_1$. Note that, by~\eqref{diagrama_super}, the following is a commutative diagram:
\begin{equation}\label{diagrama_sub}
\begin{CD}
\mathcal{R}_0 @>h_{r^{-1}} \vert>> \mathcal{R}_1 \\
@V\phi_{f,x}VV @VV\phi_{f,x}V \\
\mathcal{X}_0 @>\bar{H}_1\vert>> \mathcal{X}_1.
\end{CD}
\begin{picture}(0,0)
%\graphpaper(-90,-50)(100,75)
\put(-38,-2){$\circlearrowleft$}
\end{picture}
\end{equation}
Moreover, using a theorem of continuity of the roots for
$f(x,y)=r^{-1}(u),$ with $x \in \C\setminus \intt(\mathcal{X}_0)$ and  $u\in~[0,1],$
and since $\bar{H}_u$ is the identity on the 
second variable in $\{\diamond\} \times \C\setminus \intt(
\mathcal{X}_0)$, this diagram can be extended to $V_{\diamond}(f)$.

\begin{observacion}\label{o:inv}
Let $H_{\beta^{-1}}(\bullet,u):=H_{\beta}^{-1}(\bullet,u)$ be an
$({\bf x}^{\star},{\bf x}^{\diamond})$-special isotopy associated with $\beta^{-1}$. 
Analogously as above, one can construct the diffeomorphism $h_{r}$ and the 
commutative diagram:
\begin{equation}\label{diagrama_subnv}
\begin{CD}
\mathcal{R}_1 @>h_{r}\vert>> \mathcal{R}_0 \\
@V\phi_{f,x}VV @VV\phi_{f,x}V \\
\mathcal{X}_1 @>\bar{H}_1^{-1}\vert >> \mathcal{X}_0.
\end{CD}
\begin{picture}(0,0)
%\graphpaper(-90,-50)(100,75)
\put(-38,-2){$\circlearrowleft$}
\end{picture}
\end{equation}
\end{observacion}

\paragraph{\textit{Moving the paths $(\delta^{\tau})_{a}$ and  
$(\delta^{\tau})_{b}$.}}
By Artin Theorem $ (\mu^{\tau  \beta}_{1}, \dots, \mu^{\tau \beta}_{k})$ 
is a geometric basis of $\pi_1( X \setminus  {\bf x}^{\star}; *)$.
Let $\mathcal{Y}_{r^{-1}}$ be the braid in $N$ strings defined by the equation $f(*,y)=r^{-1}(u),u \in  [0,1]$. If we consider $\mathcal{Y}_{r^{-1}}$ as a permutation, by~\eqref{tras.originaria}, we have
\begin{equation}\label{e:nuevos_m_n}
m_{\star}(\mu_{j}^{\tau  \beta})=
m_{\star}(\mu_{j+1}^{\tau  \beta})=
(a^{\mathcal{Y}_{r^{-1}}},b^{\mathcal{Y}_{r^{-1}}})=:(\bar{a},\bar{b}).
\end{equation}

\begin{lema}\label{l:sgun_ciclo_evan_bas}
Let $\bar{l}:=l^{\mathcal{Y}_{r^{-1}}},$ with $l=1, \dots, N$. 
The liftings~$(\delta^{\tau \beta})_{\bar{l}}$ in $V_{\star}(f)$ are homologous to
the vanishing cycle $(h_{r^{-1}})_*(\mathbb{S}(0))$ if $\bar{l} = \bar{a},\bar{b}$ 
and are the zero cycle if $\bar{l} \neq \bar{a},\bar{b}$.
Moreover, $(\delta^{\tau \beta})_{\bar{a}}=-(\delta^{\tau \beta})_{\bar{b}}=
(\mu_{j+1}^{\tau \beta})_{\bar{a}} +  
(\mu_{j}^{\tau  \beta})_{\bar{b}}$.
\end{lema}
\begin{proof}
By the commutative diagram~\eqref{diagrama_sub}
and by equation~\eqref{e:nuevos_m_n} we have the following equalities in homotopy
\begin{equation*}
\array{rccccl}
h_{r^{-1}}((\delta^{\tau})_{l})&=&
(\delta^{\tau  \beta})_{\bar{l}}&=& 
(\mu_{j+1}^{\tau  \beta})_{\bar{l}} \cdot   
(\mu_{j}^{\tau  \beta})_{\bar{l}} & \ \mbox{for} \ \bar{l} \neq \bar{a}, \bar{b}, \\
h_{r^{-1}}((\delta^{\tau})_{a})&=&
(\delta^{\tau  \beta})_{\bar{a}} &=&
(\mu_{j+1}^{\tau  \beta})_{\bar{a}} \cdot   
(\mu_{j}^{\tau  \beta})_{\bar{b}} & \ \mbox{and}\\
h_{r^{-1}}((\delta^{\tau})_{b})&=&
(\delta^{\tau  \beta})_{\bar{b}}&=&
(\mu_{j+1}^{\tau  \beta})_{\bar{b}} \cdot  
(\mu_{j}^{\tau  \beta})_{\bar{a}} &
\endarray       
\end{equation*}
By Proposition~\ref{propcamhomolog}\eqref{propcamhomolog.c} the $N$ 
paths~$\delta^{\tau \beta }_i$ are elements in~$Z_1(\mathfrak{K})$ and 
hence by Proposition~\ref{propcamhomolog}\eqref{propcamhomolog.a} as cycles in $H_1(V_{\star}(f))$ are homologous to (see~\eqref{rel_grup_homolog}):
\begin{equation*}
\array{rl}
(h_{r^{-1}})_*((\delta^{\tau})_{l})&=
(\delta^{\tau  \beta})_{\bar{l}}= 
0  \ \mbox{for} \ l \neq \bar{a}, \bar{b}, \\
(h_{r^{-1}})_*((\delta^{\tau})_{a})&=
(\delta^{\tau  \beta})_{\bar{a}} =
(\mu_{j+1}^{\tau  \beta})_{\bar{a}} +   (\mu_{j}^{\tau  \beta})_{\bar{b}}
 \ \mbox{and}\\
(h_{r^{-1}})_*((\delta^{\tau})_{b})&=
(\delta^{\tau  \beta})_{\bar{b}}=
(\mu_{j+1}^{\tau  \beta})_{\bar{b}}+
(\mu_{j}^{\tau  \beta})_{\bar{a}}
=-(\mu_{j+1}^{\tau  \beta})_{\bar{a}} -
(\mu_{j}^{\tau  \beta})_{\bar{b}} 
\endarray
\end{equation*}
Therefore $(\delta^{\tau  \beta})_{\bar{b}}=-(\delta^{\tau  \beta})_{\bar{a}}$ 
in homology and by Lemma~\ref{l:ciclo_evan_bas} they are -up to orientation- 
the image of~$\mathbb{S}(0)$ in $V_{\star}(f)$.
\end{proof}

According to the commutative diagram~\eqref{diagramaprimario}, where
${\bf s}^1= {{\bf x}^{\diamond}}, {\bf s}^2= {\bf x}^{\star}, \{1, \dots, k\}={\bf x}$, 
and $\mathfrak{O}=\mathfrak{O}(\psi)$, we have that
$\Phi_{{\bf x}^{\diamond},{\bf x}^{\star}}(\mu_n^{\tau},\beta) = 
(\Psi_{\tau_{\star,\psi}} \circ \Phi_{{\bf x}, {\bf x}})  \left( \mu_n ,
\beta \right)$ or, in exponential notation,
\begin{equation}\label{igualdad_mu_j}
\mu_n^{ \tau \beta}=\mu_{n}^{ \beta  \tau_{\star,\psi}}
\in \langle \mu_1^{\tau_{\star,\psi}}, \dots, 
\mu_k^{\tau_{\star,\psi}} \rangle.
\end{equation}
In particular, by~\eqref{e:nuevos_m_n}, $m_{\star}(\mu_{j}^{ \beta \tau_{\star,\psi}})$ 
is a transposition $(\bar{a},\bar{b})$. So,  we have proved  Theorem~\ref{t:met_discrim}\eqref{met_dis.1.a}. %(2.a)
Also we have by Lemma~\ref{l:sgun_ciclo_evan_bas} that in homology	
\begin{equation}\label{igualdad_Delta_j}
(\delta^{
\beta\tau_{\star,\psi}})_{\bar{a}}
=-(\delta^{
\beta\tau_{\star,\psi}})_{\bar{b}}=(\mu_{j+1}^{\beta\tau_{\star,\psi}})_{\bar{a}} + (\mu_{j}^{\beta\tau_{\star,\psi}})_{\bar{b}}
\end{equation}
defines -up to orientation- the vanishing cycle $(h_{r^{-1}})_*(\mathbb{S}(0))$ 
and the remaining liftings of $\delta^{\beta\tau_{\star,\psi}}$ define a zero cycle. 
So, finally, if we denote by $\Delta_i$ the first non-zero lifting in homology, we 
have proved Theorem~\ref{t:met_discrim}\eqref{met_dis_Delta}	 with decomposition 
$\nabla_{ \star,\psi}(\gamma_i)= \beta^{-1}  \sigma_j \beta$
in the case $\gamma_i=\gamma$ and $\mathfrak{O}=\mathfrak{O}(\psi)$.

\subsubsection{Independence of construction.}

Finally in this last subsection we prove Theo\-rem~\ref{t:met_discrim}\eqref{met_dis.1.a} and~\eqref{met_dis_Delta} 
for any $\gamma_i$ and $\mathfrak{O}$ with a decomposition as in Theorem~\ref{t:met_discrim}~(1) (for $j=1$ according to Lemma~\ref{l:simp_1}).

\begin{lema}\label{l:simp_2}
We have the decomposition 
$\nabla_{\star,\psi}({\gamma}_i)=(\beta \alpha_{j})^{-1} \sigma_1 (\alpha_{j} \beta)$.
Note that with this decomposition Theorem~\ref{t:met_discrim}\eqref{met_dis.1.a} and~\eqref{met_dis_Delta} are 
true for the base set of braids $\mathfrak{O}(\psi)$.
\end{lema}
\begin{proof}
The first part of lemma follows from~\eqref{e:desc_nabla.2} and 
Lemma~\ref{l:reduc_1}\eqref{l:reduc_1.a}. 
For the second part, by Lemma~\ref{l:reduc_1}\eqref{l:reduc_1.b}, we have 
$\mu_1^{(\alpha_{j} \beta) } =\mu_{j}^{\beta}\ \mbox{and} \ 
\mu_2^{(\alpha_{j} \beta)}=\mu_{j+1}^{\beta}.$ The lemma follows from 
Theorem~\ref{t:met_discrim}\eqref{met_dis.1.a} and~\eqref{met_dis_Delta} for the base set of braids 
$\mathfrak{O}(\psi)$ and $\sigma_j$.
\end{proof}

\begin{observacion}\label{o:simpl}
By Lemma~\ref{l:simp_2}, we can assume that~$\theta(\psi)=\sigma_1$.
\end{observacion}

Note that Proposition~\ref{ident_trenzas}\eqref{ident_trenzas.2} and 
Remark~\ref{o:simpl} imply that
$$
\theta(\mathfrak{O})=(\tau \tau_{\diamond,\mathfrak{O}}^{-1})^{-1} \sigma_{1} 
(\tau \tau_{\diamond,\mathfrak{O}}^{-1}),
$$
so, by~\eqref{e:desc_nabla}, we have
\begin{equation}\label{des_seg}
\nabla_{\star,\mathfrak{O}}(\gamma_i)=(\tau
\tau_{\diamond,\mathfrak{O}}^{-1} \beta(\mathfrak{O}))^{-1} \sigma_{1} (\tau
\tau_{\diamond,\mathfrak{O}}^{-1} \beta(\mathfrak{O})).
\end{equation}
\begin{notacion}\label{n:simpl}
We denote by
$\nu:= \tau\tau_{\diamond,\mathfrak{O}}^{-1} \beta(\mathfrak{O}) \in \B_k$
the conjugation braid in~\eqref{des_seg}. 
\end{notacion}

\begin{lema}\label{l:centralizador}
Let $C(\sigma_{1})$ the centralizer of $\sigma_{1}$ in $\B_{k}$ and 
$\zeta \in C(\sigma_{1})$. Then one of the following is true:
\begin{enumerate}
\item $\mu_{1}^{\zeta}$ and $\mu_{2}^{\zeta}$ are conjugated 
(by the same element) of $\mu_{1}$ and $\mu_{2}$ respectively.
\item $\mu_{1}^{\zeta}$ and $\mu_{2}^{\zeta}$ are conjugated 
(by the same element) of $\mu_{2}$ and $\mu_{1}$ respectively.
\end{enumerate}
\end{lema}
\begin{proof}
By~\cite{FRZ}  
$C(\sigma_{1})= \, \langle\sigma_{2} \sigma_{1}^{2} \sigma_{2}, 
\sigma_j \, \vert \, j \neq 2 \rangle.$ 
Then by~\eqref{formulaexp}:
\begin{itemize}
\item $\mu_{1}^{\sigma_j}= \mu_{1}$ and 
$ \mu_{2}^{\sigma_j}=\mu_{2}$  for $j \neq 1,2$, 

\item $\mu_{1}^{\sigma_{1}}=\mu_{2}=\mu_{2} \cdot \mu_{2} \cdot \mu_{2}^{-1}$  and 
$ \mu_{2}^{\sigma_{1}}=\mu_{2} \cdot \mu_{1} \cdot \mu_{2}^{-1}$,

\item $\mu_{1}^{\sigma_{2} \sigma_{1}^{2} \sigma_{2}}= 
\mu_{3} \cdot \mu_{1} \cdot \mu_{3}^{-1}$ and 
$\mu_{2}^{\sigma_{2} 
\sigma_{1}^{2} \sigma_{2}}=
\mu_{3} \cdot \mu_{2} \cdot \mu_{3}^{-1}$,
\end{itemize}
and the lemma follows.
\end{proof}

By Theorem~\ref{t:met_discrim}\eqref{met_dis.2} and~\eqref{des_seg}
$\nu^{-1} \sigma_{1} \nu= \beta_i^{-1} \sigma_{1} \beta_i,$ thus
$\beta_{i}\nu^{-1} \sigma_{1}=\sigma_{1} \beta_{i}\nu^{-1}$ and hence
$\beta_i=  \zeta \nu \ \mbox{with} \ \zeta \in C(\sigma_{1}).$
By Lemma~\ref{l:centralizador} there exists 
$\eta \in \pi_1(X \setminus {\bf x};*)$ such that $\mu_{1}^{\beta_{i}}$ and 
$\mu_{2}^{\beta_{i}}$ are equal to one of the following:
\begin{enumerate}[(i)]
\item\label{prim} $(\eta^{\nu })^{-1} \cdot \mu_{1}^{\nu }  \cdot \eta^{\nu}$
and $(\eta^{\nu })^{-1} \cdot \mu_{2 }^{\nu }  \cdot \eta^{{\nu } } $,
respectively, or 
\item\label{seg}  $(\eta^{\nu })^{-1} \cdot \mu_{2 }^{\nu }  \cdot \eta^{\nu } $ and $(\eta^{\nu })^{-1} \cdot \mu_{1}^{\nu}  \cdot \eta^{\nu }$,
respectively.
\end{enumerate}

By isomorphism~\eqref{iso1} 
$\beta(\mathfrak{O})= \tau_{\diamond,\mathfrak{O}} \beta
\tau_{\star,\mathfrak{O}}^{-1}$ and according to Notation~\ref{n:simpl}
\begin{equation}\label{igual_trenz}
\nu \tau_{\star, \mathfrak{O}} =\tau \tau_{\diamond,\mathfrak{O}}^{-1} 
\beta(\mathfrak{O}) \tau_{\star,\mathfrak{O}}=
\tau \beta.
\end{equation}
Therefore by~\eqref{e:nuevos_m_n} and Remark~\ref{o:simpl} we have
\begin{equation}\label{primero_base}
m_{\star}(\mu^{\nu \tau_{\star, \mathfrak{O}}}_{1} )=
m_{\star}( \mu^{\nu \tau_{\star, \mathfrak{O}}}_{2})= 
(\bar{a},\bar{b}),
\end{equation}
and hence
\begin{equation*}%\label{segundo_base} 
m_{\star}\left(\Psi_{\tau_{\star, \mathfrak{O}}}((\eta^{\nu })^{-1} \cdot \mu_{1}^{\nu }  \cdot \eta^{\nu})\right)=
m_{\star}\left(\Psi_{\tau_{\star, \mathfrak{O}}}((\eta^{\nu })^{-1} \cdot \mu_{2}^{\nu }  \cdot \eta^{\nu})\right)
\end{equation*}
is a transposition that we will denote by $(a_i,b_i)$.

In particular this implies that regardless whether~\eqref{prim} or~\eqref{seg} 
holds above, in fact
$m_{\star}(\mu_{1}^{\beta_{i}\tau_{\star, \mathfrak{O}}})=
m_{\star}(\mu_{2}^{\beta_{i}\tau_{\star, \mathfrak{O}}})=(a_i,b_i).$
In particular we  have proved Theorem~\ref{t:met_discrim}\eqref{met_dis.1.a} (we will assume $a_i<b_i$). 

As for the path $\delta:=\mu_2 \cdot \mu_1$, case~\eqref{prim} leads to
\begin{equation}\label{delta_prim}
\delta^{\beta_{i}\tau_{\star,\mathfrak{O}}} = 
\Psi_{\tau_{\star, \mathfrak{O}}}(\delta^{ \zeta \nu})
=\Psi_{\tau_{\star, \mathfrak{O}}}(\mu_{2}^{ \zeta \nu } \cdot \mu_{1}^{ \zeta \nu }) 
=\Psi_{\tau_{\star, \mathfrak{O}}}((\eta^{\nu })^{-1} \cdot \mu_{2}^{\nu } \cdot \mu_{1}^{\nu } \cdot \eta^{\nu }),\\
\end{equation}
and analogously case~\eqref{seg} leads to
\begin{equation}\label{delta_seg}
\delta^{\beta_{i}\tau_{\star,\mathfrak{O}}} =
\Psi_{\tau_{\star, \mathfrak{O}}} ((\eta^{\nu })^{-1} \cdot \mu_{1}^{\nu } 
\cdot \mu_{2}^{\nu } \cdot \eta^{\nu }).
\end{equation}
Now we lift the path 
$\delta^{\beta_{i}\tau_{\star, \mathfrak{O}}} \subset 
\{t=\star\} \times \C$
at the point $(*,y_{l}) \in V_{\star}(f),l=1, \dots, N$. 
By Proposition~\ref{propcamhomolog}\eqref{propcamhomolog.a} 
its class in $C_1(\mathfrak{K}) /  B_1(\mathfrak{K})$ is:
\begin{itemize}
\item Case~\eqref{prim} (see~\eqref{delta_prim})
\begin{equation}\label{prim_elv}
(\delta^{\beta_{i}\tau_{\star, \mathfrak{O}}} )_{l}=
((\eta^{\nu \tau_{\star, \mathfrak{O}}})^{-1})_{l_1} +
(\mu_{2}^{\nu \tau_{\star, \mathfrak{O}}})_{l_2}  +
(\mu_{1}^{\nu \tau_{\star, \mathfrak{O}}})_{l_3} + (\eta^{\nu \tau_{\star, \mathfrak{O}}})_{l_4}
\end{equation}
where $l_1:=l$, $ l_2:=l_1^{m_{\star}((\eta^{\nu \tau_{\star, \mathfrak{O}}})^{-1})},
 l_3:=l_2^{m_{\star}(\mu_{2}^{\nu \tau_{\star, \mathfrak{O}}})}$ and 
$l_4:=l_3^{ m_{\star}(\mu_{1}^{\nu \tau_{\star, \mathfrak{O}}})}$. 
\item Case~\eqref{seg} (see~\eqref{delta_seg})
\begin{equation}\label{seg_elv}
(\delta^{\beta_{i}\tau_{\star, \mathfrak{O}}} )_{l}=
((\eta^{\nu \tau_{\star, \mathfrak{O}}})^{-1})_{l_1} +
(\mu_{1}^{\nu \tau_{\star, \mathfrak{O}}})_{l_2}  +
(\mu_{2}^{\nu \tau_{\star, \mathfrak{O}}} )_{l_3}+
(\eta^{\nu \tau_{\star, \mathfrak{O}}})_{l_4}
\end{equation}
where $l_1:=l$, $ l_2:=l_1^{m_{\star}((\eta^{\nu \tau_{\star, \mathfrak{O}}})^{-1}) },
 l_3:=l_2^{  m_{\star}(\mu_{1}^{\nu \tau_{\star, \mathfrak{O}}})}$ and 
$l_4:=l_3^{ m_{\star}(\mu_{2}^{\nu \tau_{\star, \mathfrak{O}}}) }$
\end{itemize}
In particular, for both cases we have 
\begin{equation}\label{eles}
\begin{array}{rcl}
l=a_i & \Rightarrow  & l_1=a_i, l_2=\bar{a}, l_3=\bar{b} \ \mbox{and} \ l_4=\bar{a}, \\
l=b_i & \Rightarrow &  l_1=b_i, l_2=\bar{b}, l_3=\bar{a} \ \mbox{and} \ l_4=\bar{b}, \\
l \neq a_i, b_i & \Rightarrow &  l_1 \neq a_i, b_i \ \mbox{and} \  
(\alpha:=)l_2=l_3=l_4 \neq \bar{a}, \bar{b}.
\end{array}
\end{equation}
Moreover, since the composition of the 
transpositions~$m_{\star}(\mu_{1}^{\nu \tau_{\star, \mathfrak{O}}})$ and $m_{\star}(\mu_{2}^{\nu \tau_{\star, \mathfrak{O}}})$ is the identity we have that
$l_4=l_1^{ m_{\star}((\eta^{\nu \tau_{\star, \mathfrak{O}}})^{-1})}$.
Hence $l^{ m_{\star}(\delta^{\beta_{i}\tau_{\star, \mathfrak{O}}})}=l$
and thus, by Proposition~\ref{propcamhomolog}\eqref{propcamhomolog.b} one concludes
\begin{equation}\label{opuestos}
((\eta^{\nu \tau_{\star, \mathfrak{O}}})^{-1})_{l_1}=-(\eta^{\nu \tau_{\star, \mathfrak{O}}})_{l_4}.
\end{equation} 

Finally by Proposition~\ref{propcamhomolog}\eqref{propcamhomolog.c} the lifting of
$(\delta^{\beta_{i}\tau_{\star, \mathfrak{O}}} )_{l}$ is a cycle in 
$H_1(V_{\star}(f))$. Case~\eqref{prim} implies
%\begin{footnotesize}
\begin{equation*}
(\delta^{\beta_{i}\tau_{\star, \mathfrak{O}}} )_{l} 
\stackrel{
\stackrel{
\mbox{\tiny \eqref{prim_elv},}}{\mbox{\tiny\eqref{opuestos}}}
}{=}  
(\mu_{2}^{\nu \tau_{\star, \mathfrak{O}}} )_{l_2} +  
(\mu_{1}^{\nu \tau_{\star, \mathfrak{O}}})_{l_3}\stackrel{\mbox{\tiny\eqref{eles}}}{=}
\end{equation*}
\begin{equation*}
=
\begin{cases}
(\mu_{2}^{\nu \tau_{\star,\mathfrak{O}}})_{\alpha} + 
(\mu_{1}^{\nu \tau_{\star,\mathfrak{O}}})_{\alpha}
\stackrel{
\stackrel{
\mbox{\tiny\eqref{primero_base},}}{\mbox{\tiny\eqref{rel_grup_homolog}}}
}{=} 0, \  l \neq a_i, b_i \\
(\mu_{2}^{\nu \tau_{\star, \mathfrak{O}}})_{\bar{a}} + 
(\mu_{1}^{\nu \tau_{\star,\mathfrak{O}}})_{\bar{b}} 
\stackrel{
\stackrel{\mbox{
\tiny\eqref{igual_trenz}, \eqref{igualdad_mu_j},}}{
\mbox{\tiny\eqref{igualdad_Delta_j}}}
}
{=}\hspace*{-.5cm}(\delta^{
\beta\tau_{\star,\psi}})_{\bar{a}}, \   l=a_i\\
(\mu_{2}^{\nu \tau_{\star,\mathfrak{O}}})_{\bar{b}} + 
(\mu_{1}^{\nu \tau_{\star,\mathfrak{O}}})_{\bar{a}}
\stackrel{
\stackrel{\mbox{\tiny
\eqref{igual_trenz}, \eqref{igualdad_mu_j},}}{
\mbox{\tiny\eqref{igualdad_Delta_j}}}
}{=} 
\hspace*{-.6cm} 
-(\delta^{
\beta\tau_{\star,\psi}})_{\bar{a}}, \   l=b_i.
\end{cases}
\end{equation*}
Analogously, we can use~\eqref{seg_elv} to deal with case~\eqref{seg}, one has
\begin{equation*}
(\delta^{\beta_{i}\tau_{\star, \mathfrak{O}}} )_{l}=
\begin{cases}
 0, \  l \neq a_i, b_i \\
(\mu_{1}^{\nu \tau_{\star,\mathfrak{O}}})_{\bar{a}} + 
(\mu_{2}^{\nu \tau_{\star,\mathfrak{O}}})_{\bar{b}}=-(\delta^{
\beta\tau_{\star,\psi}})_{\bar{a}} , \   l=a_i\\\
(\mu_{1}^{\nu \tau_{\star,\mathfrak{O}}})_{\bar{b}} + 
(\mu_{2}^{\nu \tau_{\star,\mathfrak{O}}})_{\bar{a}}= (\delta^{
\beta\tau_{\star,\psi}})_{\bar{a}}, \   l=b_i. \\
\end{cases}
\end{equation*}
Finally, just using the notation 
$\Delta_i:=(\delta^{\beta_{i}\tau_{\star, \mathfrak{O}}} )_{a_i}$ gives us Theorem~\ref{t:met_discrim}\eqref{met_dis_Delta}.

\subsection{The Picard-Lefschetz operators.}\label{s:monodro}
Our aim  is to prove Theorem~\ref{t:met_discrim}\eqref{met_dis_h(Delta)} (see Lemma~\ref{l:simp_1}). In order to do so 
we need the following Proposition~\ref{diag_con_en_homo_def}. 

\begin{proposicion}\label{diag_con_en_homo_def}
Let $\gamma_i,1\leq i \leq \mu(f),$ be as in Theorem~\ref{t:met_discrim}. 
Let $(*,y_l)~\in~V_{\star}(f)$, $l=1,\ldots,N$, be the preimage of $* \in \partial X$ by the
projection~\eqref{proy_prim}. Then for any $l$ there exists a monodromy $h_{\gamma_i}$ of $\gamma_i$ satisfying $h_{\gamma_i}(*,y_l)=(*,y_l)$ and
inducing the following conmutative diagram
\begin{equation*}
\begin{CD}
\pi_1(\mathcal{R}_1;(*,y_{l})) @>{(h_{{\gamma}_i} \vert)_*}>> \pi_1(\mathcal{R}_1;(*,y_{l})) \\
@V{(\tilde{\pi} \circ \phi_{f  })_*}VV @VV{(\tilde{\pi} \circ \phi_{f  })_*}V \\
 \pi_1(X \setminus {\bf x}^{\star} ;* ) 
@>>{\Psi_{(\bar{\pi} \vert_{\mathfrak{D}_{f}})_*({\gamma}_i)}}> 
\pi_1(X \setminus {\bf x}^{\star} ;* ) \\
\end{CD}
\begin{picture}(0,0)
\put(-100,3){$\circlearrowleft$}
%\put(-100,25){$\circlearrowleft$}
\end{picture}
\end{equation*}
where $\mathcal{R}_1$ is defined as in~\eqref{X_F} and $\tilde{\pi}$ is 
the projection $(t,x)\mapsto x$.
\end{proposicion}

We will first see that Theorem~\rm{\ref{t:met_discrim}\eqref{met_dis_h(Delta)}} is a consequence 
of Proposition~\ref{diag_con_en_homo_def}.

\begin{proof}[Proof of Theorem~\rm{\ref{t:met_discrim}\eqref{met_dis_h(Delta)}}]
Using the commutative diagram of Proposition~\ref{diag_con_en_homo_def} (for $l=a_s$),
the one from~\eqref{diagramaprimario} (for ${\bf s}^1={\bf s}^2=~{\bf x}^{\star}$),
and the fact that one has the equality
$\mathbb{I}_{{\bf x}^\star, {\bf x}^\star,\mathfrak{O}} 
\left(  (\bar{\pi} \vert_{\mathfrak{D}_{f}})_*(\gamma_i) \right)= \nabla_{\star,\mathfrak{O}}(\gamma_i) $ 
one obtains the following diagram
\begin{equation}\label{superdiagram}
\begin{CD}
\pi_1(\mathcal{R}_1;(*,y_{a_s})) @>{(h_{{\gamma}_i} \vert)_*}>> \pi_1(\mathcal{R}_1;(*,y_{a_s})) \\
@V{(\tilde{\pi} \circ \phi_{f  })_*}VV @VV{(\tilde{\pi} \circ \phi_{f  })_*}V \\
 \pi_1(X \setminus {\bf x}^{\star} ;* ) 
@>{\Psi_{(\bar{\pi} \vert_{\mathfrak{D}_{f}})_*({\gamma}_i)}}>> 
\pi_1(X \setminus {\bf x}^{\star} ;* ) \\
@V{\Psi^{-1}_{\tau_{\star,\mathfrak{O}}}}VV @VV{\Psi^{-1}_{\tau_{\star,\mathfrak{O}}}}V \\
\pi_1(X \setminus {\bf x} ;* )  @>{\Psi_{\nabla_{\star,\mathfrak{O}}({\gamma}_i)}}>> \pi_1(X \setminus {\bf x} ;* ).
\end{CD}
\begin{picture}(0,0)
\put(-103,-25){$\circlearrowleft$}
\put(-103,25){$\circlearrowleft$}
\end{picture}
\end{equation}
Let $\beta_{s,i}:=
\beta_s  \nabla_{\star, \mathfrak{O} }(\gamma_i) \in \B_k$ as in Theorem~\ref{t:met_discrim}\eqref{met_dis_h(Delta)}.
Then by Theorem~\ref{t:met_discrim}\eqref{met_dis_Delta} and~\eqref{superdiagram} we have
\begin{equation*}
(h_{\gamma_i})_*(\Delta_s)=(h_{\gamma_i})_*((\delta^{\beta_s  \tau_{\star,\mathfrak{O}}})_{a_s})=
(\delta^{\beta_{s,i} \tau_{\star,\mathfrak{O}}})_{a_s}.
\end{equation*}
This proves the first part of Theorem~\ref{t:met_discrim}\eqref{met_dis_h(Delta)}.

For the second part, let $\mathcal{Y}_{\gamma_i}$ be the braid with $N$ 
strands defined by the roots of $f(*,y)=\gamma_i(u),u \in  [0,1].$
This is a pure braid by~Proposition~\ref{diag_con_en_homo_def}, therefore
$$
m_{\star}(\Psi_{(\bar{\pi}\vert_{\mathfrak{D}_f})_*(\gamma_i)}
(\mu_{n}^{\beta_s \tau_{\star,\mathfrak{O}  }}))
=(a_s^{\mathcal{Y}_{\gamma_i}}, b_s^{\mathcal{Y}_{\gamma_i}})=(a_s,b_s),n\in \{1,2\}.
$$
Also note that
$
\Psi_{(\bar{\pi}\vert_{\mathfrak{D}_f})_*(\gamma_i)}(\mu_{n}^{\beta_s \tau_{\star,\mathfrak{O}}  })=
\mu_{n}^{\beta_{s,i}  \tau_{\star,\mathfrak{O}}}
$
according to~\eqref{superdiagram} and hence
\begin{equation}\label{igual_n_m}
m_{\star}(\mu_{n}^{\beta_{s,i}  \tau_{\star,\mathfrak{O}}})
=(a_s, b_s),n\in\{1,2\}.
\end{equation} 

By Artin's Theorem, the list
$(\mu_{1}^{\beta_{s,i}  \tau_{\star,\mathfrak{O}}},
\ldots, \mu_{k}^{\beta_{s,i}  \tau_{\star,\mathfrak{O}}})$ 
is a geometric basis of the group $\pi_1(X \setminus {\bf x}^{\star}; *)$. 
Therefore, by~\eqref{igual_n_m}, relations~\eqref{rel_grup_homolog} and Proposition~\ref{propcamhomolog}\eqref{propcamhomolog.a} we have the 
following equalities in homology
\begin{equation*}
\begin{aligned}
(\delta^{\beta_{s,i}  \tau_{\star,\mathfrak{O}}})_{\alpha}&=
(\mu_{2}^{\beta_{s,i}  \tau_{\star,\mathfrak{O}}})_{\alpha} +
(\mu_{1}^{\beta_{s,i}  \tau_{\star,\mathfrak{O}}})_{\alpha}=0  \ \mbox{for} \ \alpha \neq a_{s}, b_{s}, \\
(\delta^{\beta_{s,i}  \tau_{\star,\mathfrak{O}}})_{a_s}&=
(\mu_{2}^{\beta_{s,i}  \tau_{\star,\mathfrak{O}}})_{a_s} +
(\mu_{1}^{\beta_{s,i}  \tau_{\star,\mathfrak{O}}})_{b_s}  \ \mbox{and}\\
(\delta^{\beta_{s,i}  \tau_{\star,\mathfrak{O}}})_{b_s}&=
(\mu_{2}^{\beta_{s,i}  \tau_{\star,\mathfrak{O}}})_{b_s} +
(\mu_{1}^{\beta_{s,i}  \tau_{\star,\mathfrak{O}}})_{a_s} = 
-(\mu_{2}^{\beta_{s,i}  \tau_{\star,\mathfrak{O}}})_{a_s} -
(\mu_{1}^{\beta_{s,i}  \tau_{\star,\mathfrak{O}}})_{b_s}.
\end{aligned}
\end{equation*}
By Picard-Lefschetz Theorem 
$(h_{\gamma_i})_{*}(\Delta_s)=\Delta_s - n \, \Delta_i, \, n \in \Z.$
Then, since the unique non-zero liftings of 
$\delta^{\beta_{s,i}\tau_{\star,\mathfrak{O}}}$
are opposite and $\left\langle \Delta_s, \Delta_i\right\rangle $ has rank 2,
both non-zero liftings can not be written as $\Delta_s - n \, \Delta_i$.
This proves the second part of Theorem~\ref{t:met_discrim}\eqref{met_dis_h(Delta)}.
\end{proof}

We will devote the remaining part of the paper to proving
Proposition~\ref{diag_con_en_homo_def}.

\vspace*{14pt}
Since the monodromy only depends on the homotopy class of the defining path
we can assume that $\gamma_i$ is the path $\gamma$ considered in 
Section~\ref{s:meridiano}. By~\eqref{e:ele_asoc} the operator $h_{\gamma}$ can 
be factorized as $h_{r} \circ h_{S} \circ h_{r^{-1}}$ where $h_{r}$ and $h_{r^{-1}}$ 
were constructed in Section~\ref{s:cons.difeo}. Let us then study the monodromy 
$h_{S}$ more closely.

\subsubsection{The monodromy $h_{S}$.}
We will first construct a family of maps 
$$
\Gamma_u:V_{\diamond}(f) \longrightarrow V_{S(u)}(f), \, \, u \in [0,1],
$$
such that $h_{S}:= \Gamma_1$ is the identity on $\partial V_{\diamond}(f) $.
To do this we will decompose 
$$
\mathcal{S}_u:=V_{S(u)}(f) \setminus  \phi_{f,x}  ^{-1}(\{t=S(u)\}\times {\bf x}^{S(u)})
$$
into two sets~$\mathcal{S}_{u}^{(1)}$ and $\mathcal{S}_{u}^{(2)}$. Then we will 
describe the restrictions $\Gamma_{u}^{(1)}$ and $\Gamma_{u}^{(2)}$ of $\Gamma_u$ 
to $\mathcal{S}_{u}^{(1)}$ and $\mathcal{S}_{u}^{(2)}$ and finally we will reconstruct
$\Gamma_u$ by adding the set of
points~$\phi_{f,x}^{-1}(\{t=\diamond\}\times {\bf x}^{\diamond})$ continuously.

\paragraph{\textit{Construction of the sets
$\mathcal{S}_{u}^{(1)}$ and $\mathcal{S}_{u}^{(2)}$.}}
Let $\psi$ be an appropriate angle as considered in Section~\ref{s:meridiano}.
Let us consider $H$ as in~\eqref{e:H} where $[\diamond,t_i)$ is replaced 
by the positive parametrization of $S$ based at~$\diamond$ (taking $H_{\tau_{S(u),\psi}}$ as special ispotopy).
This map is defined on $[0,1]\times(\C\times[0,1])$ and we could also 
assume that it satisfies
$$H(0,x,v)=H(1,x,v), x\in\C\setminus\D_1\left( \frac{2j+1}{2} 
\right), v \in [0,1].$$
Then $\mathfrak{T}:= \bigcup_{u \in [0,1]} \mathfrak{T}_{u}$
is a topological torus whose fibers
$$ 
\mathfrak{T}_{u}:=  H\left(u,\D_1\left( \dfrac{2j+1}{2} \right) ,1 \right)\subset \Int(X), \ u \in [0,1], 
$$
are such that
$
x_{j}^{S(u)},  x_{j+1}^{S(u)} \in \mathfrak{T}_{u}$ and
$ \mathfrak{T}_{u}  \cap ( 
{\bf x}^{S(u)} \setminus \{x_{j}^{S(u)},  x_{j+1}^{S(u)} \})=\emptyset,
$
where the elements $x_{j}^{S(u)},  x_{j+1}^{S(u)}~\in~{\bf x}^{S(u)}.$

We define  
$
\mathcal{S}_{u}^{(1)}:= \phi_{f,x}^{-1}(\{t=S(u)\} \times 
( \mathfrak{T}_{u} \setminus \{x_{j}^{S(u)},  x_{j+1}^{S(u)}) \}) \subset \mathcal{S}_u $ and its complement 
$\mathcal{S}_{u}^{(2)}:=\mathcal{S}_u \setminus 
\Int(\mathcal{S}_{u}^{(1)})$, $u \in [0,1]$. By Lemma~\ref{l:proy_fa}, 
the set $\mathcal{S}_{u}^{(1)}$ has $N-1$ components. 
$N-2$ of which, say $A_{u}^{1}, \dots, A_{u}^{N-2}$, are diffeomorphic to 
$\mathfrak{T}_{u} \setminus \{x_{j}^{S(u)},  x_{j+1}^{S(u)}\}$ and one, say 
$A_{u}^{N-1}$ is diffeomorphic to 
$(\mathbb{S}^{1} \times [-1,1]) \setminus \{(-1,0),(1,0)\}$.

\paragraph{\textit{Construction of the maps
$\Gamma_{u}^{(1)}$ and~$\Gamma_{u}^{(2)}$.}}

Let us define $\varrho:\C \times [0,1]  \rightarrow \C$ as
\begin{equation*}
\varrho(x,u):=
\begin{cases}
\frac{2j+1}{2}+ re^{\sqrt{-1}( \omega +2 \pi u (1-r) )} & \mbox{if} \, \, x \, 
=\frac{2j+1}{2}+ re^{\sqrt{-1} \omega}, 0 \leq r \leq 1, \\
x & \mbox{if} \, \, x \notin  \D_1(\frac{2j+1}{2})
\end{cases}
\end{equation*}
If we denote $\varrho_u:=\varrho(\bullet,u)$, then
$\varrho$ is an isotopy from the identity $\varrho_0$ to the 
\emph{Dehn twist} $\varrho_1$ around the points $j$ and $j+1$. 
If we define the diffeomorphism $\tilde{\varrho}_u$, $u \in [0,1]$, as a composition
\begin{equation*}
\begin{CD}
\C @>\tilde{\varrho}_u>> \C \\
@VH_{\tau}(\bullet,1)VV @AAH^{-1}_{\tau_{S(u),\psi}}(\bullet,1)A \\
\C @>\varrho_u>> \C
\end{CD}
\end{equation*}
then  $\tilde{\varrho}_1$ is the identity out of~$\intt(\mathfrak{T}_0)$ 
and the diffeomorphisms $\tilde{\varrho}_u,$  $u \in [0,1]$,  define an isotopy 
$\tilde{\varrho}$ from the identity $\tilde{\varrho}_0$
to $\tilde{\varrho}_1$.   
\begin{lema}\label{rho_es_theta}
The set of $k$ strands defined by $\tilde{\varrho}({\bf x}^{\diamond},u)$, 
$u \in [0,1],$ is isotopic to the braid $\theta$ described in Notation~\ref{n:mon}.
\end{lema}
\begin{proof}
By construction, the image of $\tilde{\varrho}({\bf x}^{\diamond},u)$, $u \in [0,1]$, 
by $\mathbb{I}_{{\bf x}^{\diamond},{\bf x}^{\diamond}, \mathfrak{O}(\psi)}$
is $\sigma_{j}$. Also, by Proposition~\ref{deltaij}, the image of $\theta$ is
again~$\sigma_{j}$. Therefore both braids are equal in 
$\B({\bf x}^{\diamond},{\bf x}^{\diamond})$ since 
by Proposition~\ref{ident_trenzas}\eqref{ident_trenzas.1} 
$\mathbb{I}_{{\bf x}^{\diamond},{\bf x}^{\diamond}, \mathfrak{O}(\psi)}$ 
defines an isomorphism between $\B({\bf x}^{\diamond},{\bf x}^{\diamond})$ and~$\B_{k}$. 
\end{proof}

Consider the following diagram for $l=1, \dots, N-1, u \in [0,1],$
\begin{equation}\label{diag_cub_base}
\xymatrix{
 & & A_{u}^{l} \ar[d]^{\pi \vert} \\
A_{0}^{l}  \ar[r]_-{\pi \vert} \ar@{.>}[urr]&  \mathfrak{T}_{0} \setminus 
\{x_{j}^{\diamond},  x_{j+1}^{\diamond}\} \ar[r]_{\tilde{\varrho}_{u}} & 
\mathfrak{T}_{u} \setminus \{x_{j}^{S(u)},  x_{j+1}^{S(u)}\}.
}
\end{equation}
 We want to lift 
$\tilde{\varrho}_u \circ \pi \vert_{ A_{0}^{l}}$.
Let us consider $p_0 \in \partial \mathfrak{T}_{0} $ and 
$p_u:=\tilde{\varrho}_u(p_0) \in \partial \mathfrak{T}_{u}$. Also, we denote 
\begin{equation}\label{e:pq}
(p_u,q_{u}^{l})  \in \partial  A_{u}^{l}, 
l= 1, \dots, N-2, (p_u,q_{u}^{N-1}), (p_u,q_{u}^{N}) 
\in  \partial A_{u}^{N-1}.
\end{equation}
the preimage of $p_u$ by $\pi \vert_{\mathcal{S}_u }$.
Therefore, if $l=1, \cdots, N-2,$ the covering map $\pi \vert_{ A_{u}^{l}}$ is a 
diffeomorphism and thus its lifting is trivial. Hence the image of 
$(p_0^{},q_{0}^{l})$ is~$(p_u, q_{u}^{l})$, $u \in [0,1]$. 

\begin{observacion}\label{girodehn}
Since $\tilde{\varrho}$ is an isotopy from the identity to $\tilde{\varrho}_{1}$,  
lifting $\tilde{\varrho}_u \circ \pi \vert_{ A_{0}^{l}}$ produces an isotopy from
the identity of $A_{0}^{l}$, $l=1, \ldots, N-2$ to 
$(\pi \vert_{ A_{0}^{l}})^{-1} \circ \tilde{\rho}_1 \circ \pi \vert_{ A_{0}^{l}}$.
\end{observacion}

If $l=N-1$ let $\mu, \beta$ be generators of the free group 
$F_2 \simeq\pi_1(\mathfrak{T}_{u} \setminus \{x_{j}^{S(u)},
x_{j+1}^{S(u)}\}; p_{u}).$
The component $A_{u}^{N-1}$ is diffeomorphic to 
$(\mathbb{S}^{1} \times [-1,1]) \setminus~\{(-1,0),(1,0)\}$ and the covering map $\pi \vert_{A_{u}^{N-1}}$ is 2:1, then
$ 
(\pi \vert_{A_{u}^{N-1}})_{*} \left( \pi_1(A_{u}^{N-1};(p_u,q_{u}^{N-1}))\right) 
$
is generated by $\mu^2, \beta^2, \beta \, \cdot \, \mu$, i.e.,
words from the alphabet $\mu, \beta$ of even length.

Let $0 \leq u <1$. Since the restriction of $\tilde{\varrho}_u$ on 
$\mathfrak{T}_{0} \setminus \{x_{j}^{\diamond},  x_{j+1}^{\diamond}\}$
induces the identity on $F_2$, then 
$(\tilde{\varrho}_u \circ \pi \vert_{A_{0}^{N-1}})_{*}
\left(\pi_1(A_{0}^{N-1};(p_0,q_{0}^{N-1}))\right)$
is again the group of words in the alphabet $\mu, \beta$ of even length.
Let $u=1$. Since the restriction of $\tilde{\varrho}_1$ on 
$\mathfrak{T}_{0} \setminus \{x_{j}^{\diamond},  x_{j+1}^{\diamond}\}$
induces in homotopy the isomorphism given by
\begin{equation*}
F_2 \longrightarrow  F_2, \mu  \mapsto  \beta, \beta  \mapsto \beta  \cdot  \mu  \cdot  \beta^{-1},
\end{equation*}
then 
$(\tilde{\varrho}_1 \circ \pi \vert_{A_{0}^{N-1}})_{*}
(\pi_1(A_{0}^{N-1};(p_0,\tilde{p}_{0}^{N-1})))$
is again the group of words in the alphabet $\mu, \beta$ of even length.
Then by the Lifting Theorem there exists a unique lifting of the 
diagram~\eqref{diag_cub_base} such that the image of 
$(p_0^{},q_0^{N-1})$ is $(p_u^{},q_u^{N-1})$.

Collecting the liftings of diagram~\eqref{diag_cub_base} we can construct 
$\Gamma_{u}^{(1)}: \mathcal{S}_{0}^{(1)} \rightarrow  \mathcal{S}_{u}^{(1)},$ $u\in [0,1],$
a continuous family of maps such that the image of $(p_{0},q_{0}^{l})$ 
is $ (p_u,q_{u}^{l})$, $l=1, \dots, N-1$. Note that they also verify the 
commutative diagram
\begin{equation}\label{diag_mon_(1)}
\begin{CD}
\mathcal{S}_0^{(1)} @>\Gamma_{u}^{(1)}>> \mathcal{S}_u^{(1)} \\
@V{\pi \vert}VV @VV{\pi \vert}V \\
\mathfrak{T}_{0} \setminus \{x_{j}^{\diamond},  x_{j+1}^{\diamond}\} @>\tilde{\varrho}_{u}>> 
\mathfrak{T}_{u} \setminus \{x_{j}^{S(u)}, x_{j+1}^{S(u)}\}.
\end{CD}
\begin{picture}(0,0)
%\graphpaper(-90,-50)(100,75)
\put(-108,-2){$\circlearrowleft$}
\end{picture}
\end{equation}
Therefore, the image of $(p_0,q_{0}^{N})$ is $(p_u, q_{u}^{N})$.

To construct $\Gamma_{u}^{(2)}$ let us consider
$V_{t'}(f)^*:=V_{t'}(f) \setminus  \phi_{f,x}^{-1}(\{t=t'\}\times {\bf x}^{t'})$, 
$t' \in D,$ (see~\eqref{e:ele_asoc}). $H$ can be extended to $D$, then
\begin{small}
\begin{equation*}
\begin{CD}
\bigcup_{t' \in D} V_{t'}(f)^* \setminus \Int\left( \phi_{f,x}^{-1}\left( 
\{t=t'\} \times H\left( t',\D_1\left( \frac{2j+1}{2}\right) \setminus \{j,j+1\} ,1\right) 
 \right) 
 \right) \\
 @VV{f \vert}V \\
  D
\end{CD}
\end{equation*}
\end{small}
is a trivial fibration with fiber $\mathcal{S}_{0}^{(2)}$. Then a continuous 
family of diffeomorphisms 
$\Gamma_{u}^{(2)}: \mathcal{S}_{0}^{(2)} \rightarrow \mathcal{S}_{u}^{(2)}$ 
can be constructed, where $\Gamma_{1}^{(2)}$ is the identity and
such that the points~\eqref{e:pq} verify that the image of $(p_0^{},q_{0}^{l})$ is 
$(p_u^{},q_{u}^{l})$, $l=1, \ldots, N$.

\paragraph{\textit{Construction of the continuous family $\Gamma_{u}$.}}
We will first show that the maps $\Gamma_{u}^{(1)}$ and~$\Gamma_{u}^{(2)}$
agree on~$\mathcal{S}_{0}^{(1)} \cap \mathcal{S}_{0}^{(2)}$. Note that the
intersection $\mathcal{S}_{u}^{(1)} \cap \mathcal{S}_{u}^{(2)}$, $u \in [0,1]$, 
consists of $N$ components $B_u^{l}$, where
$B_{u}^{l}:= \partial A_{u}^{l}, l=1, \dots, N-2,$
and 
$B_{u}^{N-1} \cup B_{u}^{N}:=\partial A_{u}^{N-1}$, 
all of which are homeomorphic to $\partial \mathfrak{T}_{u}$.
Hence the following diagram
\begin{equation}\label{diag_cub_base_borde}
\xymatrix{
& & B_{u}^{l} \ar[d]^{\pi \vert}\\
B_{0}^{l} \ar[r]_-{\pi \vert} \ar@{.>}[urr]& \partial \mathfrak{T}_{0} \ar[r]_{\tilde{\varrho}_{u}} & \partial \mathfrak{T}_{u},
}
\end{equation}
$l=1, \dots, N, u \in [0,1],$ can be lifted to a unique map 
such that the points defined in~\eqref{e:pq} verify that the image of 
$(p_0^{},q_{0}^{l})$ is $(p_u^{},q_{u}^{l})$, $l=1, \ldots, N$.
Also note that the restriction of the maps $\Gamma_u^{(1)}$ and $\Gamma_u^{(2)}$ 
to $B_{0}^{l}$ coincides, by construction, with the lifting of
diagram~\eqref{diag_cub_base_borde} where the image of $(p_0^{},q_{0}^{l})$ 
is the point $(p_u^{},q_{u}^{l})$, therefore both maps are equal by uniqueness. 
Finally we have the continuous family of maps 
$\Gamma_u:\mathcal{S}_0 \rightarrow \mathcal{S}_u$, $u \in [0,1]$,
which can be extended to $V_{\diamond}(f)$ by adding points of
$\phi_{f,x}^{-1}(\{t=\diamond\}\times {\bf x}^{\diamond})$ continuously. 

Finally we define the monodromy $h_{S}$ as $\Gamma_{1}$. Since $\Gamma_1^{(2)}$ is the 
identity and $\Gamma_1^{(1)}$ verifies diagram~\eqref{diag_mon_(1)} then we have 
the following commutative diagram
\begin{equation}\label{diag_mon_def}
\begin{CD}
\mathcal{S}_0 @>h_{S} \vert>> \mathcal{S}_0 \\
@V{\phi_{f,x}}VV @VV{\phi_{f,x}}V \\
\mathcal{X}_0 @>{1 \times \tilde{\varrho}_{1}}>> \mathcal{X}_0,
\end{CD}
\begin{picture}(0,0)
%\graphpaper(-90,-50)(100,75)
\put(-39,-2){$\circlearrowleft$}
\end{picture}
\end{equation}
where $\mathcal{X}_0$ is defined as in~\eqref{X_F}.

\begin{observacion}\label{o:mon_homo}
By Remark~\ref{girodehn} and the fact that  $\Gamma_1^{(2)}$ is the identity, 
there exists an isotopy from $h_{S}$ to a monodromy which coincides with $h_{S}$ 
in $\Int(A_{0}^{N-1})$ and is the identity in 
$V_{\diamond}(f) \setminus \Int(A_{0}^{N-1})$.
\end{observacion}

\subsubsection{The monodromy $h_{\gamma}$.}

By construction
$h_{\gamma}:=  h_{r^{-1}} \circ h_{S} \circ h_{r}$
is the identity outside $h_{r^{-1}}(\Int(\mathcal{S}_{0}^{(1)}))$.
Also note that, by Remark~\ref{o:mon_homo}, there exists an isotopy from 
$h_{\gamma}$ to a monodromy which is the identity outside 
$h_{r^{-1}}(\Int(A_{0}^{N-1}))$.

Let 
$(*,y_{l}) \in V_{\star}(f)$, $l=1, \dots, N,$ be the
preimage of $* \in \partial X$ by $\pi \vert_{V_{\star}(f)}$.
By construction we have
\begin{equation}\label{nml}
\begin{matrix}
h_{\gamma}(*, y_{l})=(*, y_{l}), l=1, \dots, N.
\end{matrix}
\end{equation}

Using the diagrams~\eqref{diagrama_sub},~\eqref{diagrama_subnv},~\eqref{diag_mon_def},
the definition of $\bar{H}_1$ in~\eqref{e:Hu}, and the fact that
$\mathcal{S}_0=\mathcal{R}_0$ we have the family of commutative diagrams
\begin{equation}\label{super_diag}
\begin{CD}
\mathcal{R}_1 @>{h_{r}}>> \mathcal{R}_0 @>{h_{S}}>> \mathcal{R}_0 @>{h_{r^{-1}}}>> \mathcal{R}_1 \\
@VV{\phi_{f,x}}V @VV{\phi_{f,x}}V @VV{\phi_{f,x}}V @VV{\phi_{f,x}}V  \\
\mathcal{X}_1 @>{\bar{H}_{1}^{-1}}>>  \mathcal{X}_0 @>{1 \times \tilde{\varrho}_{1}}>> \mathcal{X}_0 @>{\bar{H}_1}>> \mathcal{X}_1 \\
@VV{\tilde{\pi}}V @VV{\tilde{\pi}}V @VV{\tilde{\pi}}V @VV{\tilde{\pi}}V  \\
X \setminus {\bf x}^{\star} @>{H_{\beta}^{-1}}>> X \setminus {\bf x}^{\diamond} 
@>{\tilde{\varrho}_{1}}>> X \setminus {\bf x}^{\diamond} @>{H_{\beta}}>> X \setminus {\bf x}^{\star},
\end{CD}
\begin{picture}(0,0)
%\graphpaper(-90,-50)(100,75)
\put(-56,-27){$\circlearrowleft$}
\put(-56,20){$\circlearrowleft$}
\put(-126,-27){$\circlearrowleft$}
\put(-126,20){$\circlearrowleft$}
\put(-196,-27){$\circlearrowleft$}
\put(-196,20){$\circlearrowleft$}
\end{picture}
\end{equation}
where $\tilde{\pi}(t,x):=x$. 

Using Lemma~\ref{rho_es_theta} 
and \eqref{n:pi_bar} one can check that 
the composition of the bottom row of maps in~\eqref{super_diag} 
induces
$
 \Psi_{\beta} \circ \Psi_{\theta} \circ \Psi_{\beta^{-1}}= 
\Psi_{(\bar{\pi} \vert_{\mathfrak{D}_{f}})_*(\gamma)} 
$
in homotopy. Therefore, by~\eqref{nml}, we have the commutative diagram
\begin{equation*}%\label{diag_con_en_homo}
\begin{CD}
\pi_1(\mathcal{R}_1;(*,y_{l})) @>{(h_{\gamma} \vert)_*}>> \pi_1(\mathcal{R}_1;(*,y_{l})) \\
@V{(\tilde{\pi} \circ \phi_{f,x })_*}VV  @VV{(\tilde{\pi} \circ \phi_{f,x })_*}V \\
\pi_1( X \setminus {\bf x}^{\star} ;* ) @>{\Psi_{(\bar{\pi} \vert_{\mathfrak{D}_{f}})_*(\gamma)}}>> \pi_1(X \setminus {\bf x}^{\star} ;* ) \\
\end{CD}
\begin{picture}(0,0)
%\graphpaper(-90,-50)(100,75)
\put(-98,0){$\circlearrowleft$}
%\put(-100,25){$\circlearrowleft$}
\end{picture}
\end{equation*}
where $l=1, \dots, N$, and hence Proposition~\ref{diag_con_en_homo_def} follows.

\def\cprime{$'$} \def\cprime{$'$} \def\cprime{$'$} \def\cprime{$'$}
  \def\cprime{$'$} \def\cprime{$'$} \def\cprime{$'$} \def\cprime{$'$}
  \def\cprime{$'$} \def\cprime{$'$}
\providecommand{\bysame}{\leavevmode\hbox to3em{\hrulefill}\thinspace}
\providecommand{\MR}{\relax\ifhmode\unskip\space\fi MR }
% \MRhref is called by the amsart/book/proc definition of \MR.
\providecommand{\MRhref}[2]{%
  \href{http://www.ams.org/mathscinet-getitem?mr=#1}{#2}
}
\providecommand{\href}[2]{#2}

\end{document}